\documentclass[11pt]{amsart} 
\usepackage{amsfonts, amsmath, amssymb, amsthm, color}
\usepackage[margin=1in]{geometry}
\usepackage{nccmath}
\usepackage[dvipsnames]{xcolor}
\usepackage[colorlinks=true]{hyperref}

\newtheorem{thm}{Theorem}[section]
\newtheorem{lemma}[thm]{Lemma}
\newtheorem{prop}[thm]{Proposition}
\newtheorem{cor}[thm]{Corollary}

\theoremstyle{definition}
\newtheorem{rmk}[thm]{Remark}
\newtheorem{defn}[thm]{Definition}

\newcommand{\ga}{\gamma}
\newcommand{\ep}{\epsilon}
\newcommand{\vep}{\varepsilon}
\newcommand{\vph}{\varphi}

\newcommand{\vsi}{\varsigma}
\newcommand{\pa}{\partial}

\newcommand{\N}{\mathbb{N}}
\newcommand{\R}{\mathbb{R}}
\renewcommand{\S}{\mathbb{S}}

\newcommand{\mba}{\mathbf{a}}

\newcommand{\mbz}{\mathbf{z}}
\newcommand{\mbone}{\mathbf{1}}

\newcommand{\mx}{{\mu,\xi}}
\newcommand{\mxe}{{\mu_{\ell},\xi_{\ell}}}
\newcommand{\bmx}{{\bmu,\xi}}
\newcommand{\mbmu}{\boldsymbol{\mu}}
\newcommand{\mbxi}{\boldsymbol{\xi}}
\newcommand{\mbmx}{\mbmu,\mbxi}

\newcommand{\mca}{\mathcal{A}}
\newcommand{\mcb}{\mathcal{B}}
\newcommand{\mcc}{\mathcal{C}}
\newcommand{\mcd}{\mathcal{D}}
\newcommand{\mce}{\mathcal{E}}
\newcommand{\mcf}{\mathcal{F}}
\newcommand{\mch}{\mathcal{H}}
\newcommand{\mcj}{\mathcal{J}}
\newcommand{\mck}{\mathcal{K}}
\newcommand{\mcm}{\mathcal{M}}

\newcommand{\mcq}{\mathcal{Q}}
\newcommand{\mcr}{\mathcal{R}}
\newcommand{\mcs}{\mathcal{S}}
\newcommand{\mct}{\mathcal{T}}
\newcommand{\mcw}{\mathcal{W}}

\newcommand{\mfB}{\mathfrak{B}}
\newcommand{\mfh}{\mathfrak{h}}
\newcommand{\mfm}{\mathfrak{m}}

\newcommand{\bc}{\bar{c}}
\renewcommand{\bf}{\bar{f}}
\newcommand{\by}{\bar{y}}
\newcommand{\bmu}{\bar{\mu}}
\newcommand{\bvp}{\bar{\vph}}
\newcommand{\ovs}{\overline{S}}
\newcommand{\ovmcw}{\overline{\mcw}}
\newcommand{\ovmfh}{\overline{\mfh}}

\newcommand{\hb}{\hat{b}}
\newcommand{\hc}{\hat{c}}
\newcommand{\hps}{\hat{\psi}}

\newcommand{\tc}{\tilde{c}}
\newcommand{\td}{\tilde{d}}
\newcommand{\tf}{\tilde{f}}

\newcommand{\tu}{\tilde{u}}
\newcommand{\ty}{\tilde{y}}

\newcommand{\tmu}{\tilde{\mu}}
\newcommand{\tbmu}{\tilde{\bar{\mu}}}

\newcommand{\tvp}{\tilde{\vph}}
\newcommand{\tps}{\tilde{\psi}}
\newcommand{\txi}{\tilde{\xi}}
\newcommand{\tyn}{\ty_{\text{north}}}

\newcommand{\wtom}{\widetilde{\Omega}}
\newcommand{\wtmch}{\widetilde{\mch}}
\newcommand{\wtmcs}{\widetilde{\mcs}}
\newcommand{\wtsh}{\widetilde{\sharp}}

\newcommand{\ar}{{\alpha,\rho}}

\newcommand{\tin}{\textnormal{in}}
\newcommand{\out}{\textnormal{out}}
\newcommand{\tpa}{\textnormal{par}}

\renewcommand{\(}{\left(}
\renewcommand{\)}{\right)}
\newcommand{\la}{\left\langle}
\newcommand{\ra}{\right\rangle}

\allowdisplaybreaks
\numberwithin{equation}{section}

\makeatletter
\@namedef{subjclassname@2020}{%
\textup{2020} Mathematics Subject Classification}
\makeatother

\begin{document}
\title[Infinite-time blowing-up solutions to perturbations of the Yamabe flow]{Infinite-time blowing-up solutions \\ to small perturbations of the Yamabe flow}

\author{Seunghyeok Kim}
\address[Seunghyeok Kim]{Department of Mathematics and Research Institute for Natural Sciences, College of Natural Sciences, Hanyang University,
222 Wangsimni-ro Seongdong-gu, Seoul 04763, Republic of Korea}
\email{shkim0401@hanyang.ac.kr shkim0401@gmail.com}

\author{Monica Musso}
\address[Monica Musso]{Department of Mathematical Sciences, University of Bath, Bath BA2 7AY, United Kingdom}
\email{m.musso@bath.ac.uk}

\begin{abstract}
Under the validity of the positive mass theorem, the Yamabe flow on a smooth compact Riemannian manifold of dimension $N \ge 3$ is known to exist for all time $t$
and converges to a solution to the Yamabe problem as $t \to \infty$.
We prove that if a suitable perturbation, which may be smooth and arbitrarily small,
is imposed on the Yamabe flow on any given Riemannian manifold $M$ of dimension $N \ge 5$,
the resulting flow may blow up at multiple points on $M$ in the infinite time.
Our proof is constructive, and indeed we construct such a flow by using solutions of the Yamabe problem on the unit sphere $\S^N$ as blow-up profiles.
We also examine the stability of the blow-up phenomena under a negativity condition on the Ricci curvature at blow-up points.
\end{abstract}

\date{\today}
\maketitle

\section{Introduction}
\subsection{History}
Let $(M,g_0)$ be a smooth compact Riemannian manifold of dimension $N \ge 3$. The Yamabe flow on $(M,g_0)$ is defined as
\begin{equation}\label{eq:Yamabef}
\begin{cases}
\displaystyle {\pa \over \pa t} g(t) = - \(S\,[g(t)] - \ovs\,[g(t)]\) g(t) &\text{for } t \in (0,T),\\
g(0) = g_0
\end{cases}
\end{equation}
where $T$ is the maximal time, $S\,[g(t)]$ is the scalar curvature of the metric $g(t)$, and $\ovs\,[g(t)]$ is the average value of $S\,[g(t)]$.
If we set $g(t) = u(t)^{4 \over N-2} g_0$ where $u(t)$ is a positive function on $M$ for each $t \in (0,T)$, then \eqref{eq:Yamabef} is reduced to the quasi-linear parabolic equation
\begin{equation}\label{eq:Yamabef2}
p u^{p-1} u_t = \frac{N+2}{4} \(\kappa_N \Delta_{g_0} u - S\,[g_0]u + \ovs\,[g(t)] u^p\) \quad \text{on } M \times (0,T)
\end{equation}
where $\kappa_N := \frac{4(N-1)}{N-2}$ and $p := \frac{N+2}{N-2}$.

Let $E_{g_0}$ be the energy corresponding to equation \eqref{eq:Yamabef2},
\[E_{g_0}(u) := \frac{\int_M \(\kappa_N |\nabla_{g_0} u|_{g_0}^2 + S\,[g_0] u^2\) dv_{g_0}}{\(\int_M |u|^{2N \over N-2} dv_{g_0}\)^{N-2 \over N}}\]
where $dv_{g_0}$ is the volume form on $(M,g_0)$, and $Y(M,g_0)$ is the Yamabe constant on $(M,g_0)$ defined by
\[Y(M,g_0) = \text{inf} \left\{E_{g_0}(u): u \in C^{\infty}(M) \setminus \{0\} \right\}.\]

If $Y(M,g_0) \le 0$, a relatively simple argument based on the maximum principle shows that the flow \eqref{eq:Yamabef} exists globally $(T = \infty)$ and converges to a metric of constant scalar curvature as $t \to \infty$.

However, treating \eqref{eq:Yamabef} in the case $Y(M,g_0) > 0$ is challenging, because one has to exclude the possibility of blowing-up phenomena.
Thanks to the series of works by Chow \cite{Ch}, Ye \cite{Ye}, Schwetlick and Struwe \cite{SS}, and Brendle \cite{Br,Br2},
it is now known that \eqref{eq:Yamabef} is always globally well-defined in time and converges to a metric of constant scalar curvature as $t \to \infty$, under the validity of the positive mass theorem;
see also \cite{CCR} in which Carlotto, Chodosh, and Rubinstein showed the existence of a Yamabe flow on a certain manifold converging at a polynomial rate.

One of the main ingredients of the proof by Brendle in \cite{Br,Br2} is to construct a suitable family
$\left\{\bar{u}_{z,\ep} \in C^{\infty}(M):\, z \in M,\, \ep > 0 \text { small}\right\}$ of test functions satisfying
\begin{equation}\label{eq:test}
\begin{cases}
\displaystyle \sup \left\{E_{g_0}(\bar{u}_{z,\ep}): z \in M,\, \ep > 0 \text { small}\right\} \le Y(\S^N,g_{\S^N}),\\
\displaystyle \lim_{\ep \to 0} \ep^{N-2 \over 2} \bar{u}_{z,\ep}(\exp_z(\ep x)) = W_{1,0}(x) \quad \text{for } z \in M,\, x \in T_zM
\end{cases}
\end{equation}
where $(\S^N,g_{\S^N})$ is the $N$-dimensional unit-sphere in $\R^{N+1}$ with the canonical metric $g_{\S^N}$, $\exp$ is the exponential map, and $W_{1,0}$ is the standard bubble in $\R^N$, that is,
\begin{equation}\label{eq:bubble}
W_{1,0}(x) = \frac{\alpha_N}{(1+|x|^2)^{N-2 \over 2}} \quad \text{for } x \in \R^N, \quad \alpha_N := (N(N-2))^{N-2 \over 4}.
\end{equation}

\medskip
Given a smooth function $h$ on $M$ such that $\max_M h > 0$, let $E_{g_0,h}$ be a perturbed energy of $E_{g_0}$ given as
\[E_{g_0,h}(u) = \frac{\int_M \left[\kappa_N |\nabla_{g_0} u|_{g_0}^2 + (S\,[g_0]+h) u^2\right] dv_{g_0}}{\(\int_M |u|^{2N \over N-2} dv_{g_0}\)^{N-2 \over N}}.\]
If $h(z) > 0$ at some $z \in M$, then there is no test function $\bar{u}_{z,\ep}$ satisfying \eqref{eq:test}.
For instance, if $N \ge 7$, $M$ is non-locally conformally flat, and $\bar{u}_{z,\ep}$ is a `bubble-like' function satisfying the second condition in \eqref{eq:test}, then
\[E_{g_0,h}(\bar{u}_{z,\ep}) \approx Y(\S^N,g_{\S^N}) + \tc_1 h(z) \ep^2 - \tc_2 \|\text{Weyl}[g_0](z)\|_{g_0}^2 \ep^4\]
as $\ep \to 0$, so
\[E_{g_0,h}(\bar{u}_{z,\ep}) > Y(\S^N,g_{\S^N}) \quad \text{provided } \ep > 0 \text{ small enough}.\]
Here, $\tc_1, \tc_2 > 0$, $\text{Weyl}[g_0]$ is the Weyl tensor on $(M,g_0)$, and $\|\cdot\|_{g_0}$ is the tensor norm in the metric $g_0$.

Now setting
\begin{equation}\label{eq:Lgh}
L_{g_0,h}u = \kappa_N \Delta_{g_0} - (S\,[g_0]+h)
\end{equation}
so that $L_{g_0} := L_{g_0,h}$ is the conformal Laplacian on $(M,g_0)$, we consider a perturbed Yamabe flow
\begin{equation}\label{eq:Yamabefp}
\begin{cases}
\displaystyle pu^{p-1}u_t = \frac{N+2}{4} (L_{g_0,h}u + \kappa_N u^p) &\text{on } M \times (0,\infty), \quad p = \frac{N+2}{N-2}, \\
u > 0 &\text{on } M \times (0,\infty),\\
u(\cdot,0) = u_0 > 0 &\text{on } M,
\end{cases}
\end{equation}
corresponding to the energy $E_{g_0,h}$. Observe that we replaced the function $\ovs\,[g(t)]$ in \eqref{eq:Yamabef} with the constant $\kappa_N > 0$,
recalling that it tends to a positive constant as $t \to \infty$ provided $Y(M,g_0) > 0$.

A natural question is what is the asymptotic behavior of a solution to \eqref{eq:Yamabefp} as $t \to T$.
As we will prove in this paper, \eqref{eq:Yamabefp} may exhibit infinite-time blow-up phenomena driven by the bubbles (also known as the Talenti-Aubin instantons)
\begin{equation}\label{eq:bubble2}
W_{\mx}(x) = \mu^{-{N-2 \over 2}} W_{1,0}(\mu^{-1}(x-\xi)), \quad x \in \R^N,\, \mu > 0,\, \xi \in \R^N;
\end{equation}
see \eqref{eq:bubble}. It is a classical result that the solution set of the Yamabe problem in $\R^N$
\begin{equation}\label{eq:YamabeRN}
-\Delta W = W^p, \quad W > 0 \quad \text{in } \R^N
\end{equation}
is precisely $\{W_{\mx}: \mu > 0,\, \xi \in \R^N\}$,
which corresponds to a family of standard metrics on $N$-dimensional spheres in $\R^{N+1}$ via the stereographic projection.

\subsection{Main theorems}
We now list the main theorems of this paper and some remarks on them.

\medskip
The following theorem precisely describes the infinite-time blow-up phenomena of the perturbed Yamabe flow \eqref{eq:Yamabefp} driven by the bubbles.
\begin{thm}\label{thm:main}
Let $(M,g_0)$ be a smooth compact Riemannian manifold of dimension $N \ge 5$ such that
\begin{equation}\label{eq:Yamabepc}
Y_h(M,g_0) := \textnormal{inf} \left\{E_{g_0,h}(u): u \in C^{\infty}(M) \setminus \{0\} \right\} > 0.
\end{equation}
Assume that $h$ is a $C^2$-function on $M$ such that $\max_M h > 0$.
Given a point $z_0 \in M$ such that $h(z_0) > 0$, there exists an initial datum $u_{z_0,0}$
such that \eqref{eq:Yamabefp} has a solution $u_{z_0}$ blowing-up at $z_0$ as $t \to \infty$.
More precisely, there is a constant $C > 0$ depending only on $(M,g_0)$, $N$, $h$, and $z_0$ such that
\begin{equation}\label{eq:approx}
\begin{cases}
\displaystyle u_{z_0}(z,t) = \(1+O\(|\exp_{z_0}^{-1}(z)|^2\)\) W_{\mu(t),\xi(t)}\(\exp_{z_0}^{-1}(z)\) &\text{if } d_{g_0}(z,z_0) \le \delta_0,\\
\displaystyle u_{z_0}(z,t) \le C \mu^{N-2 \over 2}(t) G(z,z_0) &\text{if } d_{g_0}(z,z_0) > \delta_0
\end{cases}
\end{equation}
for $t > 0$ large enough, where
\begin{itemize}
\item[-] $\exp$ is the exponential map on $(M,g_0)$, and $d_{g_0}(z,z_0)$ is the geodesic distance between $z$ and $z_0$;
\item[-] $\delta_0 > 0$ is a sufficiently small number;
\item[-] $\mu: [0,\infty) \to (0,\infty)$ and $\xi: [0,\infty) \to \R^N$ are parameters
    such that $\mu(t) \simeq t^{-\frac{1}{2}}$ and $|\xi(t)| \lesssim t^{-1+\vep}$ for all $t > 0$ large and some $\vep \in (0,1)$ small;
\item[-] $G$ is the Green's function of the perturbed conformal Laplacian $L_{g_0,h}$ defined in \eqref{eq:Green}.
\end{itemize}
\end{thm}

In fact, equation \eqref{eq:Yamabefp} also possesses solutions which blow-up at multiple points as $t \to \infty$.
\begin{thm}\label{thm:main2}
Let $(M,g_0)$ be a smooth compact Riemannian manifold of dimension $N \ge 5$ such that $Y_h(M,g_0) > 0$.
Assume that $h$ is a $C^2$-function on $M$ such that $\max_M h > 0$.
Given any $k \in \N$, choose a $k$-tuple $\mbz_0 := \{z_0^{(1)}, \ldots, z_0^{(k)}\}$ of distinct points on $M$ such that $h(z_0^{(l)}) > 0$ for $l = 1, \ldots, k$.
Then there exists an initial datum $u_{\mbz_0,0}$ such that \eqref{eq:Yamabefp} has a solution $u_{\mbz_0}$ blowing-up at each point $z_0^{(1)}, \ldots, z_0^{(k)}$ as $t \to \infty$.
More precisely, there is a constant $C > 0$ depending only on $(M,g_0)$, $N$, $h$, $k$, and $\mbz_0$ such that
\[\begin{cases}
\displaystyle u_{\mbz_0}(z,t) = \Big(1 + O\Big(\Big|\exp_{z_0^{(l)}}^{-1}(z)\Big|^2\Big)\Big)
W_{\mu^{(l)}(t),\xi^{(l)}(t)}\(\exp_{z_0^{(l)}}^{-1}(z)\) &\text{if } d_{g_0}\(z,z_0^{(l)}\) \le \delta_0,\, l = 1, \ldots, k,\\
\displaystyle u_{\mbz_0}(z,t) \le C\mu^{N-2 \over 2}(t) \sum_{l=1}^k G\(z,z_0^{(l)}\) &\text{otherwise}
\end{cases}\]
for $t > 0$ large enough, where $\mu^{(l)}: [0,\infty) \to (0,\infty)$ and $\xi^{(l)}: [0,\infty) \to \R^N$ are parameters
such that $\mu^{(l)}(t) \simeq t^{-\frac{1}{2}}$ and $|\xi^{(l)}(t)| \lesssim t^{-1+\vep}$ for all $t > 0$ large and some $\vep \in (0,1)$ small.
\end{thm}

Regarding the stability of the solution found above, we have the following result.
\begin{cor}\label{cor:main}
In the setting of Theorem \ref{thm:main2}, we further assume that the largest eigenvalue of the Ricci curvature tensor at $z_0^{(l)}$
is less than or equal to $-\frac{6}{N-4}h(z_0^{(l)})$ for each $l = 1, \ldots, k$.
Let $\sigma \in (0,1)$. There is a $k$-codimensional manifold $\mcm_{\mbz_0}$ in $C^{2,\sigma}(M)$ containing $u_{\mbz_0}$ such that
if $u_0 \in \mcm_{\mbz_0}$ is sufficiently close to $u_{\mbz_0,0}$, then \eqref{eq:Yamabefp} has a solution $u$ blowing-up at each point $z_0^{(1)}, \ldots, z_0^{(k)}$ as $t \to \infty$.
\end{cor}
\noindent In the above corollary, the technical condition on the Ricci curvature comes from the ODE system \eqref{eq:xieq} that each $\xi^{(l)}$ solves.

\begin{rmk}\label{rmk:main}
Theorems \ref{thm:main} and \ref{thm:main2} show that the Yamabe flow is an equation at the borderline guaranteeing the global existence and uniform boundness of solutions.
A couple of remarks regarding them are in order.

\medskip \noindent
(1) If $Y(M,g_0) > 0$ and $h$ is a function such that $\|h\|_{L^{N/2}(M)} < Y(M,g_0)$, then a simple application of H\"older's inequality yields that $Y_h(M,g_0) > 0$.
Therefore, if $(M,g_0)$ is a smooth compact Riemannian manifold of dimension $N \ge 5$ such that $Y(M,g_0) > 0$,
\eqref{eq:Yamabefp} exhibits infinite-time blow-up phenomena (with an arbitrary number of blow-up points) provided $h \in C^2(M)$ is a sufficiently small and $\max_M h > 0$.

\medskip \noindent
(2) Our results can be regarded as parabolic analogues of the theorems of Micheletti, Pistoia, and V\'etois \cite{MPV}, and of Esposito, Pistoia, and V\'etois \cite{EPV}
which assert the existence of blow-up solutions to slightly perturbed (elliptic) Yamabe problems.
Their results are related to $C^2(M)$-compactness property of the solution set of the Yamabe problem on $(M,g_0)$; see \cite{KMS, BM} and references therein.

In the elliptic case, the blow-up points must be a critical point of some function involving the function $h$ and a geometric quantity on $(M,g_0)$;
either the Weyl curvature or the $O(1)$-order term of the Green's function of the conformal Laplacian on $(M,g_0)$.
In our evolutionary setting, we only require that $h$ is positive at each blow-up point.

We wonder if there exist non-isolated positive blowing-up solutions (namely, clustering solutions)
or isolated non-simple positive blowing-up solutions (namely, bubble-tower solutions) to \eqref{eq:Yamabefp} as in the elliptic case \cite{RV, PV, TV, Che, Pr}.
Also, we may ask if \eqref{eq:Yamabefp} has a finite-time blowing-up solution.
We hope to examine these problems elsewhere.

\medskip \noindent
(3) Our results are in spirit close to the work of Daskalopoulos, del Pino, and Sesum \cite{DdPS} which constructed type II ancient compact solutions to the Yamabe flow on $\S^N$,
and that of Sire, Wei, and Zheng \cite{SWZ} which built finite-time extinguishing solutions to fast diffusion equations with the critical exponent in smooth bounded domains in $\R^N$.

To handle the degeneracy of the equation ($u^{p-1}$ in \eqref{eq:Yamabefp} for our setting),
the authors in \cite{SWZ} lifted the spatial domain to a subset of $\S^N$ via the stereographic projection so that the equation becomes uniformly parabolic.
However, there are limitations in employing their idea directly, because our spatial domain is a general Riemannian manifold $M$ and so it may not be embedded into $\S^N$.
We will overcome this technical difficulty with the help of a maximum principle adapted to our setting; see Lemma \ref{lemma:mp}.

\medskip \noindent
(4) Suppose that $\min_M h < 0$. In this case, our approach suggests that given a $k$-tuple of points $\mbz_0 = \{z_0^{(1)}, \ldots, z_0^{(k)}\}$ on $M$ such that $h(z_0^{(l)}) < 0$ for $l = 1, \dots, k$,
there exist ancient solutions to \eqref{eq:Yamabefp} which blow-up at each point $z_0^{(1)}, \ldots, z_0^{(k)}$ as $t \to -\infty$.
We also suspect the existence of ancient solutions to \eqref{eq:Yamabefp} which look like towers of spheres similar to ones found in \cite{DdPS}.
\end{rmk}

There have been extensive studies on type II blow-up solutions to various time-dependent energy-critical equations;
e.g. the critical nonlinear heat equations \cite{FHV, Sch, CMR, CdPM, DMW, Ha, DMW2}, the harmonic map heat flows and the nematic liquid crystal flows \cite{CDY, BHK, To, RS, DdPW, LLWWZ},
the critical wave equations \cite{KM, KST, HR, DK, DKM, Je, JM},
the wave maps equations and the Yang-Mills equations \cite{RR, RS2}, the critical Schr\"odinger equations and map equations \cite{MR, MRR}, and so on.
Our contribution towards this direction is to build solutions to energy-critical degenerate parabolic equations on general Riemannian manifolds, blowing up at an arbitrary number of points.
We believe that our method may help construct bubbling solutions to related problems on manifolds,
such as the harmonic map heat flow from a smooth closed two-dimensional Riemannian manifold to another provided the target manifold satisfies a certain geometric assumption.

Furthermore, several results on the optimal regularity and sharp extinction rates for fast diffusion equations
in bounded domains in $\R^N$ were proved recently; refer to \cite{JX, JX2, BGV, BF} among others.
In their proofs, Aronson-B\'enilan type inequalities (namely, bounds for $u^{-1}|u_t|$; see \cite{AB}) appear as one of the key tools.
In our analysis of \eqref{eq:Yamabefp}, such an inequality is derived in a very specific manner; see \eqref{eq:ArBe} below and its proof.

\subsection{Structure of the paper and comment on the proof}
In order to establish Theorems \ref{thm:main} and \ref{thm:main2} as well as Corollary \ref{cor:main},
we will apply the modulation argument combined with the inner-outer gluing procedure.
Our presentation is largely influenced by the paper \cite{CdPM}
which studied the existence of infinite-time type II multiple blowing-up solutions to the critical nonlinear heat equations in smooth bounded domains in $\R^N$.

\medskip
From Section \ref{sec:approx} to \ref{sec:inner}, we only concern the proof of Theorem \ref{thm:main}.
Necessary modifications to prove Theorem \ref{thm:main2} and the proof of Corollary \ref{cor:main} will be provided in Section \ref{sec:multi}.

In Section \ref{sec:approx}, we construct approximate solutions which behave as in \eqref{eq:approx}.
This is one of the most tricky parts of the proof, because approximate solutions must be sufficiently regular
and very close to true solutions for every point on $M$ regardless of its distance from the blow-up point.
In particular, the approximate solutions used for the elliptic analogues do not work well in our setting.
At points away from the blow-up point, we depict the approximate solution in terms of the Green's function $G$ of the perturbed conformal Laplacian $L_{g_0,h}$; cf. \cite{Sc, EPV}.
At points near the blow-up point, we deform the bubble by multiplying it by $G$ and then combining the result and a solution to the linear equation \eqref{eq:Psi} together.
During the refinement, we also determine the main order term of the dilation parameter by solving the ODE \eqref{eq:orthoN+1}.

In Section \ref{sec:inout}, we decompose equation \eqref{eq:Yamabefp} into the outer problem \eqref{eq:outer} and the inner problem \eqref{eq:inner}.

In Section \ref{sec:outer}, we prove the unique solvability of the outer problem and establish a priori estimates of the solution,
by examining the associated inhomogeneous problem \eqref{eq:inhom-o} with the maximum principle for degenerate parabolic equations.
A delicate issue is to choose suitable norms to work with.
While the authors in \cite{CdPM} worked successfully in a weighted $L^{\infty}$ setting,
we cannot do for \eqref{eq:Yamabefp} because of its degeneracy $u^{p-1}$.
To control the terms resulting from the degeneracy such as $(\psi_{\mx})_t$ in \eqref{eq:mcj2} in a pointwise sense,
we will devise various H\"older-type norms; refer to Subsection \ref{subsec:norms}.
In defining them, we must reflect that the scaling properties in the spatial variable
and in the time variable are different from each other, which makes the analysis considerably complicated.

In Section \ref{sec:inner}, we develop the existence theory for an associated inhomogeneous problem \eqref{eq:inhom-i} to the inner problem.
By lifting its spatial domain to $\S^N$, we prove that it is solvable whenever the orthogonality condition \eqref{eq:ortho1} holds.
Also, from \eqref{eq:ortho1}, we find a system of ODEs that the dilation and translation parameters satisfy, and solve it.
Finally, collecting all the information obtained so far, we find a solution to the inner problem and derive a priori estimate on it, thereby completing the proof of Theorem \ref{thm:main}.

\subsection{Notations}
We collect some notations used throughout the paper.

\medskip \noindent - The Einstein convention is used throughout the paper.
Unless otherwise stated, the indices $i$, $j$, $q$, $r$, and $s$ take values from $1$ to the dimension $N$ of the underlying manifold $M$,
and $l$ and $m$ range over values from $1$ to the number $k$ of blowing-up points of solutions.

\medskip \noindent - For a tensor field $T$ on $(M,g)$, a notation such as $T_{;a}$ stands for a covariant derivative of $T$.

\medskip \noindent - $\S^N$ is the standard unit sphere in $\R^{N+1}$, $g_{\S^N}$ is its canonical metric, and $\left|\S^N\right|$ is its surface measure.

\medskip \noindent - $\nabla_g$, $\Delta_g$, $\la \cdot, \cdot \ra_g$, and $|\cdot|_g$ are the gradient, the Laplace-Beltrami operator, the inner product, and the norm on $(M,g)$, respectively.
On $(\S^N,g_{\S^N})$, we write $\nabla_{\S^N} = \nabla_{g_{\S^N}}$, etc.
In the Euclidean space $\R^N$, we write $\nabla = \nabla_x$, $\Delta = \Delta_x$, etc, where the subscript $x$ denotes the variable in $\R^N$.

\medskip \noindent - $i(M,g_0)$ is the injectivity radius of $(M,g_0)$.

\medskip \noindent - For a surface integral, $dS$ denotes the volume form on the domain of integration.

\medskip \noindent - Given a number $\delta > 0$ and a metric $g$ on $M$, let $B^N(x,\delta) = \{y \in \R^N: |y-x| < \delta\}$ for $x \in \R^N$
and $B_g(z,\delta) = \{\xi \in M: d_g(z,\xi) < \delta\}$ for $z \in M$.

\medskip \noindent - For a function $f$ and $\ell \in \N$, we set $\pa_i f = \frac{\pa f}{\pa x_i}$ for $i = 1, \ldots, N$, $\dot{f} = f_t = \pa_t f = \frac{\pa f}{\pa t}$, and
\begin{equation}\label{eq:f_ell}
f_{[\ell]}(x) = f(x) - \sum_{|\ga| = 0}^{\ell-1} \frac{x^{\ga}}{\ga!}\, \pa_{\ga} f(0) = O(|x|^{\ell})
\end{equation}
for $|x|$ small, where $\ga$ is a multi-index.

\medskip \noindent - For $\ell \in \N \cup \{0\}$, $\sigma \in (0,1)$, $T > 0$, and a set $\Omega$,
we write $C^{\ell+\sigma}(\Omega) = C^{\ell,\sigma}(\Omega)$ and $C^{2\ell+\sigma,\ell+\sigma/2}(\Omega \times [0,T])$ to refer the H\"older space and the parabolic H\"older space, respectively.

\medskip \noindent - $\text{supp}(f)$ is the support of a function $f$.

\medskip \noindent - Let $\eta \in C^{\infty}(\R)$ be a function such that
\begin{equation}\label{eq:eta}
\begin{cases}
\eta(r) \ge 0,\ \eta'(r) \le 0 \text{ for } r \in \R,\\
\eta(r) = 1 \text{ for } r \in (-\infty,1] \text{ and } 0 \text{ for } r \in [2,\infty),
\end{cases}
\end{equation}
and $\eta_{\delta}(r) = \eta(\delta^{-1}r)$ for $r \in \R$ and $\delta > 0$.
By abuse of notation, we often write $\eta_{\delta}(x) = \eta_{\delta}(|x|)$ for $x \in \R^N$.

\medskip \noindent - $C,\, \zeta > 0$ are universal constants that may vary from line to line.

\subsection{Norms}\label{subsec:norms}
We introduce all the norms that will be used throughout the paper.
Let $\delta_0 > 0$ be the small number in the statement of Theorem \ref{thm:main}, for which we impose that $4\delta_0 < i(M,g_0)$.
Let also $t_0 > 0$ be a large number, and $\mu = \mu(t) > 0$ and $\xi = \xi(t) \in \R^N$ be small functions on $[t_0,\infty)$ tending to $0$ as $t \to \infty$.
The functions $\mu_0 = \mu_0(t)$ and $u_{\mx} = u_{\mx}^{(2)}(z,t)$ are defined in \eqref{eq:mu0} and \eqref{eq:u2} below, respectively.

\begin{defn}[Local H\"older semi-norms]\label{defn:semi-norm}
Assume that $\sigma \in (0,1)$. For a function $\lambda: [t_0,\infty) \to \R$, we set
\[[\lambda]_{C^{\sigma/2}_t}(t) = \sup\left\{\frac{|\lambda(t_1) - \lambda(t_2)|}{|t_1-t_2|^{\sigma/2}}:\, t_1,t_2 \in (\max\{t_0,t-1\},t),\, t_1 \ne t_2\right\}.\]

For a function $\psi: M \times [t_0,\infty) \to \R$, we set
\[[\psi]_{C^{\sigma}_z}(z,t) = \sup\left\{\frac{|\psi(z_1,t) - \psi(z_2,t)|}{d_{g_0}(z_1,z_2)^{\sigma}}:\, z_1, z_2 \in B_{g_0}(z,r_0(z,t)),\, z_1 \ne z_2\right\}\]
and
\[[\psi]_{C^{\sigma/2}_t}(z,t) = \sup\left\{\frac{|\psi(z,t_1) - \psi(z,t_2)|}{|t_1-t_2|^{\sigma/2}}:\, t_1,t_2 \in (\max\{t_0,t-1\},t),\, t_1 \ne t_2\right\}\]
where
\[r_0(z,t) := \begin{cases}
\displaystyle \frac{\mu_0(t)}{2} &\text{for } z \in B_{g_0}(z_0,\mu_0(t)),\\
\displaystyle \frac{d_{g_0}(z,z_0)}{2} &\text{for } z \in B_{g_0}(z_0,\delta_0) \setminus B_{g_0}(z_0,\mu_0(t)),\\
\displaystyle \frac{\delta_0}{2} &\text{for } z \in M \setminus B_{g_0}(z_0,\delta_0),
\end{cases}\]
and $t \in [t_0,\infty)$.

Let $\Omega$ be a smooth domain in $\R^N$. For a function $\psi: \Omega \times [t_0,\infty) \to \R$, we set
\[[\psi]_{C^{\sigma}_{\Omega}}(y,t) = \sup\left\{\frac{|\psi(y_1,t) - \psi(y_2,t)|}{|y_1-y_2|^{\sigma}}:\, y_1, y_2 \in B^N(y,r_1(y)) \cap \Omega,\, y_1 \ne y_2\right\}\]
where
\[r_1(y) := \begin{cases}
\displaystyle \frac{1}{2} &\text{in } \Omega \cap \{|y| < 1\},\\
\displaystyle \frac{|y|}{2} &\text{in } \Omega \cap \{|y| \ge 1\}.
\end{cases}\]
\end{defn}

\begin{defn}[A H\"older norm for functions on $[t_0,\infty)$]\label{defn:norm0}
For numbers $\nu > 0$ and $\sigma \in (0,1)$ and functions $\lambda: [t_0,\infty) \to \R$ and $\xi = (\xi_1, \ldots, \xi_N): [t_0,\infty) \to \R^N$, we define
\begin{equation}\label{eq:nunorm}
\begin{cases}
\displaystyle \|\lambda\|_{\nu;\sigma} = \sup_{t \in [t_0,\infty)} \mu_0^{-\nu}(t) \left[|\lambda(t)| + [\lambda]_{C^{\sigma/2}_t}(t)\right], \\
\displaystyle \|\xi\|_{\nu;\sigma} = \sup_{t \in [t_0,\infty)} \mu_0^{-\nu}(t) \sum_{i=1}^N \left[|\xi_i(t)| + [\xi_i]_{C^{\sigma/2}_t}(t)\right].
\end{cases}
\end{equation}
\end{defn}

\begin{defn}[Weighted integral norms for functions on $M$ or $M \times [t_0,\infty)$]\label{defn:norm1}
Let $M_{\tau} = M \times [\tau,\tau+1]$ for $\tau \ge t_0$.

\medskip \noindent - Local in time weighted $L^2$, $H^1$ and $H^2$ norms: We set
\[\|\psi_0\|_{\mch^1} = \left[\int_M (|\nabla_{g_0} \psi_0|_{g_0}^2 + \psi_0^2) dv_{g_0}\right]^{1 \over 2}\]
for a function $\psi_0$ on $M$, and
\begin{align*}
\|\psi(\cdot,\tau)\|_{L^2(M)} &= \left[\int_M \(|\psi|^2 u_{\mx}^{p-1}\)(z,\tau) dv_{g_0}\right]^{1 \over 2}, \\
\|\psi\|_{L^2(M_{\tau})} &= \left[\iint_{M_{\tau}} \(|\psi|^2 u_{\mx}^{p-1}\)(z,t) dv_{g_0} dt\right]^{1 \over 2}, \\
\|\psi\|_{H^1(M_{\tau})} &= \left\|u_{\mx}^{-{p-1 \over 2}} \nabla_{g_0} \psi\right\|_{L^2(M_{\tau})} + \|\psi\|_{L^2(M_{\tau})}, \\
\|\psi\|_{H^2(M_{\tau})} &= \|\psi_t\|_{L^2(M_{\tau})} + \left\|u_{\mx}^{1-p}\, L_{g_0,h} \psi\right\|_{L^2(M_{\tau})} + \|\psi\|_{H^1(M_{\tau})}
\end{align*}
for a function $\psi$ on $M \times [t_0,\infty)$.

\medskip \noindent - Global in time weighted $L^2$, $H^1$ and $H^2$ norms: Given two positive numbers $t_0$ and $s_0$ such that $t_0 < s_0-1$, we set
\begin{align*}
\|\psi\|_{L^2_{t_0,s_0}} = \sup_{\tau \in [t_0,s_0-1]} \|\psi\|_{L^2(M_{\tau})}, \\
\|\psi\|_{H^1_{t_0,s_0}} = \sup_{\tau \in [t_0,s_0-1]} \|\psi\|_{H^1(M_{\tau})}, \\
\|\psi\|_{H^2_{t_0,s_0}} = \sup_{\tau \in [t_0,s_0-1]} \|\psi\|_{H^2(M_{\tau})}.
\end{align*}
In the case that $s_0 = \infty$, we omit the subscript $s_0$ so that $\|\psi\|_{L^2_{t_0,\infty}} = \|\psi\|_{L^2_{t_0}}$, etc.
\end{defn}
\noindent Observe that $\|\psi\|_{H^2(M_{\tau})}$ is a norm if $Y_h(M,g_0) > 0$ holds; see \eqref{eq:Yamabepc}.

\begin{defn}[Weighted H\"older norms for functions on $M \times [t_0,\infty)$]\label{defn:norm2}
Given $\alpha \in \R$ and $\ga \ge 0$, we define
\begin{equation}\label{eq:wag}
w_{\alpha,\ga}(z,t) = \mu^{\alpha}(t) \cdot
\max \left\{\frac{\eta_{\delta_0}(d_{g_0}(z,z_0))}{\mu^{\alpha+\ga}(t) + |\exp_{z_0}^{-1}(z)-\xi(t)|^{\alpha+\ga}},\,
2^{3(\alpha+\ga)} \delta_0^{-(\alpha+\ga)}\right\}
\end{equation}
on $M \times [t_0,\infty)$. Here, $\eta_{\delta_0}$ is the cut-off function defined after \eqref{eq:eta}.

\medskip \noindent - The weighted $L^{\infty}$ norms: Given $\alpha > 0$ and $\rho \ge -\frac{N-2}{2}$, we define
\[\|f\|_{*,\ar} = \left\|\(\mu_0^{\rho} w_{\alpha,2}\)^{-1} u_{\mx}^{p-1}f\right\|_{L^{\infty}(M \times [t_0,\infty))},
\quad \|\psi_0\|_{**,\alpha} = \left\|\(\mu_0^{\rho}(t_0)w_{\alpha,0}\)^{-1} \psi_0\right\|_{L^{\infty}(M)},\]
and
\[\|\psi\|_{*',\ar} = \left\|\(\mu_0^{\rho} w_{\alpha,0}\)^{-1} \psi\right\|_{L^{\infty}(M \times [t_0,\infty))}.\]

\medskip \noindent - The weighted H\"older norms: Given $\alpha > 0$, $\rho \ge -\frac{N-2}{2}$, and $\sigma \in (0,1)$, we define
\begin{equation}\label{eq:*-norm}
\begin{medsize}
\|f\|_{*,\ar;\sigma} = \left\|\(\mu_0^{\rho} w_{\alpha,2}\)^{-1} \left|u_{\mx}^{p-1}f\right| + \(\mu_0^{\rho} w_{\alpha,2+\sigma}\)^{-1} \left[u_{\mx}^{p-1}f\right]_{C^{\sigma}_z}
+ \(\mu_0^{\rho} w_{\alpha-\sigma,2}\)^{-1} \left[u_{\mx}^{p-1}f\right]_{C^{\sigma/2}_t} \right\|_{L^{\infty}(M \times [t_0,\infty))},
\end{medsize}
\end{equation}
\begin{equation}\label{eq:**-norm}
\|\psi_0\|_{**,\alpha;\sigma} = \sum_{\ell=0}^2 \left\|\(\mu_0^{\rho}(t_0) w_{\alpha,\ell}\)^{-1} \left|\nabla_{g_0}^{\ell} \psi_0\right|
+ \(\mu_0^{\rho}(t_0) w_{\alpha,\ell+\sigma}\)^{-1} \left[\nabla_{g_0}^{\ell} \psi_0\right]_{C^{\sigma}_z} \right\|_{L^{\infty}(M)},
\end{equation}
and
\begin{equation}\label{eq:*'s-norm}
\begin{aligned}
&\, \begin{medsize}
\|\psi\|_{*',\ar;\sigma}
\end{medsize} \\
&\begin{medsize}
\displaystyle = \sum_{\ell=0}^2 \left\|\(\mu_0^{\rho} w_{\alpha,\ell}\)^{-1} \left|\nabla_{g_0}^{\ell} \psi\right|
+ \(\mu_0^{\rho} w_{\alpha,\ell+\sigma}\)^{-1} \left[\nabla_{g_0}^{\ell} \psi\right]_{C^{\sigma}_z}
+ \(\mu_0^{\rho} w_{\alpha-\sigma,\ell}\)^{-1} \left[\nabla_{g_0}^{\ell} \psi\right]_{C^{\sigma/2}_t}\right\|_{L^{\infty}(M \times [t_0,\infty))}
\end{medsize} \\
&\, \begin{medsize}
+ \left\|\(\mu_0^{\rho} w_{\alpha-2,0}\)^{-1} |\psi_t| + \(\mu_0^{\rho} w_{\alpha-2,\sigma}\)^{-1} [\psi_t]_{C^{\sigma}_z}
+ \(\mu_0^{\rho} w_{\alpha-2-\sigma,0}\)^{-1} [\psi_t]_{C^{\sigma/2}_t} \right\|_{L^{\infty}(M \times [t_0,\infty))}.
\end{medsize}
\end{aligned}
\end{equation}
\end{defn}

\begin{defn}[Weighted H\"older norms for functions on $\Omega \times [t_0,\infty)$]\label{defn:norm3}
Let $\Omega$ be a smooth domain in $\R^N$, $a \in (0,N-2)$, $b > 0$, and $s_0 > \frac{3t_0}{2} \gg 1$.

\medskip \noindent - The weighted $L^{\infty}$ norms: We set
\[\|\mch\|_{\sharp,a+2,b;t_0,s_0} = \left\|\mu_0^{-b} \(1+|\by|^{a+2}\) \(W_{1,0}^{p-1} \mch\)(\by,t)\right\|_{L^{\infty}(\R^N \times [t_0,s_0))}\]
and
\[\|\psi\|_{\sharp',a,b;t_0,s_0(\Omega)} = \left\|\mu_0^{-b} \(1+|\by|^a\) \psi(\by,t)\right\|_{L^{\infty}(\Omega \times [t_0,\infty))}.\]
In the case that $s_0 = \infty$, we write $\|\mch\|_{\sharp,a+2,b} = \|\mch\|_{\sharp,a+2,b;t_0,s_0}$
and $\|\psi\|_{\sharp',a,b(\Omega)} = \|\psi\|_{\sharp',a,b;t_0,s_0(\Omega)}$.

\medskip \noindent - The weighted H\"older norms: Given $\sigma \in (0,1)$, we set
\begin{multline}\label{eq:sharp-norm}
\begin{medsize}
\displaystyle \|\mch\|_{\sharp,a+2,b;\sigma} =
\left\|\mu_0^{-b} \left\{\(1+|\by|^{a+2}\) \left|W_{1,0}^{p-1} \mch\right|(\by,t) \right. \right.
\end{medsize} \\
\begin{medsize}
\displaystyle \left. \left. + \(1+|\by|^{a+2+\sigma}\) \left[W_{1,0}^{p-1} \mch\right]_{C^{\sigma}_{\R^N}}(\by,t)
+ \(1+|\by|^{a+2-\sigma}\) \left[W_{1,0}^{p-1} \mch\right]_{C^{\sigma/2}_t}(\by,t) \right\}\right\|_{L^{\infty}(\R^N \times [t_0,\infty))}
\end{medsize}
\end{multline}
and
\begin{equation}\label{eq:sharp'-norm}
\begin{aligned}
&\, \begin{medsize}
\displaystyle \|\psi\|_{\sharp',a,b;\sigma(\Omega)}
\end{medsize} \\
&\begin{medsize}
\displaystyle = \sum_{\ell=0}^2 \left\|\mu_0^{-b} \left\{\(1+|\by|^{a+\ell}\) \left|\nabla_{\by}^{\ell} \psi(\by,t)\right|
+ \(1+|\by|^{a+\ell+\sigma}\) \left[\nabla_{\by}^{\ell} \psi\right]_{C^{\sigma}_{\Omega}}(\by,t)
\right. \right.
\end{medsize} \\
&\hspace{227pt} \begin{medsize}
\left. \left. + \(1+|\by|^{a+\ell-\sigma}\) \left[\nabla_{\by}^{\ell} \psi\right]_{C^{\sigma/2}_t}(\by,t) \right\}\right\|_{L^{\infty}(\Omega \times [t_0,\infty))}
\end{medsize} \\
&\, \begin{medsize}
\displaystyle + \left\|\mu_0^{-b} \left\{\(1+|\by|^{a-2}\) \left|\psi_t(\by,t)\right| + \(1+|\by|^{a-2+\sigma}\) [\psi_t]_{C^{\sigma}_{\Omega}}(\by,t)
+ \(1+|\by|^{a-2-\sigma}\) [\psi_t]_{C^{\sigma/2}_t}(\by,t) \right\}\right\|_{L^{\infty}(\Omega \times [t_0,\infty))}.
\end{medsize}
\end{aligned}
\end{equation}
\end{defn}

\section{Construction of approximate solutions} \label{sec:approx}
From Section \ref{sec:approx} to Section \ref{sec:inner}, we will concern the proof of Theorem \ref{thm:main},
which asserts the existence of solutions to \eqref{eq:Yamabefp} blowing-up at a single point $z_0 \in M$ as $t \to \infty$.
During the proof, we will consider the equation
\begin{equation}\label{eq:Yamabefp2}
\begin{cases}
\displaystyle pu^{p-1}u_t = \frac{N+2}{4} (L_{g_0,h}u + \kappa_N u^p) &\text{on } M \times (t_0,\infty),\\
u > 0 &\text{on } M \times (t_0,\infty),\\
u(\cdot,t_0) = u_0 > 0 &\text{on } M
\end{cases}
\end{equation}
for $t_0 > 0$ large enough. Clearly, a time translation $u(\cdot-t_0)$ of a solution $u$ to \eqref{eq:Yamabefp2} solves \eqref{eq:Yamabefp}.

\medskip
In this section, we construct approximate solutions through two stages.

In Subsection \ref{subsec:approx1}, we define the first approximate solution which resembles to a bubble near $z_0$
and to the Green's function $G$ of the perturbed conformal Laplacian $L_{g_0,h}$ away from $z_0$.

In Subsection \ref{subsec:approx2}, we refine the first approximate solutions in a small neighborhood of $z_0$, by attaching solutions of certain linearized equations to the approximate solutions.
In order to make the linearized equations solvable, we choose the main order of the dilation factor $\mu(t) > 0$ suitably.

\subsection{First approximate solutions} \label{subsec:approx1}
Throughout the paper, we always assume that $N \ge 5$.

\medskip
Let $G$ be the Green's function of $L_{g_0,h}$, i.e.,
\begin{equation}\label{eq:Green}
- L_{g_0,h} G(z,z_0) = \delta_{z_0}(z) \quad \text{on } M
\end{equation}
where $\delta_{z_0}$ is the Dirac measure supported at $z_0 \in M$.

Given a pair of parameters $(\mx) = (\mu(t),\xi(t))$, we are going to define the first approximate solution $u^{(1)}_{\mx}$ to \eqref{eq:Yamabefp}.

\medskip \noindent
\textsc{Expression of $\mu(t)$}: Let
\begin{equation}\label{eq:ZN+1}
Z_{N+1}(y) = y \cdot \nabla W_{1,0}(y) + \frac{N-2}{2} W_{1,0}(y) = \alpha_N \(\frac{N-2}{2}\) \frac{1-|y|^2}{(1+|y|^2)^{N \over 2}} \quad \text{for } y \in \R^N
\end{equation}
and positive numbers
\begin{equation}\label{eq:c1c2}
c_1 = p \int_{\R^N} W_{1,0}^{p-1} Z_{N+1}^2 \quad \text{and} \quad c_2 = - \int_{\R^N} W_{1,0} Z_{N+1}.
\end{equation}
Reminding the hypothesis that $h(z_0) > 0$, we set
\begin{equation}\label{eq:mu0}
d_0 = \frac{1}{\sqrt{h(z_0)}} \quad \text{and} \quad \mu_0(t) = \sqrt{\frac{2c_1}{(N+2)c_2 t}} \quad \text{for } t \in [t_0,\infty)
\end{equation}
where $t_0 > 0$ is a sufficiently large number. Then we assume that $\mu(t)$ has the form
\begin{equation}\label{eq:mu}
\mu(t) = d_0\mu_0(t) + \lambda(t) =: \bmu(t) + \lambda(t) \quad \text{for } t \in [t_0,\infty)
\end{equation}
where $\lambda(t)$ is a higher-order term.

\medskip \noindent
\textsc{Expression of $\lambda(t)$ and $\xi(t)$}: The parameters $\lambda(t) \in \R$ as well as
$\xi(t) = (\xi_1, \ldots, \xi_N)(t) \in \R^N$ will be determined in Subsection \ref{subsec:inner2}. Until then, we assume that
\begin{equation}\label{eq:lx}
\begin{aligned}
\|\lambda\|_{\nu_1;\sigma_0} + \|\dot{\lambda}\|_{\nu_1+2;\sigma_0} + \|\xi\|_{\nu_2;\sigma_0} + \|\dot{\xi}\|_{\nu_2+2;\sigma_0} \le C
\end{aligned}
\end{equation}
where $\nu_1 = \nu_2 = 2-\vep_0$ for some small numbers $\vep_0 \in (0,1)$ and $\sigma_0 \in (0,1)$; refer to \eqref{eq:nunorm} for the definition of the norms.

\medskip \noindent
\textsc{Expression of approximate solutions}: From now on, we will often identify points $z \in B_{g_0}(z_0,2\delta_0)$
and $x = \exp_{z_0}^{-1}(z) \in B^N(0,2\delta_0)$ via $g_0$-normal coordinates centered at $z_0$.
Given a pair $(\mx)$ satisfying \eqref{eq:mu0}--\eqref{eq:lx}, we define
\begin{equation}\label{eq:u1}
u^{(1)}_{\mx}(z,t) = \ga_N^{-1} \kappa_N G(z,z_0) v^{(1)}_{\mx}(z,t) \quad \text{on } M \times [t_0,\infty)
\end{equation}
where
\begin{equation}\label{eq:v1}
\begin{aligned}
&\ v^{(1)}_{\mx}(z,t) \\
&= \begin{cases}
\displaystyle |x|^{N-2} \mu^{-{N-2 \over 2}} W_{1,0}(\mu^{-1}(x-\xi)) &\text{if } |x| = d_{g_0}(z,z_0) \le \delta_0, \\
\displaystyle \alpha_N \mu^{N-2 \over 2} \left[1 + \eta_{\delta_0}(x) \left\{\frac{|x|^{N-2}}{(\mu^2+|x-\xi|^2)^{(N-2)/2}} - 1\right\}\right] &\text{if } \delta_0 < |x| = d_{g_0}(z,z_0) \le 2\delta_0, \\
\displaystyle \alpha_N \mu^{N-2 \over 2} &\text{if } d_{g_0}(z,z_0) > 2\delta_0.
\end{cases}
\end{aligned}
\end{equation}
Here, $\kappa_N,\, \alpha_N > 0$ are the numbers appearing in \eqref{eq:Yamabef2} and \eqref{eq:bubble}, and $\ga_N := [(N-2) |\S^{N-1}|]^{-1}$.
\begin{lemma}
The function $u^{(1)}_{\mx}$ is of class $C^{2+\sigma_0,1+\sigma_0/2}(M \times [t_0,\infty))$.
Moreover, it is positive on $M \times [t_0,\infty)$ provided $t_0 > 0$ large enough.
\end{lemma}
\begin{proof}
Since $Y_h(M,g_0) > 0$, the principal $L^2(M)$-eigenvalue of the operator $L_{g_0,h}$ is positive, and so $L_{g_0,h}$ satisfies the maximum principle on $M$.
Therefore, the Green's function $G$ of $L_{g_0,h}$ is uniquely determined and positive on $M$.
From this fact and \eqref{eq:mu0}--\eqref{eq:lx}, we deduce the assertion.
\end{proof}
\noindent We remark that the definition of $u^{(1)}_{\mx}$ is motivated by Schoen \cite{Sc}.

\medskip
By slightly modifying the proof of \cite[Lemma 6.4]{LP} and taking $\delta_0$ smaller if needed, we see
\[\ga_N^{-1} \kappa_N |x|^{N-2} G(z,z_0) = 1 + \sum_{\ell=2}^{N-3} P_{\ell}(x) + c|x|^{N-2} \log|x| + O\(|x|^{N-2}\)\]
for $z = \exp_{g_0}(x) \in B_{g_0}(z_0,\delta_0)$ where $P_{\ell}$ is a homogeneous polynomial of degree $\ell$, and $c$ is a number that is zero for $N$ odd.
\begin{lemma}
It holds that
\begin{equation}\label{eq:P_2}
P_2(x) = \hc_1 R_{ij}(z_0)x_ix_j + \hc_2 S(z_0) |x|^2 + \hc_3 h(z_0) |x|^2
\end{equation}
where
\begin{equation}\label{eq:hc}
\hc_1 := \frac{1}{12},\quad \hc_2 := - \frac{1}{24(N-1)},\quad \hc_3 := -\frac{N-2}{8(N-1)(N-4)},
\end{equation}
$R_{ij}(z_0)$ is a component of the Ricci curvature tensor at $z_0$ on $(M,g_0)$, and $S(z_0)$ is the scalar curvature at $z_0$ on $(M,g_0)$.
\end{lemma}
\begin{proof}
We continue to identify two points $z \in B_{g_0}(z_0,2\delta_0)$ and $x = \exp_{z_0}^{-1}(z) \in B^N(0,2\delta_0)$. It is well-known that
\begin{equation}\label{eq:metric}
\begin{cases}
\displaystyle \sqrt{|g_0|}(x) = 1 - \frac{1}{6} R_{ij}(z_0) x_ix_j + \sqrt{|g_0|}_{[3]}(x), \\
\displaystyle g_0^{ij}(x) = \delta_{ij} + \frac{1}{3} R_{iqjr}(z_0) x_qx_r + (g_0^{ij})_{[3]}(x)
\end{cases}
\end{equation}
where
\begin{itemize}
\item[-] $R_{iqjr}(z_0)$ is a component of the Riemannian curvature tensor at $z_0$ on $(M,g_0)$;
\item[-] $\sqrt{|g_0|}_{[3]}$ and $(g_0^{ij})_{[3]}$ are the remainder terms in the Taylor expansions of $\sqrt{|g_0|}$ and $g_0^{ij}$, respectively, defined by \eqref{eq:f_ell}.
\end{itemize}

We set
\begin{equation}\label{eq:P_20}
P(x) = \ga_N^{-1} \kappa_N |x|^{N-2} G(x,z_0) - 1 \quad \text{for } x \in B^N(0,2\delta_0)
\end{equation}
and write $P_2(x) = c_{ij}x_ix_j$ for some $c_{ij} \in \R$ so that
\begin{equation}\label{eq:P_21}
P(x) = c_{ij}x_ix_j + P_{[3]}(x) = c_{ij}x_ix_j + O(|x|^3).
\end{equation}
Since $\ga_N \Delta |x|^{2-N} = L_{g_0,h} G(z,z_0) = \delta_{z_0}(z)$, it holds that
\begin{equation}\label{eq:P_22}
|x|^2 \Delta P - 2(N-2) x_i \pa_i P + |x|^N \(\Delta_{g_0} - \Delta\) \left\{(1+P)|x|^{2-N}\right\} - \kappa_N^{-1} |x|^2 (S+h)(1+P) = 0.
\end{equation}
Putting \eqref{eq:metric} and \eqref{eq:P_21} into \eqref{eq:P_22}, we obtain
\[-4(N-2) c_{ij} x_ix_j + 2 \sum_{i=1}^N c_{ii} |x|^2 + \frac{N-2}{3} R_{ij}(z_0) x_ix_j - \kappa_N^{-1} (S(z_0) + h(z_0))|x|^2 = O(|x|^3).\]
Solving this identity, we deduce \eqref{eq:P_2}.
\end{proof}

Define the error $\mcs(u)$ of a positive function $u$ on $M \times [t_0,\infty)$ as
\[\mcs(u) = -pu^{p-1}u_t + \frac{N+2}{4} \(L_{g_0,h}u + \kappa_N u^p\)\]
where $p = \frac{N+2}{N-2}$.
In the next lemmas, we compute the error of the first approximate solution $u^{(1)}_{\mx}$.
\begin{lemma}
For $z = \exp_{z_0}(x) = \exp_{z_0}(\mu y+\xi) \in B_{g_0}(z_0,\delta_0)$ and $t \in [t_0,\infty)$,
\begin{equation}\label{eq:S1in}
\mu^{N+2 \over 2} \mcs\(u^{(1)}_{\mx}\)(y,t) = \mce_0[\mu](y,t) + \mce_1[\mx](y,t).
\end{equation}
Here,
\begin{equation}\label{eq:E0}
\mce_0[\mu](y,t) := \mu^{-1} \dot{\mu} \(p W_{1,0}^{p-1}Z_{N+1}\)(y) + \mu^2 \mcf_0(y)
\end{equation}
and
\begin{multline}\label{eq:E1}
\mce_1[\mx](y,t) := \mu^{-1} \dot{\xi} \cdot \(p W_{1,0}^{p-1} \nabla W_{1,0}\)(y) + \mcf_1[\mx](y,t) \\
+ \mu^3 a^{\{1\}} + \mu^{\nu_2+2} a^{\{2\}} + \mu\dot{\mu} a^{\{0\}} + \mu\dot{\xi} \cdot \mba^{\{-1\}}
\end{multline}
where
\begin{equation}\label{eq:F0}
\begin{aligned}
\mcf_0(y) &:= \frac{(N+2)\kappa_N}{4} \left[\left\{\(2\hc_1 + 2N\hc_2-\kappa_N^{-1}\) S(z_0) + \(2N\hc_3-\kappa_N^{-1}\)h(z_0)\right\} W_{1,0}(y) \right. \\
&\hspace{80pt} + 4 \left\{\hc_2 S(z_0) + \hc_3 h(z_0)\right\} y_j\pa_j W_{1,0}(y) \\
&\hspace{80pt} \left. + \frac{4}{N-2} \left\{\hc_1R_{ij}(z_0)y_iy_j + \hc_2 S(z_0)|y|^2 + \hc_3 h(z_0)|y|^2\right\} W_{1,0}^p(y) \right],
\end{aligned}
\end{equation}
\begin{align}
\mcf_1[\mx](y,t) &:= (N+2)\kappa_N \mu \left[\left\{\hc_2 S(z_0) + \hc_3 h(z_0)\right\} \xi_j \pa_j W_{1,0}(y) \right. \nonumber \\
&\hspace{77pt} \left. + \frac{2}{N-2} \left\{\hc_1R_{ij}(z_0)y_i\xi_j + \hc_2 S(z_0)y_i\xi_i + \hc_3 h(z_0) y_i\xi_i\right\} W_{1,0}^p(y) \right], \label{eq:F1}
\end{align}
and $\hc_1,\, \hc_2$, and $\hc_3$ are constants defined in \eqref{eq:hc}.
Furthermore, $a^{\{\ga\}} \in \R$ and $\mba^{\{\ga\}} \in \R^N$ are $C^1$ functions of $(y,\mu,\mu^{-1}\xi)$ such that
\begin{equation}\label{eq:adecay}
\sum_{\ell=0}^1 \(1+|y|^{N-\ga+\ell}\) \(\left|\nabla_y^{\ell} a^{\{\ga\}}\right| + \left|\nabla_y^{\ell} \mba^{\{\ga\}}\right|\)
+ \mu_0^{-2} \(1+|y|^{N-\ga}\) \(\left|\pa_t a^{\{\ga\}}\right| + \left|\pa_t \mba^{\{\ga\}}\right|\)
\end{equation}
is uniformly bounded.
\end{lemma}
\begin{proof}
We recall the functions $P$, $W_{\mx}$, and $P_2$ given in \eqref{eq:P_20}, \eqref{eq:bubble2}, and \eqref{eq:P_2}, respectively.

Let us write $u^{(1)}_{\mx} = (1+P) W_{\mx}$ on $B_{g_0}(z_0,\delta_0) \times [t_0,\infty)$. Then
\begin{multline}\label{eq:error111}
\mcs\(u^{(1)}_{\mx}\) = -p(1+P)^p \pa_t W_{\mx} + \frac{N+2}{4} \left[\kappa_N \Delta_{g_0} \((1+P) W_{\mx}\) \right. \\
\left. + \kappa_N (1+P)^p W_{\mx}^p - (S+h)(1+P) W_{\mx}\right]
\end{multline}
in the $(x,t)$-variable. Also, thanks to \eqref{eq:mu}--\eqref{eq:lx}, we see
\begin{multline}\label{eq:error112}
-p(1+P)^p \pa_t W_{\mx} \\
= p \mu^{-{N+4 \over 2}} W_{1,0}^{p-1}(y) \left[\dot{\mu} Z_{N+1}(y) + \dot{\xi} \cdot \nabla W_{1,0}(y) \right] + \mu^{-{N \over 2}} \dot{\mu} a^{\{0\}} + \mu^{-{N \over 2}} \dot{\xi} \cdot \mba^{\{-1\}}.
\end{multline}

Setting
\begin{align*}
\mcr &= \left[\kappa_N \Delta_{g_0} \((1+P) W_{\mx}\) + \kappa_N (1+P)^p W_{\mx}^p - (S+h)(1+P) W_{\mx}\right] \\
&\ - \left[\kappa_N (\Delta P_2) W_{\mx} + 2\kappa_N \nabla P_2 \cdot \nabla W_{\mx} - (S+h)(z_0) W_{\mx} \right.\\
&\hspace{178pt} \left. - \frac{\kappa_N}{3} R_{ij}(z_0)x_i \pa_jW_{\mx} + \kappa_N (p-1) P_2^p W_{\mx}^p\right],
\end{align*}
we claim that
\begin{equation}\label{eq:mcr}
\mu^{N+2 \over 2} \mcr(y,t) = \mu^3 a^{\{1\}} + \mu^{\nu_2+2} a^{\{2\}}.
\end{equation}
By \eqref{eq:YamabeRN}, it holds that
\begin{equation}\label{eq:mcr2}
\begin{aligned}
\mcr &= \kappa_N (\Delta_{g_0} P - \Delta P_2) W_{\mx} + 2\kappa_N \(\la \nabla_{g_0} P, \nabla_{g_0} W_{\mx}\ra_{g_0} - \nabla P_2 \cdot \nabla W_{\mx}\) \\
&\ - \left[(S+h)(1+P) - (S+h)(z_0)\right] W_{\mx} \\
&\ + \kappa_N \left[(1+P) (\Delta_{g_0}-\Delta) W_{\mx} + \frac{1}{3} R_{ij}(z_0)x_i \pa_jW_{\mx} \right] \\
&\ + \kappa_N \left[(1+P)^p-(1+P)-(p-1)P_2\right] W_{\mx}^p.
\end{aligned}
\end{equation}
From this, we readily deduce that \eqref{eq:mcr} holds for $|y| \le 1$.
Suppose that $|y| \ge 1$. Putting \eqref{eq:bubble}, \eqref{eq:bubble2}, and the estimate
\[\left[(1+P)^p-(1+P)-(p-1)P_2\right] W_{\mx}^p = \mu^{-{N-4 \over 2}} a^{\{1\}}\]
into \eqref{eq:mcr2} yields
\begin{align*}
\(\alpha_N \mu^{N-2 \over 2}\)^{-1} \mcr &= \kappa_N (\Delta_{g_0} P - \Delta P_2) |x|^{2-N}
+ 2\kappa_N \(\la \nabla_{g_0} P, \nabla_{g_0} |x|^{2-N}\ra_{g_0} - \nabla P_2 \cdot \nabla |x|^{2-N}\) \\
&\ - \left[(S+h)(1+P) - (S+h)(z_0)\right] |x|^{2-N} \\
&\ + \kappa_N \left[(1+P) (\Delta_{g_0}-\Delta) |x|^{2-N} + \frac{1}{3} R_{ij}(z_0)x_i \pa_j|x|^{2-N} \right] \\
&\ + \mu^{-(N-3)} a^{\{1\}} + \mu^{-(N-2)+\nu_2} a^{\{2\}}.
\end{align*}
Then \eqref{eq:P_22} and \eqref{eq:P_2}--\eqref{eq:hc} imply
\begin{align*}
\mu^{N+2 \over 2} \mcr(y,t) &= - \mu^N \alpha_N \left[\kappa_N \Delta P_2 |x|^{2-N} + 2\kappa_N \nabla P_2 \cdot \nabla |x|^{2-N} - (S+h)(z_0) |x|^{2-N} \right. \\
&\hspace{55pt} \left. + \frac{(N-2)\kappa_N}{3} R_{ij}(z_0)x_ix_j |x|^{2-N}\right] + \mu^3 a^{\{1\}} + \mu^{\nu_2+2} a^{\{2\}} \\
&= - \mu^N \alpha_N \kappa_N \underbrace{\left[2\hc_1 + 2N\hc_2 - 4(N-2)\hc_2 - \kappa_N^{-1}\right]}_{=0} S(z_0) |x|^{2-N} \\
&\ - \mu^N \alpha_N \kappa_N \underbrace{\left[2N\hc_3 - 4(N-2)\hc_3 - \kappa_N^{-1}\right]}_{=0} h(z_0) |x|^{2-N} + \mu^3 a^{\{1\}} + \mu^{\nu_2+2} a^{\{2\}} \\
&= \mu^3 a^{\{1\}} + \mu^{\nu_2+2} a^{\{2\}}.
\end{align*}
We note that the term involving $R_{ij}(z_0)x_i x_j |x|^{-N}$ vanished here, because its coefficient is a multiple of $4\hc_1-\frac{1}{3} = 0$.
Therefore, the assertion \eqref{eq:mcr} is true.

Now, by applying \eqref{eq:P_2} once more, we derive
\[(\Delta P_2) W_{\mx} = \mu^{-{N-2 \over 2}} \left[\(2\hc_1 + 2N\hc_2\) S(z_0) + 2N\hc_3 h(z_0)\right] W_{1,0}(y)\]
and
\[\nabla P_2 \cdot \nabla W_{\mx} = 2\mu^{-{N \over 2}} \left[\hc_1 R_{ij}(z_0) (\mu y_i+\xi_i) + (\hc_2 S(z_0) + \hc_3 h(z_0)) (\mu y_j+\xi_j)\right] \pa_j W_{1,0}(y).\]
Accordingly,
\begin{equation}\label{eq:error113}
\begin{aligned}
&\ \Delta_{g_0} \((1+P) W_{\mx}\) + (1+P)^p W_{\mx}^p - \kappa_N^{-1} (S+h)(1+P) W_{\mx} \\
&= (\Delta P_2) W_{\mx} + 2\nabla P_2 \cdot \nabla W_{\mx} - \kappa_N^{-1} (S+h)(z_0) W_{\mx}
- \frac{1}{3} R_{ij}(z_0)x_i \pa_jW_{\mx} + (p-1)P_2W_{\mx}^p \\
&\ + \mu^{-{N-4 \over 2}} a^{\{1\}} + \mu^{\nu_2-{N-2 \over 2}} a^{\{2\}} \\
&= \mu^{-{N-2 \over 2}} \left[\(2\hc_1 + 2N\hc_2-\kappa_N^{-1}\) S(z_0) + \(2N\hc_3-\kappa_N^{-1}\)h(z_0)\right] W_{1,0}(y) \\
&\ + \mu^{-{N \over 2}} \left[4\(\hc_2 S(z_0)+\hc_3 h(z_0)\)(\mu y_j+\xi_j)\right] \pa_j W_{1,0}(y) \\
&\ + \mu^{-{N+2 \over 2}} \(\frac{4}{N-2}\) \left\{\hc_1R_{ij}(z_0)(\mu y_i+\xi_i)(\mu y_j+\xi_j) + \(\hc_2 S(z_0)+\hc_3 h(z_0)\)|\mu y+\xi|^2\right\} W_{1,0}^p(y) \\
&\ + \mu^{-{N-4 \over 2}} a^{\{1\}} + \mu^{\nu_2-{N-2 \over 2}} a^{\{2\}}.
\end{aligned}
\end{equation}

Plugging \eqref{eq:error112} and \eqref{eq:error113} into \eqref{eq:error111}, we establish the desired estimate \eqref{eq:S1in}.
\end{proof}

\begin{lemma}\label{lemma:error1out}
Let $\vep_1 \in (0,1)$ be a small number. It holds that
\begin{multline}\label{eq:S1out}
\left[1-\eta_{\mu_0^{\vep_1}}(d_{g_0}(z,z_0))\right] \mcs\(u^{(1)}_{\mx}\)(z,t) \\
= \mu^{{N \over 2}-\zeta_1\vep_1} \dot{\mu} b_1(z,t)
+ \mu^{{N+2 \over 2}-\zeta_1\vep_1} b_2(z,t) + \mu^{{N-2 \over 2}+\nu_2-\zeta_1\vep_1} b_3(z,t)
\end{multline}
for $(z,t) \in M \times [t_0,\infty)$. Here,
\begin{itemize}
\item[-] $\eta_{\mu_0^{\vep_1}} = \eta_0(\mu_0^{-\vep_1}\cdot)$ is the cut-off function introduced after \eqref{eq:eta};
\item[-] $\zeta_1 > 0$ is a number depending only on $N$, and $\nu_2 = 2-\vep_0$;
\item[-] $b_1$, $b_2$, and $b_3$ are functions on $M \times [t_0,\infty)$ such that
\begin{equation}\label{eq:b}
\sum_{\ell=0}^1 \left\|\nabla^{\ell}_z b\right\|_{L^{\infty}(M \times [t_0,\infty))} +
\sup_{(z,t) \in M \times [t_0,\infty)} \mu_0^{-2}(t) [b]_{C^{\sigma_0/2}_t}(z,t) \le C\delta_0^{-\zeta}
\end{equation}
for $b = b_1$, $b_2$ and $b_3$, where $C,\, \zeta > 0$ are constants depending only on $(M,g_0)$, $N$, $h$, and $z_0$.
\end{itemize}
\end{lemma}
\begin{proof}
The proof is split into two steps.

\medskip \noindent \textsc{Step 1.} We assume that $d_{g_0}(z,z_0) > 2\delta_0$.
Because $L_{g_0,h} G(z,z_0) = 0$, we simply obtain
\begin{align}
\mcs\(u^{(1)}_{\mx}\) &= \mcs\(\ga_N^{-1}\kappa_N\alpha_N G(z,z_0) \mu^{N-2 \over 2}\) \nonumber \\
&= - \(\frac{N+2}{2}\) \(\ga_N^{-1}\kappa_N\alpha_N G(z,z_0)\)^p \mu^{N \over 2} \dot{\mu} + \frac{(N+2)\kappa_N}{4} \(\ga_N^{-1}\kappa_N\alpha_N G(z,z_0)\)^p \mu^{N+2 \over 2} \label{eq:error12} \\
&= \mu^{{N \over 2}-\zeta_1\vep_1} \dot{\mu} b_1 + \mu^{{N+2 \over 2}-\zeta_1\vep_1} b_2. \nonumber
\end{align}
By taking $b_3 = 0$, we get \eqref{eq:S1out} for this case.

\medskip \noindent \textsc{Step 2.} We assume that $\mu_0^{\vep_1} \le d_{g_0}(z,z_0) \le 2\delta_0$.
Let $\hb_1, \hb_2, \ldots$ be functions satisfying \eqref{eq:b} for $b = \hb_1, \hb_2, \ldots$. Arguing as in \eqref{eq:error12}, we find
\begin{align*}
\mcs\(u^{(1)}_{\mx}\) &= \frac{(N+2)\kappa_N}{4 \ga_N} L_{g_0,h}\(G(z,z_0)\, v^{(1)}_{\mx}\) + \mu^{{N \over 2}-\zeta_1\vep_1} \dot{\mu} \hb_1 + \mu^{{N+2 \over 2}-\zeta_1\vep_1} \hb_2.
\end{align*}
Also, applying $L_{g_0,h} G(z,z_0) = 0$ once more, we see
\[L_{g_0,h}\(G(z,z_0)\, v^{(1)}_{\mx}\) = \kappa_N \left[2 \la \nabla_{g_0} G(z,z_0), \nabla_{g_0} v^{(1)}_{\mx} \ra_{g_0} + G(z,z_0) \Delta_{g_0} v^{(1)}_{\mx}\right].\]
On the other hand, we infer from \eqref{eq:v1} and the mean value theorem that
\begin{align*}
\pa_i v^{(1)}_{\mx}(x,t) &= \alpha_N \mu^{N-2 \over 2} \left[\frac{\pa_i \(\eta_{\delta_0}(x)\)}{(\mu^2+|x-\xi|^2)^{(N-2)/2}} \left\{|x|^{N-2} - (\mu^2+|x-\xi|^2)^{N-2 \over 2}\right\} \right. \\
&\hspace{55pt} \left. + \frac{(N-2) \eta_{\delta_0}(x)}{(\mu^2+|x-\xi|^2)^{N/2}} \left\{x_i|x|^{N-4}\(\mu^2 - 2x\cdot\xi + |\xi|^2\) + |x|^{N-2}\xi \right\}\right] \\
&= \mu^{N+2 \over 2}\hb_{3i} + \mu^{{N-2 \over 2}+\nu_2} \hb_{4i}
\end{align*}
and
\[\pa_{ij} v^{(1)}_{\mx}(x,t) = \mu^{N+2 \over 2} \hb_{3ij} + \mu^{{N-2 \over 2}+\nu_2} \hb_{4ij}.\]
Consequently, by setting
\begin{align*}
\begin{medsize}
b_1
\end{medsize}
&\begin{medsize}
\displaystyle = \left[1-\eta_{\mu_0^{\vep_1}}(d_{g_0}(\cdot,z_0))\right] \hb_1,
\end{medsize} \\
\begin{medsize}
b_2
\end{medsize}
&\begin{medsize}
\displaystyle = \left[1-\eta_{\mu_0^{\vep_1}}(d_{g_0}(\cdot,z_0))\right] \left[\hb_2 + \kappa_N \left\{2g_0^{ij} \pa_iG(\cdot,z_0)
+ \(\frac{\pa_i \sqrt{|g_0|}}{\sqrt{|g_0|}} g_0^{ij} + \pa_ig_0^{ij}\) G(\cdot,z_0)\right\} \hb_{3j} + \kappa_N g_0^{ij} G(\cdot,z_0) \hb_{3ij}\right],
\end{medsize}
\end{align*}
and so on, we deduce \eqref{eq:S1out}. Note that $\text{supp}(b_3) \subset B_{g_0}(z_0,2\delta_0) \times [t_0,\infty)$.
\end{proof}

\subsection{Second approximation solutions} \label{subsec:approx2}
In this subsection, we refine the approximate solutions to reduce their errors in the ball $B_{g_0}(z,z_0)$ significantly. To this end, we introduce a linear equation
\begin{equation}\label{eq:Psi}
\Delta \Psi(\cdot,t) + pW_{1,0}^{p-1} \Psi(\cdot,t) = - \frac{4}{(N+2)\kappa_N} \mce_0[\bmu](\cdot,t) \quad \text{in } \R^N, \quad \Psi(\cdot,t) \in \dot{W}^{1,2}(\R^N)
\end{equation}
for each $t \in [t_0,\infty)$, where $\bmu$ and $\mce_0[\bmu]$ are given in \eqref{eq:mu} and \eqref{eq:E0}, respectively.

By virtue of the definition of $\bmu$, we have
\begin{equation}\label{eq:bmu}
\bmu^{-1} \dot{\bmu} = -\frac{(N+2)c_2}{4c_1} h(z_0) \bmu^2 \quad \text{on } [t_0,\infty).
\end{equation}
Hence the function $\bmu^{-2} \mce_0[\bmu]$ is independent of $t$, and a solution $\Psi$ to \eqref{eq:Psi} is decomposed into $\Psi(y,t) = \bmu^2(t) Q(y)$ on $\R^N \times [t_0,\infty)$ where $Q$ satisfies
\begin{equation}\label{eq:Q}
\Delta Q + pW_{1,0}^{p-1} Q = - \frac{4}{(N+2)\kappa_N} \(\bmu^{-2} \mce_0[\bmu]\)(y) \quad \text{in } \R^N, \quad Q \in \dot{W}^{1,2}(\R^N).
\end{equation}

From \eqref{eq:E0} and \eqref{eq:F0}, we immediately see that $(\bmu^{-2} \mce_0[\bmu])(y) = O(|y|^{2-N})$ for $|y|$ large.
In fact, as we shall see in the next lemma, a remarkable cancellation among the terms of $\mcf_0$ occurs,
so we actually have a better decay estimate for $\bmu^{-2} \mce_0[\bmu]$.
This observation turns out to be essential in constructing the refined approximate solution with sufficiently small error.
\begin{lemma}
There exists a constant $C > 0$ depending only on $(M,g_0)$, $N$, $h$, and $z_0$ such that
\begin{equation}\label{eq:E0decay}
\left|\(\bmu^{-2} \mce_0[\bmu]\)(y)\right| \le \frac{C}{1+|y|^N} \quad \text{for } (y,t) \in \R^N \times [t_0,\infty).
\end{equation}
\end{lemma}
\begin{proof}
A direct computation shows that
\begin{multline*}
\frac{1}{\alpha_N \hc_2} \left[\(2\hc_1 + 2N\hc_2-\kappa_N^{-1}\)W_{1,0}(y) + 4\hc_2 y_j\pa_j W_{1,0}(y)\right] \\
= \frac{1}{\alpha_N \hc_3} \left[\(2N\hc_3-\kappa_N^{-1}\)W_{1,0}(y) + 4\hc_3 y_j\pa_j W_{1,0}(y)\right] = \frac{4(N-2)}{(1+|y|^2)^{N \over 2}}.
\end{multline*}
In light of \eqref{eq:E0}, \eqref{eq:F0}, \eqref{eq:bmu}, and the above identity, we obtain \eqref{eq:E0decay}.
\end{proof}

\begin{lemma}
Equation \eqref{eq:Q} has a solution.
\end{lemma}
\begin{proof}
It is well-known that the space of all bounded solutions to
\[L_0[\Psi] := W_{1,0}^{1-p} \(\Delta \Psi + pW_{1,0}^{p-1} \Psi\) = 0 \quad \text{in } \R^N\]
is spanned by the function $Z_{N+1}$ in \eqref{eq:ZN+1} and
\begin{equation}\label{eq:Zi}
Z_i(y) := \pa_i W_{1,0}(y) = -(N-2)\alpha_N \frac{y_i}{(1+|y|^2)^{N \over 2}} \quad \text{for } y \in \R^N \text{ and } i = 1, \ldots, N.
\end{equation}
This fact, \eqref{eq:E0decay}, the dimensional assumption $N \ge 5$, and the Fredholm alternative imply that \eqref{eq:Q} is solvable if and only if
\begin{equation}\label{eq:ortho}
\int_{\R^N} \(\bmu^{-2} \mce_0[\bmu]\)(y) Z_n(y) dy = 0
\end{equation}
for $n = 1, \ldots, N+1$.

By parity, \eqref{eq:ortho} holds automatically for $n = 1, \ldots, N$. Besides, setting
\[c_3 = \int_{\R^N} (y \cdot \nabla W_{1,0}(y)) Z_{N+1}(y)dy \quad \text{and} \quad c_4 = \int_{\R^N} |y|^2 \(W_{1,0}^p Z_{N+1}\)(y) dy,\]
we deduce from \eqref{eq:E0}, \eqref{eq:F0}, \eqref{eq:c1c2}, \eqref{eq:hc}, and \eqref{eq:bmu} that
\begin{equation}\label{eq:orthoN+1}
\begin{aligned}
&\ \begin{medsize}
\displaystyle \int_{\R^N} \mce_0[\bmu](y,t) Z_{N+1}(y) dy
\end{medsize} \\
&\begin{medsize}
\displaystyle = c_1 \bmu^{-1} \dot{\bmu} + \frac{(N+2)\kappa_N}{4} \underbrace{\left[-(2N\hc_3-\kappa_N^{-1})c_2 + 4\hc_3c_3 + \frac{4}{N-2}\hc_3c_4\right]}_{=\kappa_N^{-1}c_2} h(z_0) \bmu^2
\end{medsize} \\
&\ \begin{medsize}
\displaystyle + \frac{(N+2)\kappa_N}{4} \underbrace{\left[-\(2\hc_1 + 2N\hc_2 - \kappa_N^{-1}\)c_2
+ \(\frac{4\hc_1}{N} + 4\hc_2 - \frac{1}{3N}\)c_3 + \frac{4}{N-2} \(\frac{\hc_1}{N} + \hc_2\)c_4 \right]}_{=0} S(z_0) \bmu^2
\end{medsize} \\
&\begin{medsize}
\displaystyle = c_1 \bmu^{-1} \dot{\bmu} + \frac{(N+2)c_2}{4} h(z_0) \bmu^2 = 0.
\end{medsize}
\end{aligned}
\end{equation}
Thus \eqref{eq:Q} is solvable.
\end{proof}

\begin{lemma}\label{lemma:Qdecay}
There is a constant $C > 0$ depending on $Q$ such that
\begin{equation}\label{eq:Qdecay}
\left|\nabla^{\ell}_y Q(y)\right| \le \frac{C\log (2+|y|)}{|y|^{N-2+\ell}} \quad \text{for } y \in \R^N \text{ and } \ell = 0, 1, 2, 3.
\end{equation}
\end{lemma}
\begin{proof}
By employing the rescaling argument with the condition $Q \in \dot{W}^{1,2}(\R^N)$, we obtain
\[|Q(y)| \le \frac{C}{1+|y|^{N-2 \over 2}} \quad \text{for } y \in \R^N.\]
Then, having \eqref{eq:E0decay} in hand, we apply the comparison principle to \eqref{eq:Q} repeatedly, which produces
\begin{equation}\label{eq:Qdecay2}
|Q(y)| \le \frac{C}{1+|y|^{N-2-\vep}} \quad \text{for } y \in \R^N
\end{equation}
where $\vep > 0$ is any small number. By \eqref{eq:Q}, \eqref{eq:E0decay}, and \eqref{eq:Qdecay2},
\[\Delta Q(y) = - \frac{4}{(N+2)\kappa_N} \(\bmu^{-2} \mce_0[\bmu]\)(y) - p W_{1,0}^{p-1} Q = O\(\frac{1}{|y|^N}\) \quad \text{for } |y| \ge 1.\]
Consequently, from the Green's representation formula for $Q$, we get
\[\left|Q(y)\right| \le C \int_{\R^N} \frac{\left|\Delta Q(\ty)\right|}{|y-\ty|^{N-2}} d\ty
\le C \int_{\R^N} \frac{1}{|y-\ty|^{N-2}} \frac{d\ty}{1+|\ty|^N} \le \frac{C\log(2+|y|)}{|y|^{N-2}} \quad \text{for } |y| \ge 1\]
where we estimate the second integral by decomposing the domain $\R^N$ into
\[\overline{B^N(0,R_0)} \cup \overline{B^N(y,R_0)} \cup \(\R^N \setminus \(B^N(0,R_0) \cup B^N(y,R_0)\)\) \quad \text{for } R_0 := \frac{|y|}{2}.\]
Therefore, \eqref{eq:Qdecay} for $\ell = 0$ holds. 
The gradient estimate for $Q$, that is, \eqref{eq:Qdecay} for $\ell = 1, 2, 3$ follows from elliptic regularity.
\end{proof}
\noindent Let $\Psi_0 = \bmu^2 Q_0$ be the unique solution to \eqref{eq:Psi} such that
\[\int_{\R^N} \(Q_0Z_n\)(y)dy = 0 \quad \text{for } n = 1, \ldots, N+1.\]

\medskip
We now define the second (or, refined) approximate solution
\begin{equation}\label{eq:u2}
u^{(2)}_{\mx}(z,t) = \ga_N^{-1} \kappa_N G(z,z_0) v^{(2)}_{\mx}(z,t) \quad \text{on } M \times [t_0,\infty)
\end{equation}
where
\begin{equation}\label{eq:v2}
\begin{aligned}
&\ v^{(2)}_{\mx}(z,t) \\
&= \begin{cases}
\displaystyle |x|^{N-2} \mu^{-{N-2 \over 2}} \left[W_{1,0}(y) + \Psi_0(y,t)\right] &\text{if } d_{g_0}(z,z_0) \le \delta_0, \\
\begin{aligned}
& \mu^{N-2 \over 2} \left[\alpha_N + \eta_{\delta_0}(x) \left\{\alpha_N \(\frac{|x|^{N-2}}{(\mu^2+|x-\xi|^2)^{(N-2)/2}} - 1\) \right. \right. \\
&\hspace{90pt} \left.\left. + \mu^{-(N-2)} |x|^{N-2} \Psi_0(y,t)\)\right]
\end{aligned}
&\text{if } \delta_0 < d_{g_0}(z,z_0) \le 2\delta_0, \\
\displaystyle \mu^{N-2 \over 2} \alpha_N &\text{if } d_{g_0}(z,z_0) > 2\delta_0.
\end{cases}
\end{aligned}
\end{equation}
Here, $z = \exp_{z_0}(x) = \exp_{z_0}(\mu y+\xi)$ for $z \in B_{g_0}(z_0,2\delta_0)$.
By taking $t_0$ large, we see from \eqref{eq:Qdecay} that $u^{(2)}_{\mx}$ is of class $C^{2+\sigma_0,1+\sigma_0/2}(M \times [t_0,\infty))$ and positive on $M \times [t_0,\infty)$.
The next lemmas measure its error.
\begin{lemma}
It holds that
\begin{equation}\label{eq:S2in}
\mu^{N+2 \over 2} \mcs\(u^{(2)}_{\mx}\)(y,t) = \mce_2[\mx](y,t)
\end{equation}
for $y \in B^N(-\mu^{-1}\xi,\mu^{-1}\delta_0)$ and $t \in [t_0,\infty)$. Here,
\begin{align}
\mce_2[\mx](y,t) &:= \bmu^{-1} \(\dot{\lambda} - \bmu^{-1} \dot{\bmu} \lambda\) \(pW_{1,0}^{p-1}Z_{N+1}\)(y)
+ \mu^{-1} \dot{\xi} \cdot \(p W_{1,0}^{p-1} \nabla W_{1,0}\)(y) + 2\bmu\lambda \mcf_0(y) \nonumber \\
&\ + \mcf_1[\mx](y,t) + \mu^3 a^{\{1\}} + \mu^{2\nu_1} a^{\{2\}} + \mu\dot{\mu} a^{\{0\}}
+ \mu^{\nu_1-2} \dot{\lambda}\, a^{\{-2\}} + \mu\dot{\xi} \cdot \mba^{\{-1\}} \label{eq:E2}
\end{align}
where
\begin{itemize}
\item[-] $\mcf_0$, $\mcf_1$, and $P$ are the functions in \eqref{eq:F0}, \eqref{eq:F1}, and \eqref{eq:P_20}, respectively, and $\nu_1$ and $\nu_2$ are numbers in \eqref{eq:lx};
\item[-] $a^{\{\ga\}} \in \R$ and $\mba^{\{\ga\}} \in \R^N$ are $C^1$ functions of $(y,\mu,\mu^{-1}\xi,\mu^{-1}\lambda)$ satisfying \eqref{eq:adecay}.
\end{itemize}
\end{lemma}
\begin{proof}
We see from \eqref{eq:u1}--\eqref{eq:v1}, \eqref{eq:u2}--\eqref{eq:v2}, \eqref{eq:S1in}, and \eqref{eq:Psi} that
\begin{equation}\label{eq:S21}
\begin{aligned}
\mu^{N+2 \over 2} \mcs\(u^{(2)}_{\mx}\) &= \mu^{N+2 \over 2} \mcs\(u^{(1)}_{\mx} + (1+P) \mu^{-{N-2 \over 2}} \Psi_0\) \\
&= \frac{(N+2)\kappa_N}{4} \(\Delta_{g_0(x)} \Psi_0 + p W_{1,0}^{p-1} \Psi_0\) + \mu^{N+2 \over 2} \mcs\(u^{(1)}_{\mx}\) + \mca(\Psi_0)\\
&= \frac{(N+2)\kappa_N}{4} \(\Delta_{g_0(x)}-\Delta\) \Psi_0 + (\mce_0[\mu]-\mce_0[\bmu]) + \mce_1[\mx] + \mca(\Psi_0)
\end{aligned}
\end{equation}
in the $(y,t)$-variable. Here,
\begin{align*}
\mca(\Psi_0) &:= \frac{(N+2)\kappa_N}{4} \left[\(W_{1,0}+\Psi_0\)^p - W_{1,0}^p - pW_{1,0}^{p-1}\Psi_0\right] \\
&\ + p \mu^{-1} \dot{\mu} \left[\(W_{1,0}+\Psi_0\)^{p-1} \left\{y \cdot \nabla(W_{1,0}+\Psi_0) + \frac{N-2}{2} \(W_{1,0}+\Psi_0\) \right\} - W_{1,0}^{p-1}Z_{N+1} \right] \\
&\ + p \mu^{-1} \dot{\xi} \cdot \left[\(W_{1,0}+\Psi_0\)^{p-1} \nabla(W_{1,0}+\Psi_0) - W_{1,0}^{p-1}\nabla W_{1,0} \right] + \mu^4 \log(2+|y|) a^{\{2\}}.
\end{align*}
By Taylor's theorem and \eqref{eq:Qdecay},
\begin{equation}\label{eq:S22}
\mca(\Psi_0) = \log(2+|y|) \left[\mu^4 a^{\{2\}} + \mu \dot{\mu} a^{\{-2\}} + \mu \dot{\xi} \cdot \mba^{\{-3\}}\right].
\end{equation}
Also, in view of \eqref{eq:lx}, \eqref{eq:metric}, \eqref{eq:Qdecay}, it holds that
\begin{equation}\label{eq:S23}
\(\Delta_{g_0(x)}-\Delta\) \Psi_0 = \mu^4 \log(2+|y|) a^{\{2\}}
\end{equation}
and
\begin{multline}\label{eq:S24}
\mce_0[\mu] - \mce_0[\bmu] = \bmu^{-1} \(\dot{\lambda} - \bmu^{-1} \dot{\bmu} \lambda\) \(pW_{1,0}^{p-1}Z_{N+1}\)(y) + 2\bmu\lambda \mcf_0(y) \\
+ \mu^{2\nu_1} a^{\{2\}} + \mu^{\nu_1-2} \dot{\lambda}\, a^{\{-2\}}.
\end{multline}

Inserting \eqref{eq:S22}--\eqref{eq:S24} and \eqref{eq:E1} into \eqref{eq:S21}, and then arranging the resulting terms, we deduce the desired equality \eqref{eq:S2in}.
\end{proof}

\begin{lemma}\label{lemma:error2out}
It holds that
\begin{multline}\label{eq:S2out}
\left[1-\eta_{\mu_0^{\vep_1}}(d_{g_0}(z,z_0))\right] \mcs\(u^{(2)}_{\mx}\)(z,t) \\
= \mu^{{N \over 2}-\zeta_2\vep_1} \dot{\mu} b_4(z,t) + \mu^{{N+2 \over 2}-\zeta_2\vep_1} b_5(z,t) + \mu^{{N-2 \over 2}+\nu_2-\zeta_2\vep_1} b_6(z,t)
\end{multline}
for $(z,t) \in M \times [t_0,\infty)$. Here,
\begin{itemize}
\item[-] $\vep_1 \in (0,1)$ is the small number in Lemma \ref{lemma:error1out},
    and $\eta_{\mu_0^{\vep_1}}$ is the cut-off function introduced after \eqref{eq:eta};
\item[-] $\zeta_2 > 0$ is a number depending only on $N$, and $\nu_2 = 2-\vep_0$;
\item[-] $b_4$, $b_5$, and $b_6$ are functions on $M \times [t_0,\infty)$ such that \eqref{eq:b} holds for $b = b_4$, $b_5$, and $b_6$.
\end{itemize}
\end{lemma}
\begin{proof}
If $d_{g_0}(z,z_0) > 2\delta_0$, then by \eqref{eq:u1}, \eqref{eq:u2}, and \eqref{eq:error12}, we have
\[\mcs\(u^{(2)}_{\mx}\) = \mcs\(u^{(1)}_{\mx}\) = \mu^{{N \over 2}-\zeta_1\vep_1} \dot{\mu} b_1 + \mu^{{N+2 \over 2}-\zeta_1\vep_1} b_2.\]

Suppose that $\mu_0^{\vep_1} \le d_{g_0}(z,z_0) \le 2\delta_0$. Setting
\[\tu^{(2)}_{\mx} = \ga_N^{-1} \kappa_N G(z,z_0) \cdot \mu^{-{N-2 \over 2}} \bmu^2 \eta_{\delta_0}(x) |x|^{N-2} Q_0(y),\]
we write
\begin{multline}\label{eq:S26}
\mcs\(u^{(2)}_{\mx}\) = \mcs\(u^{(1)}_{\mx}\) + \frac{N+2}{4} L_{g_0,h} \tu^{(2)}_{\mx} \\
- \left\{\pa_t\(u^{(1)}_{\mx}+\tu^{(2)}_{\mx}\)^p - \pa_t \(u^{(1)}_{\mx}\)^p\right\} + \frac{(N+2)\kappa_N}{4} \left\{\(u^{(1)}_{\mx}+\tu^{(2)}_{\mx}\)^p - \(u^{(1)}_{\mx}\)^p\right\}.
\end{multline}
Let $\hb_5$, $\hb_6$, and $\hb_7$ be functions satisfying \eqref{eq:b}. Applying \eqref{eq:Qdecay}, we argue as in Step 2 of the proof of Lemma \ref{lemma:error1out}. Then we obtain
\begin{equation}\label{eq:S27}
L_{g_0,h} \tu^{(2)}_{\mx} = \mu^{{N+2 \over 2}-\vep_1\zeta_2} \hb_5.
\end{equation}
Also, we compute
\begin{multline}\label{eq:S28}
\pa_t\(u^{(1)}_{\mx}+\tu^{(2)}_{\mx}\)^p - \pa_t \(u^{(1)}_{\mx}\)^p \\
= p \left[\left\{\(u^{(1)}_{\mx}+\tu^{(2)}_{\mx}\)^{p-1} - \(u^{(1)}_{\mx}\)^{p-1} \right\} \pa_t u^{(1)}_{\mx} + \(u^{(1)}_{\mx}+\tu^{(2)}_{\mx}\)^{p-1} \pa_t \tu^{(2)}_{\mx}\right]
= \mu^{{N+4 \over 2}-\vep_1\zeta_2} \dot{\mu} \hb_6
\end{multline}
and
\begin{equation}\label{eq:S29}
\(u^{(1)}_{\mx}+\tu^{(2)}_{\mx}\)^p - \(u^{(1)}_{\mx}\)^p = \mu^{{N+6 \over 2}-\vep_1\zeta_2} \hb_7.
\end{equation}
By substituting \eqref{eq:S27}--\eqref{eq:S29}, and \eqref{eq:S1out} into \eqref{eq:S26}, and then defining $b_4$, $b_5$, and $b_6$ suitably, we obtain \eqref{eq:S2out}.
\end{proof}

\section{Inner-outer gluing procedure} \label{sec:inout}
For the sake of brevity, we write $u_{\mx} = u_{\mx}^{(2)}$ in the sequel.

In the rest of the paper, we construct a remainder term $\psi_{\mx}$ such that $u = u_{\mx} + \psi_{\mx}$ is a solution to \eqref{eq:Yamabefp2}.
To this end, we apply the inner-outer gluing procedure as in \cite{CdPM}.
It amounts to decomposing $\psi_{\mx}$ into two parts
\begin{equation}\label{eq:psimxo}
\psi_{\mx} = \psi_{\mx}^{\out} + \eta_{\mu_0^{\vep_1}}(d_{g_0}(z,z_0))\, \hps_{\mx}^{\, \tin}\(\exp_{z_0}^{-1}(z),t\) \quad \text{on } M \times [t_0,\infty)
\end{equation}
and determining $\psi_{\mx}^{\out}$ and $\hps_{\mx}^{\, \tin}$ by solving so-called outer and inner problems.

\medskip
\noindent \textsc{Outer problem}: Let $\psi_{\mx}^{\out}$ be a function on $M \times [t_0,\infty)$ solving
\begin{multline}\label{eq:outer}
p u_{\mx}^{p-1} \(\psi_{\mx}^{\out}\)_t = \frac{N+2}{4} L_{g_0,h} \psi_{\mx}^{\out} + \mcw_{\mx} \psi_{\mx}^{\out} \\
+ \(1-\eta_{\mu_0^{\vep_1}}\)\mcs(u_{\mx}) + \mcj_1\left[\hps_{\mx}^{\, \tin}\right] + \mcj_2\left[\psi_{\mx}^{\out}, \hps_{\mx}^{\, \tin}\right].
\end{multline}
Here, $L_{g_0,h}$ is the differential operator defined in \eqref{eq:Lgh},
\[\mcw_{\mx} := \frac{(N+2)\kappa_Np}{4} \left[u_{\mx}^{p-1} - \eta_{\mu_0^{\vep_1}} \(u_{\mx}^{(1)}\)^{p-1}\right],\]
\begin{equation}\label{eq:mcj1}
\begin{aligned}
\mcj_1\left[\hps_{\mx}^{\, \tin}\right]
&:= \frac{(N+2)\kappa_Np}{4} \left[u_{\mx}^{p-1} - \(u_{\mx}^{(1)}\)^{p-1}\right] \eta_{\mu_0^{\vep_1}} \hps_{\mx}^{\, \tin}
- p \left[u_{\mx}^{p-1} - \(u_{\mx}^{(1)}\)^{p-1}\right] \eta_{\mu_0^{\vep_1}} \(\hps_{\mx}^{\, \tin}\)_t \\
&\ + \frac{(N+2)\kappa_N}{4} \left[\(\Delta_{g_0}\eta_{\mu_0^{\vep_1}}\) \hps_{\mx}^{\, \tin}
+ 2 \la \nabla_{g_0} \eta_{\mu_0^{\vep_1}}, \nabla_{g_0} \hps_{\mx}^{\, \tin} \ra_{g_0} \right]
\end{aligned}
\end{equation}
and
\begin{equation}\label{eq:mcj2}
\begin{aligned}
\mcj_2\left[\psi_{\mx}^{\out}, \hps_{\mx}^{\, \tin}\right]
&:= \frac{(N+2)\kappa_N}{4} \left[\(u_{\mx} + \psi_{\mx}\)^p - u_{\mx}^p - pu_{\mx}^{p-1} \psi_{\mx}\right] \\
&\ - p \left[\(u_{\mx} + \psi_{\mx}\)^{p-1} - u_{\mx}^{p-1}\right] \left[\(1-\eta_{\mu_0^{\vep_1}/2}\) (u_{\mx})_t + \(\psi_{\mx}\)_t\right].
\end{aligned}
\end{equation}

\medskip
\noindent \textsc{Inner problem}: As before, we identify points $z \in B_{g_0}(z_0,\delta_0)$ and $x \in B^N(0,\delta_0)$ via $g_0$-normal coordinates at $z_0 \in M$.
Let $\hps_{\mx}^{\, \tin}$ be a function on $B^N(0,2\mu_0^{\vep_1}) \times [t_0,\infty)$ satisfying
\begin{align*}
p \(u_{\mx}^{(1)}\)^{p-1} \(\hps_{\mx}^{\, \tin}\)_t &= \frac{N+2}{4} \left[L_{g_0,h} \hps_{\mx}^{\, \tin}
+ \kappa_N p\(u_{\mx}^{(1)}\)^{p-1} \(\hps_{\mx}^{\, \tin} + \psi_{\mx}^{\out}\) \right] \\
&\ + \mu^{-{N+2 \over 2}} \mce_2[\mx](\mu^{-1}(x-\xi),t) \\
&\ - p \left[\(u_{\mx} + \hps_{\mx}^{\, \tin} + \psi_{\mx}^{\out}\)^{p-1} - u_{\mx}^{p-1}\right] \eta_{\mu_0^{\vep_1}/2} \(u_{\mx}\)_t.
\end{align*}
We write $x = \bmu\by+\xi$ and define $\psi_{\mx}^{\tin}$ by the relation
\begin{equation}\label{eq:psi0}
\begin{aligned}
\hps_{\mx}^{\, \tin}(x,t) &= \ga_N^{-1} \kappa_N |x|^{N-2} G(x,z_0) \cdot \bmu^{-{N-2 \over 2}} \psi_{\mx}^{\tin}\(\bmu^{-1}(x-\xi), t\) \\
&= (1+P(x)) \bmu^{-{N-2 \over 2}} \psi_{\mx}^{\tin}\(\by, t\)
\end{aligned}
\end{equation}
for $(\by,t) \in B^N(-\bmu^{-1}\xi, 2\bmu^{-1}\mu_0^{\vep_1})$ where $P$ is the function in \eqref{eq:P_20}. Then $\psi_{\mx}^{\tin}$ solves
\begin{equation}\label{eq:inner}
\begin{aligned}
p W_{1,0}^{p-1} \(\psi_{\mx}^{\tin}\)_t &= \frac{(N+2)\kappa_N}{4} \(\Delta \psi_{\mx}^{\tin} + p W_{1,0}^{p-1} \psi_{\mx}^{\tin}\) \\
&\ + \({\bmu \over \mu}\)^{N-2 \over 2} (1+P)^{-p}(x) \mce_2[\mx](y,t) + \mck_1\left[\psi_{\mx}^{\tin}\right] + \mck_2\left[\psi_{\mx}^{\tin}, \psi_{\mx}^{\out}\right]
\end{aligned}
\end{equation}
in the $(\by,t)$-variable, where $y = \mu^{-1}\bmu\by$,
\begin{equation}\label{eq:mck1}
\begin{aligned}
\begin{medsize}
\displaystyle \mck_1\left[\psi_{\mx}^{\tin}\right]
\end{medsize}
&\begin{medsize}
\displaystyle := \frac{(N+2)\kappa_Np}{4} \left[W_{1,0}^{p-1}(y) - W_{1,0}^{p-1}(\by)\right] \psi_{\mx}^{\tin}
- p \left[W_{1,0}^{p-1}(y) - W_{1,0}^{p-1}(\by)\right] \(\psi_{\mx}^{\tin}\)_t
\end{medsize} \\
&\, \begin{medsize}
\displaystyle + p W_{1,0}^{p-1}(y) \left[\bmu^{-1} \dot{\bmu} \(\by \cdot \nabla \psi_{\mx}^{\tin} + \frac{N-2}{2} \psi_{\mx}^{\tin}\) + \bmu^{-1} \dot{\xi} \cdot \nabla \psi_{\mx}^{\tin}\right]
\end{medsize} \\
&\, \begin{medsize}
\displaystyle + \frac{(N+2)\kappa_N}{4} \left[\({\mu \over \bmu}\)^2 (1+P)^{1-p}(x) \Delta_{g_0(x)} \psi_{\mx}^{\tin} - \Delta \psi_{\mx}^{\tin}\right]
- \frac{N+2}{4} \mu^2 (1+P(x)) (S+h)(x) \psi_{\mx}^{\tin}
\end{medsize} \\
&\, \begin{medsize}
\displaystyle + \frac{(N+2)\kappa_N}{4} \({\mu \over \bmu}\)^2 (1+P)^{-p}(x)
\left[2\bmu g_0^{ij}(x) (\pa_{x_i}P)(x) \pa_j\psi_{\mx}^{\tin} + \bmu^2 (\Delta_{g_0}P)(x) \psi_{\mx}^{\tin}\right]
\end{medsize}
\end{aligned}
\end{equation}
and
\begin{equation}\label{eq:mck2}
\begin{aligned}
\begin{medsize}
\displaystyle \mck_2\left[\psi_{\mx}^{\tin}, \psi_{\mx}^{\out}\right]
\end{medsize}
&\begin{medsize}
\displaystyle := \frac{(N+2)\kappa_Np}{4} \({\bmu \over \mu}\)^{N-2 \over 2} (1+P)^{-1}(x) W_{1,0}^{p-1}(y) \mu^{N-2 \over 2} \psi_{\mx}^{\out}(x,t)
\end{medsize} \\
&\, \begin{medsize}
\displaystyle - p \({\bmu \over \mu}\)^{N-2 \over 2} \left[\left\{\(W_{1,0}+\bmu^2 Q_0\)(y)
+ \({\mu \over \bmu}\)^{N-2 \over 2} \psi_{\mx}^{\tin} + (1+P)^{-1}(x) \mu^{N-2 \over 2} \psi_{\mx}^{\out}(x,t) \right\}^{p-1} \right.
\end{medsize} \\
&\hspace{70pt} \begin{medsize}
\displaystyle \left. - \(W_{1,0}+\bmu^2 Q_0\)^{p-1}(y) \right] \eta_{\mu_0^{\vep_1}/2}(x) \mu^{N-2 \over 2} \left\{\mu^{-{N-2 \over 2}} \(W_{1,0}+\bmu^2 Q_0\)(y) \right\}_t.
\end{medsize}
\end{aligned}
\end{equation}

\medskip
If $\psi_{\mx}^{\out}$ and $\psi_{\mx}^{\tin}$ solve \eqref{eq:outer} and \eqref{eq:inner}, respectively,
then $u = u_{\mx} + \psi_{\mx}$ will satisfy \eqref{eq:Yamabefp2} provided it is positive on $M \times [t_0,\infty)$.
In Sections \ref{sec:outer} and \ref{sec:inner}, we will look for solutions to \eqref{eq:outer} and \eqref{eq:inner}
with appropriate choices of parameters $(\mx)$ and initial conditions $u_0$, and verify that $u$ is indeed positive.

\section{Outer problem} \label{sec:outer}
This section is devoted to the analysis of the outer problem \eqref{eq:outer}.

In Subsections \ref{subsec:outer1} and \ref{subsec:outer2}, we develop existence theory and a priori estimates on a solution to an associated inhomogeneous problem \eqref{eq:inhom-o}.
The main technical point is to control the degenerate factor $u_{\mx}^{p-1}$.

In Subsections \ref{subsec:outer3} and \ref{subsec:outer4}, we apply the results for \eqref{eq:inhom-o}
and the contraction mapping theorem to derive the unique solvability of \eqref{eq:outer} as well as several a priori estimates on the solution.

\subsection{Inhomogeneous problem associated to \eqref{eq:outer}: Weighted $H^2$ estimate} \label{subsec:outer1}
In this subsection, we will prove the unique existence of a solution to the inhomogeneous problem
\begin{equation}\label{eq:inhom-o}
\begin{cases}
\displaystyle p u_{\mx}^{p-1} \psi_t = \frac{N+2}{4} L_{g_0,h} \psi + \mcw_{\mx} \psi + u_{\mx}^{p-1} f &\text{on } M \times (t_0,\infty), \\
\psi(\cdot,t_0) = \psi_0 &\text{on } M
\end{cases}
\end{equation}
in a weighted $H^2$ space by deriving a priori estimate for solutions. Here and after, we mean by a solution to \eqref{eq:inhom-o} and other related equations a function that satisfies them in a weak sense.

\medskip
We first establish the following version of a priori weighted $H^2$ estimate. Refer to Definition \ref{defn:norm1} for the definition of the norms.
\begin{lemma}\label{lemma:outer1}
Suppose that $Y_h(M,g_0) > 0$, $t_0 < s_0-1$, $\|f\|_{L^2_{t_0,s_0}} + \|\psi_0\|_{\mch^1} < \infty$ and \eqref{eq:mu0}--\eqref{eq:lx} hold.
Then there is a constant $C > 0$ depending only on $(M,g_0)$, $N$, $h$, and $z_0$ such that every solution $\psi$ to
\begin{equation}\label{eq:inhom-o2}
\begin{cases}
\displaystyle p u_{\mx}^{p-1} \psi_t = \frac{N+2}{4} L_{g_0,h} \psi + \mcw_{\mx} \psi + u_{\mx}^{p-1} f &\text{on } M \times (t_0,s_0), \\
\psi(\cdot,t_0) = \psi_0 &\text{on } M
\end{cases}
\end{equation}
satisfies
\begin{equation}\label{eq:outer1}
\|\psi\|_{H^2_{t_0,s_0}} \le C\(\|\psi\|_{L^2_{t_0,s_0}} + \|f\|_{L^2_{t_0,s_0}} + \|\psi_0\|_{\mch^1}\)
\end{equation}
provided $t_0 > 0$ is large enough.
\end{lemma}
\begin{proof}
The proof is inspired by that of \cite[Lemma 3.2]{DdPS}. We will divide it into four steps.

Throughout the proof, we assume that $C$ depends only on $(M,g_0)$, $N$, $h$, and $z_0$, and in particular, is independent of $s_0$.

\medskip \noindent \textsc{Step 1.} We claim that
\begin{equation}\label{eq:ArBe}
u_{\mx}^{-1} |(u_{\mx})_t| \le C \mu_0^2.
\end{equation}

If $z \in B_{g_0}(z_0,\delta_0)$, we infer from \eqref{eq:u2}--\eqref{eq:v2}, \eqref{eq:Qdecay}, and \eqref{eq:mu0}--\eqref{eq:lx} that
\begin{align*}
u_{\mx}^{-1} (u_{\mx})_t &= \frac{N-2}{2} \mu^{-1} \dot{\mu} + \left[ \frac{1}{(1+|y|^2)^{N-2 \over 2}} + \bmu^2 Q_0(y)\right]^{-1} \\
&\hspace{75pt} \times \left[\frac{(N-2) \mu^{-1} (\dot{\xi} \cdot y-\dot{\mu})}{(1+|y|^2)^{N \over 2}}
+ O(\mu_0^4) Q_0(y) + O(\mu_0) \(\dot{\mu}y-\dot{\xi}\) \cdot \nabla Q_0(y)\right] \\
&= O\(\mu^{-1} \dot{\mu} + \frac{\mu^{-1} \dot{\xi}}{1+|y|}\) = O(\mu_0^2)
\end{align*}
in the $(y,t)$-variable.

If $z \in B_{g_0}(z_0,2\delta_0) \setminus \overline{B_{g_0}(z_0,\delta_0)}$, then
\[u_{\mx}^{-1} (u_{\mx})_t = \frac{N-2}{2} \mu^{-1} \dot{\mu} + O\(\mu_0^2 \(\mu_0^2 \log\mu_0 + \delta_0^{-1} \mu_0^{\nu_2}\)\) = O(\mu_0^2).\]

If $z \in M \setminus \overline{B_{g_0}(z_0,2\delta_0)}$, then
\[u_{\mx}^{-1} (u_{\mx})_t = \frac{N-2}{2} \mu^{-1} \dot{\mu} = O(\mu_0^2).\]

Consequently, the claim follows.

\medskip \noindent \textsc{Step 2.} Fixing any $t \in [t_0,s_0]$, we multiply \eqref{eq:inhom-o} by $\psi$ and integrate the resultant equality over $M$. Then we obtain
\begin{multline}\label{eq:outer11}
\frac{p}{2}\, \pa_t \int_M \psi^2 u_{\mx}^{p-1} dv_{g_0} - {N+2 \over 4} \int_M \psi L_{g_0,h} \psi \, dv_{g_0} \\
= \int_M \mcw_{\mx} \psi^2 dv_{g_0} + \int_M f \psi u_{\mx}^{p-1} dv_{g_0} + \frac{p(p-1)}{2} \int_M \left[u_{\mx}^{-1} (u_{\mx})_t\right] \psi^2 u_{\mx}^{p-1} dv_{g_0}.
\end{multline}
On the other hand, H\"older's inequality yields
\[\int_M \psi^2 u_{\mx}^{p-1} dv_{g_0} \le C \(\int_M u_{\mx}^{p+1} dv_{g_0}\)^{p-1 \over p+1} \(\int_M |\psi|^{p+1} dv_{g_0}\)^{2 \over p+1} \le C \(\int_M |\psi|^{p+1} dv_{g_0}\)^{2 \over p+1}.\]
Hence, together with the condition $Y_h(M,g_0) > 0$ and the Sobolev inequality, we deduce \begin{equation}\label{eq:outer12}
\int_M \psi^2 u_{\mx}^{p-1} dv_{g_0} 
\le C \int_M (|\nabla_{g_0} \psi|_{g_0}^2 + \psi^2)\, dv_{g_0} \le - C \int_M \psi L_{g_0,h} \psi \, dv_{g_0}.
\end{equation}
Owing to \eqref{eq:ArBe}, \eqref{eq:outer12}, the bound $|\mcw_{\mx}| \le Cu_{\mx}^{p-1}$ on $M \times [t_0,s_0]$, and H\"older's inequality, we see from \eqref{eq:outer11} that
\begin{equation}\label{eq:outer13}
\begin{aligned}
&\ \pa_t \int_M \psi^2 u_{\mx}^{p-1} dv_{g_0} + \int_M (|\nabla_{g_0} \psi|_{g_0}^2 + \psi^2)\, dv_{g_0} \\
&\le C\(\int_M \mcw_{\mx} \psi^2 dv_{g_0} + \int_M f \psi u_{\mx}^{p-1} dv_{g_0} + \mu_0^2(t) \int_M \psi^2 u_{\mx}^{p-1} dv_{g_0}\) \\
&\le 
C\(\|\psi(\cdot,t)\|_{L^2(M)}^2 + \|f(\cdot,t)\|_{L^2(M)}^2\).
\end{aligned}
\end{equation}
Given any $\tau \in (t_0,s_0-1]$, we set a function $\chi(t) = t-\tau \in [0,1]$ for $t \in [\tau,\tau+1]$. For such $t$,
\[\pa_t \(\chi(t) \int_M \psi^2 u_{\mx}^{p-1} dv_{g_0}\) + \chi(t) \int_M (|\nabla_{g_0} \psi|_{g_0}^2 + \psi^2)\, dv_{g_0}
\le C\(\|\psi(\cdot,t)\|_{L^2(M)}^2 + \|f(\cdot,t)\|_{L^2(M)}^2\)\]
by \eqref{eq:outer13}. Integrating it over $t \in [\tau,\tau+1]$, we arrive at
\begin{multline}\label{eq:outer14}
\int_M \(\psi^2 u_{\mx}^{p-1}\) (z,\tau+1)\, dv_{g_0} + \iint_{M_{\tau}} \chi(t) (|\nabla_{g_0} \psi|_{g_0}^2 + \psi^2)\, dv_{g_0} dt \\
\le C\(\|\psi\|_{L^2(M_{\tau})}^2 + \|f\|_{L^2(M_{\tau})}^2\).
\end{multline}
From \eqref{eq:outer13}, we also derive
\begin{multline}\label{eq:outer15}
\int_M \(\psi^2 u_{\mx}^{p-1}\) (z,t_0+1)\, dv_{g_0} + \iint_{M_{t_0}} (|\nabla_{g_0} \psi|_{g_0}^2 + \psi^2)\, dv_{g_0} dt \\
\le C\(\|\psi\|_{L^2(M_{t_0})}^2 + \|f\|_{L^2(M_{t_0})}^2 + \|\psi_0\|_{\mch^1}^2\).
\end{multline}
To get this inequality, we use the initial datum $\psi_0$ in \eqref{eq:inhom-o} instead of introducing the cut-off function $\chi(t)$, and \eqref{eq:outer12}.

\medskip \noindent \textsc{Step 3.} Fixing $t \in [t_0,s_0]$ again, we multiply \eqref{eq:inhom-o} by $\psi_t$, integrate the result over $M$,
and then apply H\"older's inequality and $|\pa_t \mcw_{\mx}| \le C\mu_0^2 u_{\mx}^{p-1}$ on $M \times [t_0,s_0]$. Then we observe
\begin{multline*}
\pa_t \left[- {N+2 \over 4} \int_M \psi L_{g_0,h} \psi\, dv_{g_0} - \int_M \mcw_{\mx} \psi^2 dv_{g_0}\right] + \|\psi_t\|_{L^2(M)}^2 \\
\le C\(\int_M |\pa_t \mcw_{\mx}| \psi^2 dv_{g_0} + \|f\|_{L^2(M)}^2\) \le C\(\mu_0^2(t) \|\psi\|_{L^2(M)}^2 + \|f\|_{L^2(M)}^2\).
\end{multline*}
Using $\chi^2(t)$ as a cut-off function and employing \eqref{eq:outer14}, we find
\begin{align}
\iint_{M_{\tau}} \chi^2(t) \psi_t^2 u_{\mx}^{p-1} dv_{g_0} dt
&\le C \left[\int_M (\mcw_{\mx} \psi^2)(z,\tau+1)\, dv_{g_0} + \iint_{M_{\tau}} \chi(t) (|\nabla_{g_0} \psi|_{g_0}^2 + \psi^2)\, dv_{g_0} dt \right. \nonumber \\
&\hspace{160pt} \left. + \mu_0^2(t) \|\psi\|_{L^2(M_{\tau})}^2 + \|f\|_{L^2(M_{\tau})}^2 \right] \label{eq:outer16} \\
&\le C\(\|\psi\|_{L^2(M_{\tau})}^2 + \|f\|_{L^2(M_{\tau})}^2\) \nonumber
\end{align}
for all $\tau \in (t_0,s_0-1]$. Furthermore, we have
\begin{equation}\label{eq:outer17}
\|\psi_t\|_{L^2(M_{t_0})}^2 \le C\(\|\psi\|_{L^2(M_{t_0})}^2 + \|f\|_{L^2(M_{t_0})}^2 + \|\psi_0\|_{\mch^1}^2\)
\end{equation}
as in \eqref{eq:outer15}.

\medskip \noindent \textsc{Step 4.} Combining \eqref{eq:outer14}--\eqref{eq:outer17}, using \eqref{eq:outer12}, and performing algebraic manipulations, we conclude that
\[\|\psi_t\|_{L^2_{t_0,s_0}} + \|\psi\|_{H^1_{t_0,s_0}} \le C\(\|\psi\|_{L^2_{t_0,s_0}} + \|f\|_{L^2_{t_0,s_0}} + \|\psi_0\|_{\mch^1}\).\]
By \eqref{eq:inhom-o} and the previous estimate,
\begin{align*}
\|u_{\mx}^{1-p} L_{g_0,h} \psi\|_{L^2_{t_0,s_0}} &\le C\(\|\psi_t\|_{L^2_{t_0,s_0}} + \|\psi\|_{L^2_{t_0,s_0}} + \|f\|_{L^2_{t_0,s_0}}\) \\
&\le C\(\|\psi\|_{L^2_{t_0,s_0}} + \|f\|_{L^2_{t_0,s_0}} + \|\psi_0\|_{\mch^1}\).
\end{align*}
Adding the above two estimates, we immediately obtain \eqref{eq:outer1}.
\end{proof}

Next, we improve Lemma \ref{lemma:outer1} by dropping the norm of $\psi$ in the right-hand side of \eqref{eq:outer1}. We need some preliminary definitions.
\begin{defn}\label{defn:stereo}
Let $\Pi$ be the inverse of the stereographic projection
\[\Pi(y) = \(\frac{2y}{1+|y|^2}, -\frac{1-|y|^2}{1+|y|^2}\) \in \S^N_n := \S^N \setminus \{(0,\ldots,0,1)\} \subset \R^{N+1}\]
for $y \in \R^N$, and $\Pi_*f: \S^N_n \to \R$ be a weighted push-forward of a function $f: \R^N \to \R$ given by
\[(\Pi_*f)(\ty) = \(\frac{1+|y|^2}{2}\)^{N-2 \over 2} f(y) \quad \text{for } \ty := \Pi(y) \in \S^N_n.\]
\end{defn}
\begin{lemma}\label{lemma:outer1a}
Suppose that all the conditions of Lemma \ref{lemma:outer1} hold, and $s_0 > \frac{3t_0}{2}$.
Then there is a constant $C > 0$ depending only on $(M,g_0)$, $N$, $h$, and $z_0$ such that every solution $\psi$ to \eqref{eq:inhom-o2} satisfies
\begin{equation}\label{eq:outer1a}
\|\psi\|_{H^2_{t_0,s_0}} \le C\(\|f\|_{L^2_{t_0,s_0}} + \|\psi_0\|_{\mch^1}\)
\end{equation}
provided $t_0 > 0$ is large enough.
\end{lemma}
\begin{proof}
Thanks to \eqref{eq:outer1}, it suffices to show that
\begin{equation}\label{eq:outer1a1}
\|\psi\|_{L^2_{t_0,s_0}} \le C\(\|f\|_{L^2_{t_0,s_0}} + \|\psi_0\|_{\mch^1}\).
\end{equation}
Throughout the proof, we assume that $C > 0$ depends only on $(M,g_0)$, $N$, $h$, and $z_0$.

\medskip \noindent \textsc{Step 1.} To establish \eqref{eq:outer1a1}, we argue by contradiction.
Suppose that there are increasing sequences $\{t_{\ell}\}_{\ell \in \N}$, $\{s_{\ell}\}_{\ell \in \N}$ of positive numbers such that
\begin{equation}\label{eq:outer1a11}
s_{\ell} > \frac{3t_{\ell}}{2} \quad \text{and} \quad t_{\ell},\, s_{\ell} \to \infty \quad \text{as } \ell \to \infty,
\end{equation}
sequences $\{\mu_{\ell}\}_{\ell \in \N}$ and $\{\xi_{\ell}\}_{\ell \in \N}$ of parameters satisfying \eqref{eq:mu0}--\eqref{eq:lx},
and sequences $\{\psi_{\ell}\}_{\ell \in \N}$, $\{f_{\ell}\}_{\ell \in \N}$, $\{\psi_{0\ell}\}_{\ell \in \N}$ of functions which satisfy
\begin{equation}\label{eq:outer1a12}
\begin{cases}
\displaystyle p u_{\mxe}^{p-1} (\psi_{\ell})_t = \frac{N+2}{4} L_{g_0,h} \psi_{\ell} + \mcw_{\mxe} \psi_{\ell} + u_{\mxe}^{p-1} f_{\ell} &\text{on } M \times (t_{\ell},s_{\ell}), \\
\psi_{\ell}(\cdot,t_{\ell}) = \psi_{0\ell} &\text{on } M
\end{cases}
\end{equation}
and
\begin{equation}\label{eq:outer1a2}
\|\psi_{\ell}\|_{L^2_{t_{\ell},s_{\ell}}} = 1, \quad \|f_{\ell}\|_{L^2_{t_{\ell},s_{\ell}}} + \|\psi_{0\ell}\|_{\mch^1} \to 0 \quad \text{as } \ell \to \infty.
\end{equation}
There exists $\tau_{\ell} \in [t_{\ell},s_{\ell}-1]$ such that
\begin{equation}\label{eq:outer1a3}
\frac{1}{2} \le \iint_{M_{\tau_{\ell}}} \psi_{\ell}^2 u_{\mxe}^{p-1} dv_{g_0} dt \le 1
\quad \text{and} \quad
\|\psi_{\ell}\|_{H^1_{t_{\ell},s_{\ell}}} \le C
\end{equation}
where the latter inequality is valid by virtue of Lemma \ref{lemma:outer1}.
By passing to a subsequence if necessary, we may assume that $\{\tau_{\ell}\}_{\ell \in \N}$ is increasing.

\medskip
The main assertion of this step is
\begin{equation}\label{eq:outer1a4}
\liminf_{\ell \to \infty} (\tau_{\ell} - t_{\ell}) = \infty.
\end{equation}

Let us prove it. We define
\[\Psi_{\ell}(\tau) := \int_M \(\psi_{\ell}^2 u_{\mxe}^{p-1}\)(\cdot,\tau) dv_{g_0} \quad \text{for } \tau \in [t_{\ell},s_{\ell}-1].\]
If we integrate \eqref{eq:outer13} over $s \in [t_{\ell},\tau]$, we find
\[\Psi_{\ell}(\tau) \le C\left[\int_{t_{\ell}}^{\tau} \Psi_{\ell}(s)ds + (\tau-t_{\ell}+1)\(\|f_{\ell}\|_{L^2_{t_{\ell},s_{\ell}}}^2 + \|\psi_{0\ell}\|_{\mch^1}^2\)\right].\]
By Gr\"onwall's inequality and \eqref{eq:outer1a2}, it follows that
\[\Psi_{\ell}(\tau) \le o(1) \cdot (\tau-t_{\ell}+1)e^{C(\tau-t_{\ell})}.\]
This and \eqref{eq:outer1a3} yield the assertion \eqref{eq:outer1a4}.

\medskip \noindent \textsc{Step 2.} Set
\begin{align*}
\phi_{\ell}(z,t) = \psi_{\ell}(z,t+\tau_{\ell}), \quad u_{\ell}(z,t) = u_{\mxe}(z,t+\tau_{\ell}),\\
\mcw_{\ell}(z,t) = \mcw_{\mxe}(z,t+\tau_{\ell}), \quad \tf_{\ell}(z,t) = f_{\ell}(z,t+\tau_{\ell}).
\end{align*}
We infer from \eqref{eq:outer1a12}, \eqref{eq:outer1a3} and \eqref{eq:outer1a4} that
\[
p u_{\ell}^{p-1} (\phi_{\ell})_t = \frac{N+2}{4} L_{g_0,h} \phi_{\ell} + \mcw_{\ell}\, \phi_{\ell} + u_{\ell}^{p-1} \tf_{\ell} \quad \text{on } M \times (t_{\ell}-\tau_{\ell},1)\]
where $t_{\ell}-\tau_{\ell} \to -\infty$ as $\ell \to \infty$, and
\begin{equation}\label{eq:outer1a6}
\frac{1}{2} \le \iint_{M_0} \phi_{\ell}^2 u_{\ell}^{p-1} dv_{g_0} dt \le 1 \quad \text{and} \quad \|\psi_{\ell}\|_{H^1_{t_{\ell}-\tau_{\ell},1}} \le C.
\end{equation}
Moreover, by \eqref{eq:u2}, \eqref{eq:outer1a3} and the Sobolev inequality, there exists a large number $R_1 > 0$ such that
\begin{equation}\label{eq:outer1a7}
\begin{aligned}
&\ \int_0^1 \int_{M \setminus B_{g_0}(z_0,R_1\tmu_{\ell})} \phi_{\ell}^2 u_{\ell}^{p-1} dv_{g_0} dt \\
&\le C \int_0^1 \left[\int_{M \setminus B_{g_0}(z_0,\delta_0)} \delta_0^{-2N} \mu_0^N(\tau_{\ell}) dv_{g_0}
+ \int_{\{R_1 \le |y-\txi_{\ell}| \le \delta_0 \tmu_{\ell}^{-1}\}} W_{1,0}^{p+1}(y)dy \right]^{2 \over N} dt \\
&\le C \(\delta_0^{-4} \mu_0^2(\tau_{\ell}) + R_1^{-2}\) \le \frac{1}{4}
\end{aligned}
\end{equation}
for all large $\ell \in \N$. Here, $\tmu_{\ell} := \mu_{\ell}(\cdot+\tau_{\ell})$ and $\txi_{\ell} := \xi_{\ell}(\cdot+\tau_{\ell})$.

\medskip \noindent \textsc{Step 3.} According to \eqref{eq:outer1a6} and \eqref{eq:outer1a7}, \begin{equation}\label{eq:outer1a8}
\int_0^1 \int_{B_{g_0}(z_0,R_1\tmu_{\ell})} \phi_{\ell}^2 u_{\ell}^{p-1} dv_{g_0} dt \ge \frac{1}{4}.
\end{equation}
Identifying $z \in B_{g_0}(z_0,\frac{\delta_0}{4})$ and $x = \tmu_{\ell} y+\txi_{\ell} \in B^N(0,\frac{\delta_0}{4})$ as before, we define
\[\vph_{\ell}(y,t) = \tmu_{\ell}^{N-2 \over 2} \phi_{\ell}(\tmu_{\ell} y+\txi_{\ell},t) \quad \text{for } (y,t) \in B_{\ell} \times (t_{\ell}-\tau_{\ell},1)\]
where $B_{\ell} := B^N\(-\tmu_{\ell}^{-1}\,\txi_{\ell}, \frac{\delta_0}{4} \tmu_{\ell}^{-1}\)$. It solves \begin{align*}
&\ p (1+P)^{p-1}(x) \(W_{1,0} + \tbmu_{\ell}^2 Q_0\)^{p-1} (\vph_{\ell})_t \\
&= \frac{N+2}{4} \left[\kappa_N \Delta_{g_0(x)} \vph_{\ell} - \tmu_{\ell}^2 (S+h) \vph_{\ell}\right] + \tmu_{\ell}^2 \ovmcw_{\ell} \vph_{\ell}
+ (1+P)^{p-1}(x) \(W_{1,0} + \tbmu_{\ell}^2 Q_0\)^{p-1} \bf_{\ell} \\
&\ + p (1+P)^{p-1}(x) \(W_{1,0} + \tbmu_{\ell}^2 Q_0\)^{p-1}
\left[\tmu_{\ell}^{-1} \dot{\tmu}_{\ell} \(y \cdot \nabla \vph_{\ell} + \frac{N-2}{2} \vph_{\ell}\) + \tmu_{\ell}^{-1} \dot{\txi}_{\ell} \cdot \nabla \vph_{\ell}\right]
\end{align*}
in $B_{\ell} \times (t_{\ell}-\tau_{\ell},1)$, where $P$ is the function in \eqref{eq:P_20}, $\tbmu_{\ell} := \bmu_{\ell}(\cdot+\tau_{\ell})$,
\[\ovmcw_{\ell}(y,t) := \mcw_{\ell}(\tmu_{\ell} y+\txi_{\ell},t) \quad \text{and} \quad \bf_{\ell}(y,t) := \tf_{\ell}(\tmu_{\ell} y+\txi_{\ell},t).\]
Also, in light of \eqref{eq:outer1a6} and \eqref{eq:outer1a8},
\[\int_0^1 \int_{B^N(0,2R_1)} \vph_{\ell}^2 W_{1,0}^{p-1} dy dt \ge C \quad \text{and} \quad
\sup_{\tau \in [t_{\ell}-\tau_{\ell},0]} \int_{\tau}^{\tau+1} \int_{B_{\ell}} |\nabla \vph_{\ell}|^2 dydt \le C.\]
Hence, there exists a function $\vph_{\infty}$ in $\R^N \times (-\infty,1)$ such that for each $\tau \in [t_{\ell}-\tau_{\ell},0]$,
\[\vph_{\ell} \to \vph_{\infty} \ \begin{cases}
\text{weakly in } \dot{W}^{1,2}(\R^N \times (\tau,\tau+1)) \\
\text{strongly in } L^2_{\text{loc}}(\R^N \times (\tau,\tau+1)) \\
\text{a.e. in } \R^N \times (\tau,\tau+1)
\end{cases}
\text{as } \ell \to \infty,\]
up to a subsequence. In particular, exploiting \eqref{eq:outer1a2} and \eqref{eq:mu0}--\eqref{eq:lx}, we see that
\begin{equation}\label{eq:outer1a9}
pW_{1,0}^{p-1}(\vph_{\infty})_t = \frac{(N+2)\kappa_N}{4} \Delta \vph_{\infty} \quad \text{in } \R^N \times (-\infty,1)
\end{equation}
as well as
\begin{equation}\label{eq:outer1a91}
\int_0^1 \int_{B^N(0,2R_1)} \vph_{\infty}^2 W_{1,0}^{p-1} dy dt \ge C \quad \text{and} \quad \sup_{\tau \in (-\infty,0]} \int_{\tau}^{\tau+1} \int_{\R^N} |\nabla \vph_{\infty}|^2 dy dt \le C.
\end{equation}
In the next step, we will deduce that
\begin{equation}\label{eq:outer1a92}
\vph_{\infty} = 0 \quad \text{in } \R^N \times (-\infty,1),
\end{equation}
a contradiction to the first inequality in \eqref{eq:outer1a91}. Therefore, \eqref{eq:outer1a1} is valid.

\medskip \noindent \textsc{Step 4.}
Referring to Definition \ref{defn:stereo}, we write $\bvp_{\infty} = \Pi_* \vph_{\infty}$ on $\S^N_n$.
From \eqref{eq:outer1a9}, the second inequality in \eqref{eq:outer1a91}, and the conformal covariance of conformal Laplacians, we observe
\begin{equation}\label{eq:outer1a93}
(\bvp_{\infty})_t = \frac{\kappa_N}{N} \left[\Delta_{\S^N} \bvp_{\infty} - \frac{N(N-2)}{4} \bvp_{\infty}\right] \quad \text{in } \S^N_n \times (-\infty,1)
\end{equation}
and
\begin{equation}\label{eq:outer1a94}
\sup_{\tau \in (-\infty,0]} \int_{\tau}^{\tau+1} \int_{\S^N} \(|\nabla_{\S^N} \bvp_{\infty}|_{\S^N}^2 + \bvp_{\infty}^2\) dS_{\ty} dt \le C.
\end{equation}
Parabolic regularity theory and \eqref{eq:outer1a94} imply that $(0,\ldots,0,1) \in \S^N$ is a removable singularity of $\bvp_{\infty}$ and $\|\bvp_{\infty}\|_{C^{\infty}(\S^N \times (-\infty,1))} \le C$.

Let us apply a well-known argument to prove that $\bvp_{\infty} = 0$ on $\S^N \times (-\infty,1)$. By using \eqref{eq:outer1a93} and Bochner's formula, we obtain
\begin{align*}
\(|\nabla \bvp|^2\)_t &= 2 \la \nabla \bvp, \nabla \bvp_t \ra = \frac{2\kappa_N}{N} \left[\la \nabla \bvp, \nabla \Delta \bvp \ra - \frac{N(N-2)}{4} |\nabla \bvp|^2\right] \\
&= \frac{\kappa_N}{N} \left[\Delta |\nabla \bvp|^2 - 2\left|\nabla^2 \bvp\right|^2 - 2\text{Ric}(\nabla \bvp, \nabla \bvp) - \frac{N(N-2)}{2} |\nabla \bvp|^2\right]
\le \frac{\kappa_N}{N} \Delta |\nabla \bvp|^2
\end{align*}
where $\text{Ric} = (N-1)g_{\S^N}$ is the Ricci curvature on $\S^N$, and we wrote $\bvp = \bvp_{\infty}$, $\nabla = \nabla_{\S^N}$, etc. Hence, for any $\tau_0 \in (-\infty,1)$ fixed and $t > \tau_0$,
\begin{align*}
\(\bvp^2 + \frac{2\kappa_N}{N} (t-\tau_0) |\nabla \bvp|^2\)_t
&\le \frac{2\kappa_N}{N} \left[\bvp \left\{\Delta \bvp - \frac{N(N-2)}{4} \bvp\right\} + |\nabla \bvp|^2 + \frac{\kappa_N}{N}(t-\tau_0) \Delta |\nabla \bvp|^2 \right] \\
&\le \frac{\kappa_N}{N} \Delta \left[\bvp^2 + \frac{2\kappa_N}{N}(t-\tau_0) |\nabla \bvp|^2\right].
\end{align*}
By the maximum principle,
\[\bvp^2 + \frac{2\kappa_N}{N} (t-\tau_0) |\nabla \bvp|^2 \le \max_{\ty \in \S^N} \bvp^2(\ty,\tau_0) \le C.\]
Thus
\[|\nabla \bvp(\ty,t)| \le \frac{C}{\sqrt{t-\tau_0}} \quad \text{for all } (\ty,t) \in \S^N \times [\tau_0,1).\]
Taking $\tau_0 \to -\infty$ shows that $\bvp_{\infty}$ depends only on $t$. In view of \eqref{eq:outer1a93},
\[\bvp(t) = \bc_1 e^{-\frac{\kappa_N(N-2)}{4}\, t} + \bc_2 \quad \text{in } (-\infty,1) \quad \text{for some } \bc_1, \bc_2 \in \R.\]
Since it must be bounded in $(-\infty,1)$, we have that $\bc_1 = 0$. From \eqref{eq:outer1a93} again, we conclude that $\bc_2 = 0$.

Now, \eqref{eq:outer1a92} follows at once.
\end{proof}

Using a priori estimate \eqref{eq:outer1a}, we deduce the unique existence of a solution to \eqref{eq:inhom-o}.
\begin{cor}\label{cor:outer1}
Suppose that $Y_h(M,g_0) > 0$, $\|f\|_{L^2_{t_0}} + \|\psi_0\|_{\mch^1} < \infty$ and \eqref{eq:mu0}--\eqref{eq:lx} hold.
Then there exists a unique solution $\psi$ to \eqref{eq:inhom-o} satisfying
\begin{equation}\label{eq:outer1b}
\|\psi\|_{H^2_{t_0}} \le C\(\|f\|_{L^2_{t_0}} + \|\psi_0\|_{\mch^1}\).
\end{equation}
\end{cor}
\begin{proof}
Pick an increasing sequence $\{s_{\ell}\}_{\ell \in \N}$ such that $\frac{3t_0}{2} < s_{\ell} \to \infty$ as $\ell \to \infty$,
and let $\psi_{\ell}$ be a unique solution to a uniformly parabolic equation \eqref{eq:inhom-o2} with $s_0 = s_{\ell}$.
Then a priori estimate \eqref{eq:outer1a} with $s_0 = s_{\ell}$ is true for $\psi_{\ell}$, and the constant $C > 0$ in \eqref{eq:outer1a} is independent of $\ell \in \N$.
Thus, passing to a subsequence, $\psi_{\ell}$ converges weakly to a function $\psi_{\infty}$ in a Banach space equipped with the norm $\|\cdot\|_{H^2_{t_0,t}}$ for each $t > \frac{3t_0}{2}$.
As a result, $\psi_{\infty}$ is the only function satisfying \eqref{eq:inhom-o} and \eqref{eq:outer1b}.
Setting $\psi = \psi_{\infty}$, we conclude the proof.
\end{proof}

\subsection{Inhomogeneous problem associated to \eqref{eq:outer}: Pointwise estimate} \label{subsec:outer2}
We next derive a pointwise estimate for a solution to the inhomogeneous problem \eqref{eq:inhom-o}.
A main tool is a version of the maximum principle presented in the following lemma.
\begin{lemma}\label{lemma:mp}
Suppose that $\psi$ satisfies
\[\begin{cases}
\displaystyle p u_{\mx}^{p-1} \psi_t \ge \frac{N+2}{4} L_{g_0,h} \psi + \mcw_{\mx} \psi &\text{on } M \times (t_0,\infty), \\
\psi(\cdot,t_0) \ge 0 &\text{on } M
\end{cases}\]
in a weak sense. Then $\psi \ge 0$ on $M \times [t_0,\infty)$.
\end{lemma}
\begin{proof}
Fix $s_0 > t_0$ and let $M(s_0) = M \times (t_0,s_0)$. We set
\[\phi(z,t) = e^{-C_0t} \psi(z,t) \quad \text{for } (z,t) \in M \times [t_0,\infty)\]
where
\[C_0 := p^{-1} \(\min_{M(s_0)} u_{\mx}^{p-1}\)^{-1} \left[\frac{N+2}{4} \max_M (|S|+|h|) + \max_{M(s_0)} |\mcw_{\mx}| + 1\right].\]
Then
\begin{equation}\label{eq:mp1}
\begin{aligned}
0 &\ge - \frac{N+2}{4} \left[\kappa_N \iint_{M(s_0)} \la \nabla_{g_0} \phi, \nabla_{g_0} \omega \ra_{g_0}
+ \iint_{M(s_0)} (S+h) \phi \omega \right] + \iint_{M(s_0)} \mcw_{\mx} \phi \omega \\
&\ + p(p-1) \iint_{M(s_0)} u_{\mx}^{-1} (u_{\mx})_t \phi \omega u_{\mx}^{p-1}
- C_0 p \iint_{M(s_0)} \phi \omega u_{\mx}^{p-1} + p \iint_{M(s_0)} \phi \omega_t u_{\mx}^{p-1} \\
&\ - p \int_M (\phi \omega u_{\mx}^{p-1})(\cdot,s_0) + p \int_M (\phi \omega u_{\mx}^{p-1})(\cdot,t_0)
\end{aligned}
\end{equation}
for any nonnegative regular function $\omega$ on $M(s_0)$, and
\[f_0 := \frac{N+2}{4} (S+h) - \mcw_{\mx} + p(p-1) u_{\mx}^{-1} (u_{\mx})_t \cdot u_{\mx}^{p-1} + C_0 p u_{\mx}^{p-1} > 0 \quad \text{on } M(s_0)\]
by \eqref{eq:ArBe}.

Fix a function $\omega_0 \in C^{\infty}(M)$ such that $\omega_0 \ge 0$ on $M$.
By standard parabolic theory, the backward uniformly parabolic equation
\[\begin{cases}
\displaystyle p u_{\mx}^{p-1} \omega_t + \frac{N+2}{4} \kappa_N \Delta_{g_0} \omega - f_0 \omega = 0 &\text{on } M(s_0), \\
\omega(\cdot,s_0) = \omega_0 &\text{on } M
\end{cases}\]
has a unique solution $\omega \ge 0$ on $M(s_0)$. Taking it as a test function for \eqref{eq:mp1} and employing $\phi(\cdot,t_0) \ge 0$ on $M$, we find
\[0 \ge - p \int_M (\phi u_{\mx}^{p-1})(\cdot,s_0) \omega_0 dv_{g_0}.\]
Since $\omega_0$ and $s_0$ were chosen arbitrarily, we must have that $\phi(\cdot,t_0) \ge 0$, namely, $\psi(\cdot,t_0) \ge 0$ on $M \times [t_0,\infty)$.
\end{proof}

In the rest of this subsection, we will deduce a priori pointwise estimates for \eqref{eq:inhom-o} by means of barrier and rescaling arguments.
For the definition of the associated norms, refer to Definition \ref{defn:norm2}.
\begin{lemma}\label{lemma:outer2}
Suppose that $Y_h(M,g_0) > 0$, $\|f\|_{L^2_{t_0}} + \|\psi_0\|_{\mch^1} < \infty$, and \eqref{eq:mu0}--\eqref{eq:lx} hold
so that equation \eqref{eq:inhom-o} admits the unique solution $\psi$; refer to Corollary \ref{cor:outer1}.
Assume further that $\|f\|_{*,\ar} + \|\psi_0\|_{**,\alpha} < \infty$ for some $\alpha \in (0,N-2)$ and $\rho \ge -\frac{N-2}{2}$.
Then there exist constants $C > 0$ large and $\delta_1 > 0$ small depending only on $(M,g_0)$, $N$, $h$, $z_0$, $\alpha$, and $\rho$ such that
\begin{equation}\label{eq:outer2}
|\psi(z,t)| \le C \left[\delta_0^{-2} \|f\|_{*,\ar} \mu_0^{\rho}(t) + \|\psi_0\|_{**,\alpha} \mu_0^{\rho}(t_0) e^{-\delta_1(t-t_0)}\right] w_{\alpha,0}(z,t)
\end{equation}
for $(z,t) \in M \times [t_0,\infty)$.
\end{lemma}
\begin{proof}
Throughout the proof, we assume that $C > 0$ depends only on $(M,g_0)$, $N$, $h$, $z_0$, $\alpha$, and $\rho$.
Moreover, we write $\psi = \psi[f,\psi_0]$, emphasizing the dependence of $\psi$ on $f$ and $\psi_0$.
Let $\psi_1 = \psi[0,\psi_0]$ and $\psi_2 = \psi[f,0]$ so that $\psi = \psi_1 + \psi_2$.

\medskip \noindent \textsc{Step 1: A bound for $\psi_1$.}
First of all, we assume that $d_{g_0}(z,z_0) < \frac{\delta_0}{4}$ so that $u(x,t) = (1+P(x)) \mu^{-{N-2 \over 2}} [W_{1,0}(y) + \bmu^2 Q_0(y)]$ for $x = \exp_{z_0}^{-1}(z) = \mu y + \xi \in B^N(0,\frac{\delta_0}{4})$ and $t \in [t_0,\infty)$.
Given a small number $\delta_1 > 0$, we set
\[\tps_{11}(z,t) = \|\psi_0\|_{**,\alpha} \mu_0^{\rho}(t_0) e^{-\delta_1(t-t_0)} \times
\begin{cases}
2-|y|^2 &\text{for } |y| \le 1,\\
|y|^{-\alpha} &\text{for } |y| > 1.
\end{cases}\]

If $|y| \le 1$, then
\begin{align*}
pu_{\mx}^{p-1} (\tps_{11})_t &= p \mu^{-2} W_{1,0}^{p-1}(y) \(1+O(\mu_0^2)\) \|\psi_0\|_{**,\alpha} \mu_0^{\rho}(t_0) e^{-\delta_1(t-t_0)} \left[-\delta_1 (2-|y|^2) + O(\mu_0^2)\right] \\
&\ge -4 p \alpha_N^{p-1} \delta_1 \|\psi_0\|_{**,\alpha} \mu_0^{\rho}(t_0) e^{-\delta_1(t-t_0)} \mu^{-2},
\end{align*}
\[-L_{g_0,h} \tps_{11} = \|\psi_0\|_{**,\alpha} \mu_0^{\rho}(t_0) e^{-\delta_1(t-t_0)} L_{g_0,h}|y|^2 \ge \frac{3N}{2} \|\psi_0\|_{**,\alpha} \mu_0^{\rho}(t_0) e^{-\delta_1(t-t_0)} \mu^{-2},\]
and
\[\left|\mcw_{\mx} \tps_{11}\right| = o(\mu_0^{-2}) \|\psi_0\|_{**,\alpha} \mu_0^{\rho}(t_0) e^{-\delta_1(t-t_0)}.\]
Thus we obtain
\begin{equation}\label{eq:outer21}
pu_{\mx}^{p-1} (\tps_{11})_t - \frac{N+2}{4} L_{g_0,h} \tps_{11} - \mcw_{\mx} \tps_{11} \ge \frac{N(N+2)}{4} \|\psi_0\|_{**,\alpha} \mu_0^{\rho}(t_0) e^{-\delta_1(t-t_0)} \mu^{-2} > 0
\end{equation}
for sufficiently small $\delta_1 > 0$.

If $|y| > 1$, then
\begin{align*}
pu_{\mx}^{p-1} (\tps_{11})_t &= p \mu^{-2} W_{1,0}^{p-1}(y) \(1+O(\delta_0^2)\) \|\psi_0\|_{**,\alpha} \mu_0^{\rho}(t_0) e^{-\delta_1(t-t_0)} \left[-\delta_1 + O(\mu_0^2)\right] |y|^{-\alpha} \\
&\ge -2 p \alpha_N^{p-1} \delta_1 \|\psi_0\|_{**,\alpha} \mu_0^{\rho}(t_0) e^{-\delta_1(t-t_0)} \mu^{-2} |y|^{-(\alpha+4)}
\end{align*}
and
\[- \frac{N+2}{4} L_{g_0,h} \tps_{11} - \mcw_{\mx} \tps_{11} \ge \left[C + O\(\mu_0^{8 \over N-2} \log\mu_0\)\right] \|\psi_0\|_{**,\alpha} \mu_0^{\rho}(t_0) e^{-\delta_1(t-t_0)} \mu^{-2} |y|^{-(\alpha+2)}\]
where $C > 0$ is a constant depending only on $N$ and $\alpha$. Therefore,
\begin{equation}\label{eq:outer22}
pu_{\mx}^{p-1} (\tps_{11})_t - \frac{N+2}{4} L_{g_0,h} \tps_{11} - \mcw_{\mx} \tps_{11} \ge C \|\psi_0\|_{**,\alpha} \mu_0^{\rho}(t_0) e^{-\delta_1(t-t_0)} \mu^{-2} |y|^{-(\alpha+2)} > 0
\end{equation}
for small $\delta_0,\, \delta_1 > 0$. On the other hand, there exists a large constant $c > 0$ depending only on $(M,g_0)$, $N$, and $\alpha$ such that
\begin{equation}\label{eq:outer22IC}
(c \tps_{11} \pm \psi_1)(z,t_0) = (c\tps_{11})(z,t_0) \pm \psi_0(z) \ge 0
\end{equation}
for $z \in B_{g_0}(z_0,\frac{\delta_0}{4})$.

\medskip
Next, we handle the case that $d_{g_0}(z,z_0) > \frac{\delta_0}{6}$.
Let $v$ be a solution to $-L_{g_0,h}v = 1$ on $M$.
Owing to the condition $Y_h(M,g_0) > 0$, the strong maximum principle holds for the elliptic operator $L_{g_0,h}$.
Hence, such $v$ exists uniquely in $C^2(M)$ and is positive on $M$. Let us set
\[\tps_{12}(z,t) = 2^{3\alpha} \(\min_{z \in M} v(z)\)^{-1} \delta_0^{-\alpha} \|\psi_0\|_{**,\alpha}
\mu_0^{\rho}(t_0) e^{-\delta_1(t-t_0)} \mu^{\alpha}(t) v(z) \quad \text{on } M \times [t_0,\infty).\]
A simple calculation shows that if $\frac{\delta_0}{6} < d_{g_0}(z,z_0) < \frac{\delta_0}{4}$ and $t \in [t_0,\infty)$, then $\tps_{12}(z,t) > \tps_{11}(z,t)$. Besides,
\[pu_{\mx}^{p-1} (\tps_{12})_t \ge - C \delta_1 \delta_0^{-(\alpha+4)} \|\psi_0\|_{**,\alpha} \mu_0^{\rho}(t_0) e^{-\delta_1(t-t_0)} \mu^{\alpha+2}\]
and
\[- \frac{N+2}{4} L_{g_0,h} \tps_{12} - \mcw_{\mx} \tps_{12} = \left[C + O\(\delta_0^{-4} \mu_0^2\)\right] \delta_0^{-\alpha} \|\psi_0\|_{**,\alpha} \mu_0^{\rho}(t_0) e^{-\delta_1(t-t_0)} \mu^{\alpha}.\]
Consequently, we have
\begin{equation}\label{eq:outer23}
pu_{\mx}^{p-1} (\tps_{12})_t - \frac{N+2}{4} L_{g_0,h} \tps_{12} - \mcw_{\mx} \tps_{12} \ge C \delta_0^{-\alpha} \|\psi_0\|_{**,\alpha} \mu_0^{\rho}(t_0) e^{-\delta_1(t-t_0)} \mu^{\alpha} > 0
\end{equation}
taking $\delta_1 > 0$ smaller if needed. In addition, it holds that
\begin{equation}\label{eq:outer23IC}
(c\tps_{12} \pm \psi_1)(z,t_0) \ge 0 \quad \text{on } M
\end{equation}
for some large $c > 0$.

\medskip
Now, letting $\tps_{11}(z,t) = 0$ for $z \in M \setminus B_{g_0}(z,\frac{\delta_0}{4})$ and $t \in [t_0,\infty)$, we define a function
\[\tps = \max\left\{\tps_{11}, \tps_{12}\right\} \quad \text{on } M \times [t_0,\infty).\]
By applying Lemma \ref{lemma:mp} for functions $c\tps \pm \psi_1$ and using \eqref{eq:outer21}--\eqref{eq:outer23IC}, we deduce
\begin{equation}\label{eq:psi1}
|\psi_1(z,t)| \le C \|\psi_0\|_{**,\alpha} \mu_0^{\rho}(t_0) e^{-\delta_1(t-t_0)} w_{\alpha,0}(z,t) \quad \text{on } M \times [t_0,\infty).
\end{equation}

\medskip \noindent \textsc{Step 2: A bound for $\psi_2$.}
For $z \in B_{g_0}(z_0,\frac{\delta_0}{4})$, we redefine
\[\tps_{21}(z,t) = \delta_0^{-2} \|f\|_{*,\ar} \mu^{\rho}(t) \times \begin{cases}
2-|y|^2 &\text{for } |y| \le 1,\\
|y|^{-\alpha} &\text{for } |y| > 1
\end{cases}\]
where $x = \exp^{-1}(z) = \mu y + \xi \in B^N(0,\frac{\delta_0}{4})$ and $t \in [t_0,\infty)$.

To treat the case that $d_{g_0}(z,z_0) > \frac{\delta_0}{6}$, we also set
\[\tps_{22}(z,t) = 2^{3\alpha} \(\min_{z \in M} v(z)\)^{-1} \delta_0^{-(\alpha+2)} \|f\|_{*,\ar} \mu^{\alpha+\rho}(t) v(z) \quad \text{on } M \times [t_0,\infty).\]

Then, arguing as in Step 1, we arrive at
\begin{equation}\label{eq:psi2}
|\psi_2(z,t)| \le C \delta_0^{-2} \|f\|_{*,\ar} \mu_0^{\rho}(t) w_{\alpha,0}(z,t) \quad \text{on } M \times [t_0,\infty).
\end{equation}

\medskip
Summing \eqref{eq:psi1} and \eqref{eq:psi2} up, we get \eqref{eq:outer2}.
\end{proof}

\begin{cor}\label{cor:outer2}
Suppose that all the conditions of Lemma \ref{lemma:outer2} hold,
and $\|f\|_{*,\ar;\sigma_0} + \|\psi_0\|_{**,\alpha;\sigma_0} < \infty$ where $\sigma_0 \in (0,1)$ is the number appearing in \eqref{eq:lx}.
There exist constants $C > 0$ large and $\delta_1 > 0$ small depending only on $(M,g_0)$, $N$, $h$, $z_0$, $\alpha$, $\rho$, and $\sigma_0$ such that
\begin{equation}\label{eq:outer2a}
\begin{cases}
\left|\nabla_{g_0}^{\ell} \psi(z,t)\right| \le C \left[\delta_0^{-2} \|f\|_{*,\ar;\sigma_0} \mu_0^{\rho}(t) + \|\psi_0\|_{**,\alpha;\sigma_0} \mu_0^{\rho}(t_0) e^{-\delta_1(t-t_0)}\right] w_{\alpha,\ell}(z,t), \\
\left[\nabla_{g_0}^{\ell} \psi\right]_{C^{\sigma_0}_z}(z,t) \le C \left[\delta_0^{-2} \|f\|_{*,\ar;\sigma_0} \mu_0^{\rho}(t) + \|\psi_0\|_{**,\alpha;\sigma_0} \mu_0^{\rho}(t_0) e^{-\delta_1(t-t_0)}\right] w_{\alpha,\ell+\sigma_0}(z,t), \\
\left[\nabla_{g_0}^{\ell} \psi\right]_{C^{\sigma_0/2}_t}(z,t) \le C \left[\delta_0^{-2} \|f\|_{*,\ar;\sigma_0} \mu_0^{\rho}(t) + \|\psi_0\|_{**,\alpha;\sigma_0} \mu_0^{\rho}(t_0) e^{-\delta_1(t-t_0)}\right] w_{\alpha-\sigma_0,\ell}(z,t),
\end{cases}
\end{equation}
and
\begin{equation}\label{eq:outer2b}
\begin{cases}
|\psi_t(z,t)| \le C \left[\delta_0^{-2} \|f\|_{*,\ar;\sigma_0} \mu_0^{\rho}(t) + \|\psi_0\|_{**,\alpha;\sigma_0} \mu_0^{\rho}(t_0) e^{-\delta_1(t-t_0)}\right] w_{\alpha-2,0}(z,t), \\
[\psi_t]_{C^{\sigma_0}_z}(z,t) \le C \left[\delta_0^{-2} \|f\|_{*,\ar;\sigma_0} \mu_0^{\rho}(t)
+ \|\psi_0\|_{**,\alpha;\sigma_0} \mu_0^{\rho}(t_0) e^{-\delta_1(t-t_0)}\right] w_{\alpha-2,\sigma_0}(z,t), \\
[\psi_t]_{C^{\sigma_0/2}_t}(z,t) \le C \left[\delta_0^{-2} \|f\|_{*,\ar;\sigma_0} \mu_0^{\rho}(t)
+ \|\psi_0\|_{**,\alpha;\sigma_0} \mu_0^{\rho}(t_0) e^{-\delta_1(t-t_0)}\right] w_{\alpha-2-\sigma_0,0}(z,t)
\end{cases}
\end{equation}
for $(z,t) \in M \times [t_0,\infty)$ and $\ell \in \{0,1,2\}$. In particular, we have that
\begin{equation}\label{eq:outer2a0}
\|\psi\|_{*',\ar;\sigma_0} \le C_1\delta_0^{-2}\(\|f\|_{*,\ar;\sigma_0} + \|\psi_0\|_{**,\alpha;\sigma_0}\)
\end{equation}
where $C_1 > 0$ is a constant depending only on $(M,g_0)$, $N$, $h$, $z_0$, $\alpha$, $\rho$, and $\sigma_0$.
\end{cor}
\begin{proof}
The proof is decomposed into two steps.

\medskip \noindent \textsc{Step 1.} We examine the case that $d_{g_0}(z,z_0) \le \frac{\delta_0}{12}$.

We set
\[\phi(y,t) = \psi(z,t) \quad \text{for } y 
\in B^N\(-\xi,\frac{\delta_0}{4\mu}\) \text{ and } t \in [t_0,\infty).\]
It solves
\begin{multline}\label{eq:outer24}
\phi_t = \frac{N-2}{4} \mu^{-2} u_{\mx}^{1-p}(x,t) \Delta_{g_0(x)} \phi + \mu^{-1}\(\dot{\mu} y + \dot{\xi}\) \cdot \nabla \phi \\
+ p^{-1} u_{\mx}^{1-p}(x,t) \left[\frac{N+2}{4}(S+h)(x) + \mcw_{\mx}(x,t)\right] \phi + p^{-1} f(x,t).
\end{multline}

For a number $t > t_0+2$ and a set $\Omega \subset \R^{N+1}$, we define
\[\Omega_{11} = \{|y| \le 6\} \times [t-1,t], \quad \Omega_{12} = \{|y| \le 8\} \times [t-2,t],\]
and a standard parabolic H\"older norm
\[\|\phi\|_{C^{\sigma_0,\sigma_0/2}(\Omega)} = \|\phi\|_{L^{\infty}(\Omega)} + \sup_{(y_1,t_1) \ne (y_2,t_2) \in \Omega}
\frac{|\phi(y_1,t_1)-\phi(y_2,t_2)|}{|y_1-y_2|^{\sigma_0} + |t_1-t_2|^{\sigma_0/2}}.\]
A straightforward computation shows that the $C^{\sigma_0,\sigma_0/2}(\Omega_{12})$-norms of the coefficients in \eqref{eq:outer24} are uniformly bounded in $t$.
Therefore, the interior parabolic Schauder estimate is applicable, which gives
\begin{equation}\label{eq:outer25}
\sum_{\ell=0}^2 \left\|\nabla^{\ell} \phi\right\|_{C^{\sigma_0,\sigma_0/2}(\Omega_{11})} + \|\phi_t\|_{C^{\sigma_0,\sigma_0/2}(\Omega_{11})}
\le C \(\|\phi\|_{L^{\infty}(\Omega_{12})} + \|f(x,t)\|_{C^{\sigma_0,\sigma_0/2}(\Omega_{12})}\).
\end{equation}
Given any $(y_1,t_1),\, (y_2,t_2) \in \Omega_{12}$, we have
\[\frac{|f(x_1,t_1)-f(x_2,t_2)|}{|y_1-y_2|^{\sigma_0} + |t_1-t_2|^{\sigma_0/2}}
\le (1+O(\delta_0^2)) \mu^{\sigma_0}(t_1) \frac{|f(x_1,t_1)-f(x_{12},t_1)|}{d_{g_0}(x_1,x_{12})^{\sigma_0}} + \frac{|f(x_{12},t_1)-f(x_2,t_2)|}{|t_1-t_2|^{\sigma_0/2}}\]
where
\[x_1 := \mu(t_1)y_1+\xi(t_1), \quad x_2 := \mu(t_2)y_2+\xi(t_2) \quad \text{and} \quad x_{12} := \mu(t_1)y_2+\xi(t_1).\]
Moreover, if $|y_1-y_2| < \delta$ for a number $\delta > 0$ small enough, then $d_{g_0}(x_1,x_2) < \mu_0(t_1)$ and so
\begin{align*}
|f(x_1,t_1)-f(x_{12},t_1)|
&\le u_{\mx}^{1-p}(x_1,t_1) \left[u_{\mx}^{p-1}f\right]_{C^{\sigma_0}_z}(x_1,t_1)\, d_{g_0}(x_1,x_{12})^{\sigma_0} \\
&\ + \sup_{\theta \in (0,1)} \left|\nabla_x (u_{\mx}^{1-p})(x_1+\theta(x_{12}-x_1),t_1)\right| \left|(u_{\mx}^{p-1}f)(x_{12},t_1)\right| \\
&\le C\left[\mu^{-\sigma_0}(t_1) \|f\|_{*,\ar;\sigma_0} + d_{g_0}(x_1,x_{12})^{1-\sigma_0} \|f\|_{*,\ar}\right] \mu_0^{\rho}(t_1) d_{g_0}(x_1,x_{12})^{\sigma_0}.
\end{align*}
If $|t_1-t_2| < 1$, then
\[|x_{12}-x_2| \le C\mu_0^3(t_1)(|y_2|+o(1))|t_1-t_2| \le C\mu_0^3(t_1)|t_1-t_2|,\]
which implies
\begin{align*}
|f(x_{12},t_1)-f(x_2,t_2)| &\le |f(x_{12},t_1)-f(x_2,t_1)| + |f(x_2,t_1)-f(x_2,t_2)| \\
&\le C\left[\mu_0^{2\sigma_0}(t_1) |t_1-t_2|^{\sigma_0/2} \|f\|_{*,\ar;\sigma_0} + \|f\|_{*,\ar}\right] \mu_0^{\rho}(t_1)|t_1-t_2|^{\sigma_0/2} \\
&\le C \|f\|_{*,\ar;\sigma_0} \mu_0^{\rho}(t_1)|t_1-t_2|^{\sigma_0/2}.
\end{align*}
It follows that
\begin{equation}\label{eq:outer26}
\|f(x,t)\|_{C^{\sigma_0,\sigma_0/2}(\Omega_{12})} \le C \|f\|_{*,\ar;\sigma_0}\mu_0^{\rho}(t).
\end{equation}
Substituting \eqref{eq:outer2} and \eqref{eq:outer26} into \eqref{eq:outer25} gives us that
\begin{equation}\label{eq:outer2a1}
\begin{cases}
\left|\nabla^{\ell}\psi(x,t)\right| = \mu^{-\ell} \left|\nabla^{\ell}\phi(y,t)\right|
\le C \left[\delta_0^{-2} \|f\|_{*,\ar;\sigma_0} \mu_0^{\rho} + \|\psi_0\|_{**,\alpha} \mu_0^{\rho}(t_0) e^{-\delta_1(t-t_0)}\right] \mu_0^{-\ell}, \\
\begin{aligned}
&\mu_0^{\sigma_0} \left[\nabla^{\ell} \psi\right]_{C^{\sigma_0}_z}(x,t) + \left[\nabla^{\ell} \psi\right]_{C^{\sigma_0/2}_t}(x,t) \\
&\hspace{135pt} \le C \left[\delta_0^{-2} \|f\|_{*,\ar;\sigma_0} \mu_0^{\rho} + \|\psi_0\|_{**,\alpha} \mu_0^{\rho}(t_0) e^{-\delta_1(t-t_0)}\right] \mu_0^{-\ell}
\end{aligned}
\end{cases}
\end{equation}
and
\begin{equation}\label{eq:outer2a2}
\begin{cases}
\begin{aligned}
|\psi_t(x,t)| &\le |\phi_t(y,t)| + \mu^{-1}\(\dot{\mu}|y|+\dot{\xi}\) |\nabla\psi(x,t)| \\
&\le C \left[\delta_0^{-2} \|f\|_{*,\ar;\sigma_0} \mu_0^{\rho} + \|\psi_0\|_{**,\alpha} \mu_0^{\rho}(t_0) e^{-\delta_1(t-t_0)}\right],
\end{aligned} \\
\mu_0^{\sigma_0} [\psi_t]_{C^{\sigma_0}_z}(x,t) + [\psi_t]_{C^{\sigma_0/2}_t}(x,t)
\le C \left[\delta_0^{-2} \|f\|_{*,\ar;\sigma_0} \mu_0^{\rho} + \|\psi_0\|_{**,\alpha} \mu_0^{\rho}(t_0) e^{-\delta_1(t-t_0)}\right]
\end{cases}
\end{equation}
for $\ell \in \{0,1,2\}$, $|x-\xi| < 4\mu$, and $t > t_0+2$.

Let $n_0$ be a natural number such that $2^{n_0} \le \frac{\delta_0}{24} \mu^{-1}(t) < 2^{n_0+1}$.
In the rest of this step, we will estimate $\phi$ on each set
\[\Omega_{n1} := \left\{2^n < |y| < 2^{n+1}\right\} \times \left[t-2^{-2n},t\right] \quad \text{for } n = 1, \ldots, n_0.\]
With this aim, we introduce
\[\wtom_{n1} := \{1<|Y|<2\} \times [-1,0], \quad \wtom_{n2} := \left\{\frac{1}{2} < |Y| < \frac{5}{2}\right\} \times [-2,0]\]
and
\[\vph(Y,\tau) = \phi\(2^n\, Y, t+2^{-2n}\tau\) \quad \text{for } (Y,\tau) \in \wtom_{n2}.\]
By \eqref{eq:outer24}, it is a solution to
\begin{align*}
\vph_{\tau} &= \frac{N-2}{4} 2^{-4n} \mu^{-2} u_{\mx}^{1-p}\(x,t+2^{-2n}\tau\) \Delta_{g_0(x)} \vph + 2^{-2n} \mu^{-1} \(\dot{\mu} y + \dot{\xi}\) \cdot \nabla \vph \\
&\ + p^{-1} 2^{-2n} u_{\mx}^{1-p}\(x,t+2^{-2n}\tau\) \left[\frac{N+2}{4}(S+h)(x) + \mcw_{\mx}\(x,t+2^{-2n}\tau\)\right] \vph \\
&\ + 2^{-2n} p^{-1} f\(x,t+2^{-2n}\tau\)
\end{align*}
in $\wtom_{n2}$. A direct computation gives that the above equation is uniformly parabolic
and all the $C^{\sigma_0,\sigma_0/2}(\wtom_{n2})$-norms of the coefficients in it are uniformly bounded in $n$ and $t$. Hence
\begin{multline}\label{eq:outer27}
\sum_{\ell=0}^2 \left\|\nabla^{\ell} \vph\right\|_{C^{\sigma_0,\sigma_0/2}(\wtom_{n1})} + \|\vph_{\tau}\|_{C^{\sigma_0,\sigma_0/2}(\wtom_{n1})} \\
\le C \(\|\vph\|_{L^{\infty}(\wtom_{n2})} + \left\|2^{-2n} f\(x,t+2^{-2n}\tau\)\right\|_{C^{\sigma_0,\sigma_0/2}(\wtom_{n2})}\),
\end{multline}
which can be interpreted as
\begin{equation}\label{eq:outer2a3}
\begin{cases}
\begin{aligned}
\left|\nabla^{\ell} \psi(x,t)\right| &= \(\mu 2^n\)^{-\ell} \left|\nabla^{\ell} \vph(Y,0)\right| \\
&\le C \left[\delta_0^{-2} \|f\|_{*,\ar;\sigma_0} \mu_0^{\rho} + \|\psi_0\|_{**,\alpha} \mu_0^{\rho}(t_0) e^{-\delta_1(t-t_0)}\right] \frac{\mu_0^{-\ell}}{|y|^{\alpha+\ell}},
\end{aligned}
\\
\begin{aligned}
& (\mu_0|y|)^{\sigma_0} \left[\nabla^{\ell} \psi\right]_{C^{\sigma_0}_z}(x,t) + |y|^{-\sigma_0} \left[\nabla^{\ell} \psi\right]_{C^{\sigma_0}_t}(x,t) \\
&\hspace{52pt} \le C \left[\delta_0^{-2} \|f\|_{*,\ar;\sigma_0} \mu_0^{\rho} + \|\psi_0\|_{**,\alpha} \mu_0^{\rho}(t_0) e^{-\delta_1(t-t_0)}\right] \frac{\mu_0^{-\ell}}{|y|^{\alpha+\ell}}
\end{aligned}
\end{cases}
\end{equation}
and
\begin{equation}\label{eq:outer2a4}
\begin{cases}
\begin{aligned}
|\psi_t(x,t)| &\le C\mu^2 |\nabla \vph(Y,\tau)| + 2^{2n} |\vph_{\tau}(Y,\tau)| \\
&\le C \left[\delta_0^{-2} \|f\|_{*,\ar;\sigma_0} \mu_0^{\rho} + \|\psi_0\|_{**,\alpha} \mu_0^{\rho}(t_0) e^{-\delta_1(t-t_0)}\right] \frac{1}{|y|^{\alpha-2}},
\end{aligned} \\
\begin{aligned}
&(\mu_0|y|)^{\sigma_0} [\psi_t]_{C^{\sigma_0}_z}(x,t) + |y|^{-\sigma_0} [\psi_t]_{C^{\sigma_0/2}_t}(x,t) \\
&\hspace{40pt} \le C \left[\delta_0^{-2} \|f\|_{*,\ar;\sigma_0} \mu_0^{\rho} + \|\psi_0\|_{**,\alpha} \mu_0^{\rho}(t_0) e^{-\delta_1(t-t_0)}\right] \frac{1}{|y|^{\alpha-2}}
\end{aligned}
\end{cases}
\end{equation}
for $\ell \in \{0,1,2\}$, $2\mu < |x-\xi| < \frac{\delta_0}{12}$, and $t \ge t_0+2$.
To get the estimates for the $C^{\sigma_0}_t$-norms in \eqref{eq:outer2a3}--\eqref{eq:outer2a4},
we need \eqref{eq:outer27} with $t$ replaced by any number in the interval $[t-1+2^{2n}, t]$.

To treat the case that $t \in [t_0,t_0+2]$, we repeat the above argument with the interior Schauder estimate for Cauchy problems. For example, instead of \eqref{eq:outer25}, we use
\begin{multline*}
\sum_{\ell=0}^2 \left\|\nabla^{\ell} \phi\right\|_{C^{\sigma_0,\sigma_0/2}(\Omega_{11}')} + \|\phi_t\|_{C^{\sigma_0,\sigma_0/2}(\Omega_{11}')} \\
\le C \(\|\phi\|_{L^{\infty}(\Omega_{12}')} + \|f(x,t)\|_{C^{\sigma_0,\sigma_0/2}(\Omega_{12}')} + \sum_{\ell=0}^2 \left\|\nabla^{\ell} \phi_0\right\|_{C^{\sigma_0}(\{|y| \le 8\})}\)
\end{multline*}
where $\phi_0(y) := \psi_0(x)$ for $y \in B^N\(-\xi,\frac{\delta_0}{4\mu}\)$,
\[\Omega_{11}' := \{|y| \le 6\} \times [t_0,t_0+2] \quad \text{and} \quad \Omega_{12}' := \{|y| \le 8\} \times [t_0,t_0+2].\]
Then, the corresponding estimate to \eqref{eq:outer2a1} is
\[\left|\nabla^{\ell}\psi(x,t)\right| + \mu_0^{\sigma_0} \left[\nabla^{\ell} \psi\right]_{C^{\sigma_0}_z}(x,t) + \left[\nabla^{\ell} \psi\right]_{C^{\sigma_0/2}_t}(x,t)
\le C \left[\delta_0^{-2} \|f\|_{*,\ar;\sigma_0} \mu_0^{\rho} + \|\psi_0\|_{**,\alpha;\sigma_0} \mu_0^{\rho}(t_0)\right] \mu_0^{-\ell},\]
and we can deduce an analogous estimate to \eqref{eq:outer2a2}, \eqref{eq:outer2a3}, and \eqref{eq:outer2a4}.

\medskip \noindent \textsc{Step 2.} We examine the case that $d_{g_0}(z,z_0) \ge \frac{\delta_0}{12}$.

Define a function $\tps$ on $M \times [t_0,\infty)$ by
\[\psi(z,t) = \tps(z,s(t)) \quad \text{where } s(t) := \delta_0^4 t^2.\]
It solves
\[\tps_s = \frac{N-2}{8\delta_0^4t} u_{\mx}^{1-p} \Delta_{g_0}\tps + \frac{1}{2p\, \delta_0^4t} u_{\mx}^{1-p} \left[\frac{N+2}{4}(S+h) + \mcw_{\mx}\right] \tps + \frac{1}{2p\, \delta_0^4t} f(z,t).\]
If $d_{g_0}(z,z_0) \ge \frac{\delta_0}{24}$, a coefficient function $(\delta_0^4 t)^{-1} u_{\mx}^{1-p}$ is bounded from above, bounded away from 0
(because $(\delta_0^4 t)^{-1} u_{\mx}^{1-p} \simeq t^{-1}\mu_0^{-2} \simeq 1$), and uniformly Lipschitz.
Therefore, we conclude from the parabolic Schauder estimate (for Cauchy problems) that
\begin{equation}\label{eq:outer2a5}
\begin{cases}
\left|\nabla^{\ell} \psi(z,t)\right| = \left|\nabla^{\ell} \tps(z,s)\right|
\le C \left[\delta_0^{-2} \|f\|_{*,\ar;\sigma_0} \mu_0^{\rho} + \|\psi_0\|_{**,\alpha} \mu_0^{\rho}(t_0) e^{-\delta_1(t-t_0)}\right] \delta_0^{-(\alpha+\ell)} \mu_0^{\alpha}, \\
\begin{aligned}
&\delta_0^{\sigma_0} \left[\nabla^{\ell} \psi\right]_{C^{\sigma_0}_z}(x,t) + (\delta_0^{-1} \mu_0)^{\sigma_0} \left[\nabla^{\ell} \psi\right]_{C^{\sigma_0}_t}(x,t) \\
&\hspace{118pt} \le C \left[\delta_0^{-2} \|f\|_{*,\ar;\sigma_0} \mu_0^{\rho} + \|\psi_0\|_{**,\alpha} \mu_0^{\rho}(t_0) e^{-\delta_1(t-t_0)}\right] \delta_0^{-(\alpha+\ell)} \mu_0^{\alpha}
\end{aligned}
\end{cases}
\end{equation}
for $\ell \in \{0,1,2\}$ and
\begin{equation}\label{eq:outer2a6}
\begin{cases}
|\psi_t(z,t)| \le 2\delta_0^4t|\tps_s(z,s)| \le C \left[\delta_0^{-2} \|f\|_{*,\ar;\sigma_0} \mu_0^{\rho}
+ \|\psi_0\|_{**,\alpha} \mu_0^{\rho}(t_0) e^{-\delta_1(t-t_0)}\right] \delta_0^{2-\alpha} \mu_0^{\alpha-2} \\
\begin{aligned}
&\delta_0^{\sigma_0} [\psi_t]_{C^{\sigma_0}_z}(x,t) + (\delta_0^{-1}\mu_0)^{\sigma_0} [\psi_t]_{C^{\sigma_0}_t}(x,t) \\
&\hspace{105pt} \le C \left[\delta_0^{-2} \|f\|_{*,\ar;\sigma_0} \mu_0^{\rho}
+ \|\psi_0\|_{**,\alpha} \mu_0^{\rho}(t_0) e^{-\delta_1(t-t_0)}\right] \delta_0^{2-\alpha} \mu_0^{\alpha-2}
\end{aligned}
\end{cases}
\end{equation}
provided $d_{g_0}(z,z_0) \ge \frac{\delta_0}{12}$.

\medskip
Combining \eqref{eq:outer2a1}--\eqref{eq:outer2a6}, we find the desired gradient estimates \eqref{eq:outer2a} and \eqref{eq:outer2b}.
\end{proof}

\subsection{Unique solvability of \eqref{eq:outer} and a priori estimate of the solution} \label{subsec:outer3}
Denote $\mcb_{\bmx} := B^N(-\bmu^{-1}\xi, 2\bmu^{-1}\mu_0^{\vep_1})$.
Let $\hps^{\, \tin}$ be a function on $B^N(-0,2\mu_0^{\vep_1}) \times [t_0,\infty)$
and $\psi^{\tin}$ be a function on $\mcb_{\bmx} \times [t_0,\infty)$ such that
\begin{equation}\label{eq:psi00}
\hps^{\, \tin}(x,t) = (1+P(x)) \bmu^{-{N-2 \over 2}} \psi^{\tin}\(\by, t\);
\end{equation}
cf. \eqref{eq:psi0}. We remark that $\psi^{\tin}$ needs not to be a solution to \eqref{eq:inner} at this moment.
Letting also $\psi_0$ be a function on $M$, we assume that
\begin{equation}\label{eq:outer31}
\left\|\psi^{\tin}\right\|_{\sharp',a,b;\sigma_0\(\mcb_{\bmx}\)} \le C_{21}
\quad \text{and} \quad \|\psi_0\|_{**,\alpha;\sigma_0} \le C_{22} \mu_0^{\delta_2}(t_0)
\end{equation}
for some $C_{21},\, C_{22} > 0$, $a \in (\sigma_0,N-2)$, $b \simeq 3$, $\alpha \in (0,a)$ to be determined later, and an arbitrarily chosen $\delta_2 > 0$;
the $\sharp'$- and $**$-norms are given in \eqref{eq:sharp'-norm} and \eqref{eq:**-norm}, respectively,
and $\sigma_0 \in (0,1)$ is the number appearing in \eqref{eq:lx}.
By employing the existence theory and a priori estimates for the inhomogeneous problem \eqref{eq:inhom-o},
we shall prove the unique solvability of the outer problem \eqref{eq:outer}, i.e.,
\begin{equation}\label{eq:outer3}
\begin{cases}
\displaystyle p u_{\mx}^{p-1} \(\psi^{\out}\)_t = \frac{N+2}{4} L_{g_0,h} \psi^{\out} + \mcw_{\mx} \psi^{\out}
+ u_{\mx}^{p-1} f^{\out}_{\mx}\left[\psi^{\out}, \hps^{\, \tin}\right] &\text{on } M \times (t_0,\infty), \\
\psi^{\out}(\cdot,t_0) = \psi_0 &\text{on } M
\end{cases}
\end{equation}
where
\begin{equation}\label{eq:fout}
f^{\out}_{\mx}\left[\psi^{\out}, \hps^{\, \tin}\right] := u_{\mx}^{1-p} \left[\(1-\eta_{\mu_0^{\vep_1}}\)\mcs(u_{\mx})
+ \mcj_1\left[\hps^{\tin}\right] + \mcj_2\left[\psi^{\out}, \hps^{\tin}\right]\right].
\end{equation}

\begin{prop}\label{prop:outer3}
Assume that \eqref{eq:outer31} and \eqref{eq:mu0}--\eqref{eq:lx} hold, and $t_0 > 0$ is sufficiently large.
Suppose also that $\alpha \in (0,a)$ and $\beta \in (0,b]$ satisfy $\alpha+\beta = N-2$,
and set $\rho = -\frac{N-2}{2}+\beta = \frac{N-2}{2}-\alpha$. Choose any
\begin{equation}\label{eq:delta3}
\delta_3 \in \(0, \min\{2,\, \beta,\, \beta(p-1),\, a+b-(N-2),\, (p-1)(a+b-(N-2))\}\right].
\end{equation}
Then \eqref{eq:outer3} possesses a unique solution $\psi^{\out} = \psi^{\out}[\lambda, \xi, \dot{\lambda}, \dot{\xi}, \hps^{\, \tin}, \psi_0]$ such that
\begin{equation}\label{eq:outer32}
\left\|\psi^{\out}\right\|_{*',\ar;\sigma_0} \le C \delta_0^{-2} \mu_0^{\delta_4}(t_0)
\end{equation}
where $C > 0$ is a constant depending only on $(M,g_0)$, $N$, $h$, $z_0$, $\alpha$, and $\sigma_0$,
\begin{equation}\label{eq:delta4}
\delta_4 \in \(0, \min\left\{\nu_2-\zeta_2\vep_1,\, 2(p-1), \delta_3\right\}\) \cap (0,\min\{(1-\vep_1)(a-\alpha), \delta_2\}],
\end{equation}
and the $*'$-norm is defined in \eqref{eq:*'s-norm}.
\end{prop}
\noindent The proof of Proposition \ref{prop:outer3} consists of several lemmas.
Let $\mbone_{N=5}$ be a quantity equal to $1$ if $N = 5$ and $0$ if $N \ge 6$.
\begin{lemma}
We have
\begin{equation}\label{eq:oest1}
\left\|u_{\mx}^{1-p} \(1-\eta_{\mu_0^{\vep_1}}\) \mcs(u_{\mx})\right\|_{*,\ar;\sigma_0}
\le C_{31} \delta_0^{-\zeta_{31}} \mu_0^{\nu_2-\zeta_2\vep_1}(t_0)
\end{equation}
where
\begin{itemize}
\item[-] $C_{31},\, \zeta_{31} > 0$ are constants depending only on $(M,g_0)$, $N$, $h$, $z_0$, $\alpha$, and $\sigma_0$;
\item[-] $\vep_1 \in (0,1)$ is the small number in Lemma \ref{lemma:error1out},
    $\zeta_2 > 0$ is the number depending only on $N$ in Lemma \ref{lemma:error2out},
    and $\nu_2 = 2-\vep_0$ is the number in \eqref{eq:lx};
\item[-] the $*$-norm is defined in \eqref{eq:*-norm}.
\end{itemize}
\end{lemma}
\begin{proof}
Write $\wtmcs = \(1-\eta_{\mu_0^{\vep_1}}\) \mcs(u_{\mx})$. Applying \eqref{eq:S2out} and \eqref{eq:mu0}--\eqref{eq:lx}, we observe
\[\sum_{\ell=0}^1 \left|\nabla_z^{\ell} \wtmcs(z,t)\right|
+ \left[\(1-\eta_{\mu_0^{\vep_1}}\) \wtmcs\, \right]_{C^{\sigma_0/2}_t}(z,t)
\le C \delta_0^{-\zeta} \(\mu_0^{{N+2 \over 2}-\zeta_2\vep_1} + \mu_0^{{N-2 \over 2}+\nu_2-\zeta_2\vep_1}\).\]
Hence
\begin{multline*}
\left[\(\mu_0^{\rho} w_{\alpha,2}\)^{-1} \left|\wtmcs(z,t)\right| + \(\mu_0^{\rho} w_{\alpha,2+\sigma_0}\)^{-1} \left[\wtmcs(z,t)\right]_{C^{\sigma_0}_z}
+ \(\mu_0^{\rho} w_{\alpha-\sigma_0,2}\)^{-1} \left[\wtmcs(z,t)\right]_{C^{\sigma_0/2}_t}\right](z,t) \\
\le C \delta_0^{-\zeta} \mu_0^{-\alpha-\rho} \mu_0^{{N-2 \over 2}+\nu_2-\zeta_2\vep_1}
= C \delta_0^{-\zeta} \mu_0^{\nu_2-\zeta_2\vep_1},
\end{multline*}
which reads \eqref{eq:oest1}.
\end{proof}

\begin{lemma}\label{lemma:oest2}
We have
\begin{equation}\label{eq:oest2}
\left\|u_{\mx}^{1-p} \mcj_1\left[\psi^{\out}, \hps^{\tin}\right]\right\|_{*,\ar;\sigma_0}
\le C_{32} \mu_0^{\delta_4}(t_0) \left\|\psi^{\tin}\right\|_{\sharp',a,b;\sigma_0\(\mcb_{\bmx}\)}
\end{equation}
where $C_{32} > 0$ depends only on $(M,g_0)$, $N$, $h$, $z_0$, $\alpha$, and $\sigma_0$, and $\delta_4 > 0$ is a number satisfying \eqref{eq:delta4}.
\end{lemma}
\begin{proof}
It suffices to consider only points in the set $\{x \in \R^N: |x| \le 2\mu_0^{\vep_1}\} \supset \text{supp}(\eta_{\mu_0^{\vep_1}})$.

In view of \eqref{eq:u1}--\eqref{eq:v1}, \eqref{eq:u2}--\eqref{eq:v2}, \eqref{eq:P_20}, and \eqref{eq:Qdecay}, we have
\begin{align*}
\left|u_{\mx}^{p-1} - \(u_{\mx}^{(1)}\)^{p-1}\right|(x,t) &= (1+P(x))^{p-1} \mu^{-2} \left|(W_{1,0} + \Psi_0)^{p-1} - W_{1,0}^{p-1}\right|(y,t) \\
&\le C \mu_0^{-2} \left[\mbone_{N=5}\, \bmu^2 \(W_{1,0}^{p-2} Q_0\)(y) + \(\bmu^2Q_0(y)\)^{p-1}\right] \\
&\le C \left[\mbone_{N=5}\, \log(2+|y|) + \mu_0^{-2+\frac{8}{N-2}} \log(2+|y|)^{p-1}\right] \frac{1}{1+|y|^4}.
\end{align*}
Besides, by \eqref{eq:psi00},
\begin{equation}\label{eq:oest21}
\left|\hps^{\, \tin}(x,t)\right| \le C \bmu^{-{N-2 \over 2}} \left|\psi^{\tin}(\by,t)\right|
\le C \mu_0^{-{N-2 \over 2}} \left\|\psi^{\tin}\right\|_{\sharp',a,b\(\mcb_{\bmx}\)} {\mu_0^b \over 1+|\by|^a}
\end{equation}
and
\[\left|\hps^{\, \tin}_t(x,t)\right| \le C \mu_0^{-{N-2 \over 2}} \left\|\psi^{\tin}\right\|_{\sharp',a,b\(\mcb_{\bmx}\)} {\mu_0^b \over 1+|\by|^{a-2}}\]
where $x = \mu y+\xi = \bmu\by+\xi$. Accordingly,
\begin{equation}\label{eq:oest22}
\begin{aligned}
&\ \left\|p u_{\mx}^{1-p} \eta_{\mu_0^{\vep_1}} \left[u_{\mx}^{p-1} - \(u_{\mx}^{(1)}\)^{p-1}\right]
\left[\frac{(N+2)\kappa_N}{4} \left|\hps^{\, \tin}\right| + \left|\hps^{\, \tin}_t\right|\right] \right\|_{*,\ar} \\
&\le C \left[\mbone_{N=5}\, \mu_0^2|\log\mu_0| + (\mu_0^2|\log\mu_0|)^{p-1}\right]
\left\|\psi^{\tin}\right\|_{\sharp',a,b\(\mcb_{\bmx}\)} \sup_{(\by,t) \in \mcb_{\bmx} \times [t_0,\infty)} {\mu_0^{b-\beta} \over 1+|\by|^{a-\alpha}} \\
&\le C \left[\mbone_{N=5}\, \mu_0^2|\log\mu_0| + (\mu_0^2|\log\mu_0|)^{p-1}\right](t_0) \left\|\psi^{\tin}\right\|_{\sharp',a,b\(\mcb_{\bmx}\)}
\end{aligned}
\end{equation}
where the second inequality comes from the assumption that $\alpha \in (0,a)$ and $\beta \in (0,b]$.

On the other hand, we see
\begin{align*}
&\ \left|\(\Delta_{g_0} \eta_{\mu_0^{\vep_1}}\) \hps^{\, \tin}\right|(x,t)
+ \left[\left|\nabla_{g_0} \eta_{\mu_0^{\vep_1}}\right| \left|\nabla_{g_0} \hps^{\, \tin}\right|\right](x,t) \\
&\le C \mu_0^{-{N-2 \over 2}} \left\|\psi^{\tin}\right\|_{\sharp',a,b\(\mcb_{\bmx}\)}
\(\mu_0^{-2\vep_1} {\mu_0^b \over 1+|\by|^a} + \mu_0^{-\vep_1} {\mu_0^b \over 1+|\by|^{a+1}}\) \\
&\le C \mu_0^{-{N+2 \over 2}} \left\|\psi^{\tin}\right\|_{\sharp',a,b\(\mcb_{\bmx}\)} {\mu_0^b \over 1+|\by|^{a+2}}.
\end{align*}
Therefore,
\begin{equation}\label{eq:oest23}
\begin{aligned}
&\ \left\|u_{\mx}^{1-p} \cdot \frac{(N+2)\kappa_N}{4} \left[\(\Delta_{g_0} \eta_{\mu_0^{\vep_1}}\) \hps^{\, \tin}
+ 2 \la \nabla_{g_0} \eta_{\mu_0^{\vep_1}}, \nabla_{g_0} \hps^{\, \tin} \ra_{g_0} \right]\right\|_{*,\ar} \\
&\le C \left\|\psi^{\tin}\right\|_{\sharp',a,b\(\mcb_{\bmx}\)}
\sup_{\left\{(\by,t): \bmu^{-1}\mu_0^{\vep_1} \le |\by+\bmu^{-1}\xi| \le 2\bmu^{-1}\mu_0^{\vep_1},\, t \ge t_0\right\}} {\mu_0^{b-\beta} \over 1+|\by|^{a-\alpha}} \\
&\le C \mu_0^{(1-\vep_1)(a-\alpha)}(t_0) \left\|\psi^{\tin}\right\|_{\sharp',a,b\(\mcb_{\bmx}\)}.
\end{aligned}
\end{equation}
Combining \eqref{eq:oest22} and \eqref{eq:oest23} gives
\[\left\|u_{\mx}^{1-p} \mcj_1\left[\psi^{\out}, \hps^{\tin}\right]\right\|_{*,\ar} \le C \mu_0^{\delta_4}(t_0) \left\|\psi^{\tin}\right\|_{\sharp',a,b\(\mcb_{\bmx}\)}.\]

\medskip
Further computations yield a weighted H\"older bound of $\mcj_1$ and so \eqref{eq:oest2}. We omit the details.
\end{proof}

\begin{lemma}
We have
\begin{equation}\label{eq:oest3}
\begin{aligned}
\begin{medsize}
\displaystyle \left\|u_{\mx}^{1-p} \mcj_2\left[\psi^{\out}, \hps^{\tin}\right]\right\|_{*,\ar;\sigma_0}
\end{medsize}
&\begin{medsize}
\displaystyle \le C_{33} \delta_0^{-\zeta_{33}} \left[\mbone_{N=5} \left\{\(\left\|\psi^{\out}\right\|_{*',\ar;\sigma_0}
+ \left\|\psi^{\tin}\right\|_{\sharp',a,b;\sigma_0\(\mcb_{\bmx}\)}\) \mu_0^{4-3\vep_1}(t_0) \right. \right.
\end{medsize} \\
&\hspace{75pt} \begin{medsize}
\displaystyle \left. + \(\left\|\psi^{\out}\right\|_{*',\ar;\sigma_0}^2 + \left\|\psi^{\tin}\right\|_{\sharp',a,b;\sigma_0\(\mcb_{\bmx}\)}^2 \mu_0^{\delta_3}(t_0)\)\right\}
\end{medsize} \\
&\hspace{50pt} \begin{medsize}
\displaystyle + \left\{\(\left\|\psi^{\out}\right\|_{*',\ar;\sigma_0}^{p-1}
+ \left\|\psi^{\tin}\right\|_{\sharp',a,b;\sigma_0\(\mcb_{\bmx}\)}^{p-1}\) \mu_0^{4-(N-2)\vep_1}(t_0) \right.
\end{medsize} \\
&\hspace{75pt} \begin{medsize}
\displaystyle \left. \left.+ \(\left\|\psi^{\out}\right\|_{*',\ar;\sigma_0}^p
+ \left\|\psi^{\tin}\right\|_{\sharp',a,b;\sigma_0\(\mcb_{\bmx}\)}^p \mu_0^{\delta_3}(t_0)\)\right\} \right]
\end{medsize}
\end{aligned}
\end{equation}
where $C_{33},\, \zeta_{33} > 0$ depends only on $(M,g_0)$, $N$, $h$, $z_0$, $\alpha$, and $\sigma_0$, and $\delta_3 > 0$ is the number satisfying \eqref{eq:delta3}.
\end{lemma}
\begin{proof}
By \eqref{eq:mcj2}, \eqref{eq:psimxo}, and \eqref{eq:ArBe},
\[\mcj_2\left[\psi^{\out}, \hps^{\, \tin}\right] = A_1 + A_2 \quad \text{on } M \times [t_0,\infty)\]
where
\[|A_1| \le C \left[\mbone_{N=5} u_{\mx}^{p-2} \left\{\left|\psi^{\out}\right|^2 + \eta_{\mu_0^{\vep_1}}^2 \left|\hps^{\, \tin}\right|^2\right\}
+ \left\{\left|\psi^{\out}\right|^p + \eta_{\mu_0^{\vep_1}}^p \left|\hps^{\, \tin}\right|^p\right\}\right]\]
and
\begin{multline*}
|A_2| \le C \left[\mbone_{N=5} u_{\mx}^{p-2} \left\{\left|\psi^{\out}\right| + \eta_{\mu_0^{\vep_1}} \left|\hps^{\, \tin}\right|\right\}
+ \left\{\left|\psi^{\out}\right|^{p-1} + \eta_{\mu_0^{\vep_1}}^{p-1} \left|\hps^{\, \tin}\right|^{p-1}\right\}\right] \\
\times \left[\mu_0^{{N+2 \over 2}-(N-2)\vep_1} + \left|\psi_t^{\out}\right| + \eta_{\mu_0^{\vep_1}} \left|\hps_t^{\, \tin}\right|\right].
\end{multline*}

Let us estimate the term $A_1$. We have
\begin{align*}
&\ \(\mu_0^{\rho} w_{\alpha,2}\)^{-1} \left[\mbone_{N=5} u_{\mx}^{p-2} \left|\psi^{\out}\right|^2 + \left|\psi^{\out}\right|^p\right] \\
&\le C \left[\mbone_{N=5} \left\|\psi^{\out}\right\|_{*',\ar}^2 \mu_0^{\rho} u_{\mx}^{p-2} \(w_{\alpha,2}^{-1} w_{\alpha,0}^2\)
+ \left\|\psi^{\out}\right\|_{*',\ar}^p \mu_0^{\rho(p-1)} \(w_{\alpha,2}^{-1} w_{\alpha,0}^p\)\right] \\
&\le C \delta_0^{-\zeta} \left[\mbone_{N=5} \left\|\psi^{\out}\right\|_{*',\ar}^2 \(\mu_0^{\beta} + \mu_0^2\)
+ \left\|\psi^{\out}\right\|_{*',\ar}^p \(\mu_0^{\beta(p-1)} + \mu_0^2\)\right]
\end{align*}
on $M \times [t_0,\infty)$. In addition, for $x = \mu y+\xi \in \text{supp}(\eta_{\mu_0^{\vep_1}})$, it holds that
\begin{align*}
&\ \(\mu_0^{\rho} w_{\alpha,2}\)^{-1} \left[\mbone_{N=5} u_{\mx}^{p-2} \eta_{\mu_0^{\vep_1}}^2 \left|\hps^{\, \tin}\right|^2
+ \eta_{\mu_0^{\vep_1}}^p \left|\hps^{\, \tin}\right|^p \right](x,t) \\
&\le C \left[\mbone_{N=5} \left\|\psi^{\tin}\right\|_{\sharp',a,b\(\mcb_{\bmx}\)}^2 \mu_0^{2b+\alpha-3} \(1+|y|^{\alpha+1-2a}\)
+ \left\|\psi^{\tin}\right\|_{\sharp',a,b\(\mcb_{\bmx}\)}^p \mu_0^{bp-\beta} \(1+|y|^{\alpha+2-ap}\)\right] \\
&\le C \left[\mbone_{N=5} \left\|\psi^{\tin}\right\|_{\sharp',a,b\(\mcb_{\bmx}\)}^2 \(\mu_0^{2b-\beta} + \mu_0^{2(a+b-2)}\)
+ \left\|\psi^{\tin}\right\|_{\sharp',a,b\(\mcb_{\bmx}\)}^p \(\mu_0^{bp-\beta} + \mu_0^{(a+b)p-N}\)\right]
\end{align*}
where we applied \eqref{eq:oest21} for the first inequality.

Arguing as above, we can also estimate
\begin{align*}
&\, \begin{medsize}
\(\mu_0^{\rho} w_{\alpha,2}\)^{-1} |A_2|
\end{medsize} \\
&\begin{medsize}
\le C \delta_0^{-\zeta} \left[\mbone_{N=5} \left\{\(\left\|\psi^{\out}\right\|_{*',\ar} + \left\|\psi^{\tin}\right\|_{\sharp',a,b\(\mcb_{\bmx}\)}\) \mu_0^{4-3\vep_1}
+ \(\left\|\psi^{\out}\right\|_{*',\ar}^2 + \left\|\psi^{\tin}\right\|_{\sharp',a,b\(\mcb_{\bmx}\)}^2 \mu_0^{\delta_3}(t_0)\)\right\} \right.
\end{medsize} \\
&\begin{medsize}
\hspace{35pt} \left. + \left\{\(\left\|\psi^{\out}\right\|_{*',\ar}^{p-1} + \left\|\psi^{\tin}\right\|_{\sharp',a,b\(\mcb_{\bmx}\)}^{p-1}\) \mu_0^{4-(N-2)\vep_1}
+ \(\left\|\psi^{\out}\right\|_{*',\ar}^p + \left\|\psi^{\tin}\right\|_{\sharp',a,b\(\mcb_{\bmx}\)}^p \mu_0^{\delta_3}(t_0)\)\right\} \right]
\end{medsize}
\end{align*}
on $M \times [t_0,\infty)$. To deduce the inequality, we repeatedly use the relations $\alpha+\beta = N-2$, $\alpha \in (0,a)$, and $\beta \in (0,b]$.

Combining all the computations together, we establish a weighted $L^{\infty}$-estimate for $\mcj_2$,
that is, \eqref{eq:oest3} where the weighted H\"older norms are replaced with the associated weighted $L^{\infty}$-norms.

\medskip
Further computations yield a weighted H\"older bound of $\mcj_2$ and so \eqref{eq:oest3}. We omit the details.
\end{proof}

\begin{proof}[Completion of the proof of Proposition \ref{prop:outer3}]
For the parameters $(\mx)$ satisfying \eqref{eq:mu0}--\eqref{eq:lx} and the function $\psi_0$ satisfying \eqref{eq:outer31},
let $\mct^{\out}_{\mx}[f,\psi_0]$ be a unique solution to \eqref{eq:inhom-o} obtained in Corollary \ref{cor:outer1}.
A function $\psi^{\out}$ is a solution to \eqref{eq:outer3} if and only if it satisfies
\begin{equation}\label{eq:outer34}
\psi^{\out} = \mct^{\out}_{\mx}\left[f^{\out}_{\mx}\left[\psi^{\out}, \hps^{\, \tin}\right],\psi_0\right]
\end{equation}
where $f^{\out}_{\mx}$ is the map given in \eqref{eq:fout}. We shall prove the existence of a fixed point of
the operator $\mct^{\out}_{\mx}[f^{\out}_{\mx}[\cdot, \hps^{\, \tin}],\psi_0]$
by applying the contraction mapping theorem on the set
\begin{equation}\label{eq:outer33}
\mcd^{\out} := \left\{\psi^{\out}: \|\psi^{\out}\|_{*',\ar;\sigma_0}
\le 2C_1(C_{21}C_{32}+ C_{22}) \delta_0^{-2} \mu_0^{\delta_4}(t_0)\right\}
\end{equation}
where $C_1,\, C_{21},\, C_{22},\, C_{32},\, C_{33} > 0$ are the numbers appearing in \eqref{eq:outer2a0}, \eqref{eq:outer31}, \eqref{eq:oest2}, and \eqref{eq:oest3}.

By employing \eqref{eq:outer2a0}, \eqref{eq:mcj1}, \eqref{eq:mcj2}, \eqref{eq:oest1}, \eqref{eq:oest2}, \eqref{eq:oest3}, \eqref{eq:outer31},
and \eqref{eq:delta4}, and taking a larger $t_0 > 0$ if required, we obtain
\begin{align*}
&\ \left\|\mct^{\out}_{\mx}\left[f^{\out}_{\mx}\left[\psi^{\out}, \hps^{\, \tin}\right],\psi_0\right]\right\|_{*',\ar;\sigma_0} \\
&\le C_1\delta_0^{-2} \(\left\|f^{\out}_{\mx}\left[\psi^{\out}, \hps^{\, \tin}\right]\right\|_{*,\ar;\sigma_0} + \|\psi_0\|_{**,\alpha;\sigma_0}\) \\
&\le C_1\delta_0^{-2} \left[C_{31} \delta_0^{-\zeta_{31}} \mu_0^{\nu_2-\zeta_2\vep_1}(t_0)
+ C_{21} C_{32} \mu_0^{\delta_4}(t_0) + C_{22} \mu_0^{\delta_2}(t_0) \right. \\
&\hspace{45pt} + C_{33} \delta_0^{-\zeta_{33}}
\left\{\mbone_{N=5} \(2C_{21} \mu_0^{4-3\vep_1}(t_0) + \left\|\psi^{\out}\right\|_{*',\ar;\sigma_0}^2 + C_{21}^2 \mu_0^{\delta_3}(t_0) \) \right. \\
&\hspace{100pt} \left. \left. + \(2C_{21}^{p-1} \mu_0^{4-(N-2)\vep_1}(t_0) + \left\|\psi^{\out}\right\|_{*',\ar;\sigma_0}^p
+ C_{21}^p \mu_0^{\delta_3}(t_0)\) \right\}\right] \\
&\le 2C_1 (C_{21} C_{32}+ C_{22}) \delta_0^{-2} \mu_0^{\delta_4}(t_0)
\end{align*}
for any $\psi \in \mcd^{\out}$. Furthermore,
\begin{equation}\label{eq:outer37}
\begin{aligned}
&\ \left\|\mct^{\out}_{\mx}\left[f^{\out}_{\mx}\left[\psi^{\out}_1, \hps^{\, \tin}\right],\psi_0\right]
- \mct^{\out}_{\mx}\left[f^{\out}_{\mx}\left[\psi^{\out}_2, \hps^{\, \tin}\right],\psi_0\right]\right\|_{*',\ar;\sigma_0} \\
&\le C_1 \delta_0^{-2} \left\|\mcj_2\left[\psi^{\out}_1, \hps^{\, \tin}\right] - \mcj_2\left[\psi^{\out}_2, \hps^{\, \tin}\right]\right\|_{*,\ar;\sigma_0} \\
&\le C \delta_0^{-2} \left[\left\|\(u_{\mx}+\psi^{\out}_1\)^p_+ - \(u_{\mx}+\psi^{\out}_2\)^p_+
- pu_{\mx}^{p-1}\(\psi^{\out}_1-\psi^{\out}_2\)\right\|_{*,\ar;\sigma_0} \right. \\
&\hspace{40pt} + \left\|\left\{\(u_{\mx}+\psi^{\out}_1\)^{p-1}_+ - \(u_{\mx}+\psi^{\out}_2\)^{p-1}_+\right\}
\left\{\(1-\eta_{\mu_0^{\vep_1}}\) (u_{\mx})_t + \(\psi^{\out}_1\)_t\right\}\right\|_{*,\ar;\sigma_0} \\
&\hspace{40pt} \left. + \left\|\left\{\(u_{\mx}+\psi^{\out}_2\)^{p-1}_+ - (u_{\mx})^{p-1}_+\right\} \(\psi^{\out}_1-\psi^{\out}_2\)_t\right\|_{*,\ar;\sigma_0} \right] \\
&\le C \delta_0^{-\zeta} \left[\mbone_{N=5} \mu_0^{\delta_4}(t_0) + \mu_0^{\delta_4(p-1)}(t_0)\right] \left\|\psi^{\out}_1 - \psi^{\out}_2\right\|_{*',\ar;\sigma_0}
\end{aligned}
\end{equation}
for any $\psi_1, \psi_2 \in \mcd^{\out}$.
Therefore, the operator $\mct^{\out}_{\mx}[f^{\out}_{\mx}[\cdot, \hps^{\, \tin}],\psi_0]$ is a contraction map on $\mcd^{\out}$,
and it has a fixed point in $\mcd^{\out}$, namely, a point satisfying \eqref{eq:outer34}.
Finally, \eqref{eq:outer32} follows directly from \eqref{eq:outer33}.
\end{proof}

\subsection{Derivatives of the solution to \eqref{eq:outer} with respect to parameters} \label{subsec:outer4}
In the following proposition, we provide quantitative estimates on the derivatives of $\psi^{\out}$ with respect to its parameters.

\begin{prop}
Suppose that all the conditions of Proposition \ref{prop:outer3} hold.
There exists a constant $C,\, \zeta > 0$ depending only on $(M,g_0)$, $N$, $h$, $z_0$, $\alpha$, and $\sigma_0$ such that
\begin{equation}\label{eq:outer35}
\begin{cases}
\left\|\pa_{\lambda}\psi^{\out}\left[\bar{\lambda}\right]\right\|_{*',\ar;\sigma_0}
\le C\delta_0^{-\zeta} \mu_0^{\delta_4}(t_0) \mu_0^{\nu_1-1} \left\|\bar{\lambda}\right\|_{\nu_1},\\
\left\|\pa_{\xi}\psi^{\out}\left[\bar{\xi}\right]\right\|_{*',\ar;\sigma_0}
\le C\delta_0^{-\zeta} \mu_0^{\delta_4}(t_0) \mu_0^{\nu_2-1} \left\|\bar{\xi}\right\|_{\nu_2},\\
\left\|\pa_{\dot{\lambda}}\psi^{\out}\big[\dot{\bar{\lambda}}\big]\right\|_{*',\ar;\sigma_0}
\le C\delta_0^{-\zeta} \mu_0^{\nu_1+3-\zeta \vep_1} \left\|\dot{\bar{\lambda}}\right\|_{\nu_1+2},
\\
\left\|\pa_{\dot{\xi}}\psi^{\out}\big[\dot{\bar{\xi}}\big]\right\|_{*',\ar;\sigma_0}
\le C\delta_0^{-\zeta} \mu_0^{\nu_2+4-\zeta \vep_1} \left\|\dot{\bar{\xi}}\right\|_{\nu_2+2},
\\
\left\|\pa_{\hps^{\, \tin}}\psi^{\out}\left[\bar{\psi}\right]\right\|_{*',\ar;\sigma_0}
\le C\delta_0^{-\zeta} \mu_0^{\delta_4}(t_0) \left\|\mu_0^{-{N-2 \over 2}} \bar{\psi}\right\|_{\sharp',a,b;\sigma_0}
\end{cases}
\end{equation}
and
\begin{equation}\label{eq:outer36}
\left\|\psi^{\out}[\psi_{01}] - \psi^{\out}[\psi_{02}]\right\|_{*',\ar;\sigma_0}
\le C\delta_0^{-2}\left\|\psi_{01}-\psi_{02}\right\|_{**,\alpha;\sigma_0}
\end{equation}
Here, we write
\[\pa_{\lambda}\psi^{\out}[\bar{\lambda}] := \left.\pa_s \psi^{\out} \left[\lambda+s\bar{\lambda}, \xi, \dot{\lambda}, \dot{\xi}, \hps^{\, \tin}, \psi_0\right]\right|_{s=0}, \quad
\psi^{\out}[\psi_0] := \psi^{\out}[\lambda, \xi, \dot{\lambda}, \dot{\xi}, \hps^{\, \tin}, \psi_0],\]
and so on. Also, the $\nu$- and $\sharp'$-norms were defined in \eqref{eq:nunorm} and \eqref{eq:sharp'-norm}, respectively.
\end{prop}
\begin{proof}
The implicit function theorem implies that the functions $\pa_{\lambda}\psi^{\out}[\bar{\lambda}], \ldots, \pa_{\hps^{\, \tin}}\psi^{\out}\left[\bar{\psi}\right]$
are well-defined in a neighborhood of given data $(\lambda, \xi, \dot{\lambda}, \dot{\xi}, \hps^{\, \tin}, \psi_0)$.
Examining the equations that the functions solve, we can deduce \eqref{eq:outer35}.
We skip the details, because the necessary computations are similar to those in the proof of Proposition \ref{prop:outer3};
refer to \cite[Proposition 4.2]{CdPM} where an analogous result was obtained for critical nonlinear heat equation.

Furthermore, by \eqref{eq:outer2a0} and \eqref{eq:outer37},
\begin{align*}
&\ \left\|\psi^{\out}[\psi_{01}] - \psi^{\out}[\psi_{02}]\right\|_{*',\ar;\sigma_0} \\
&\le C\delta_0^{-2}\(\left\|\mcj_2\left[\psi^{\out}[\psi_{01}], \hps^{\tin}\right] - \mcj_2\left[\psi^{\out}[\psi_{02}], \hps^{\tin}\right]\right\|_{*,\ar;\sigma_0} + \|\psi_{01} - \psi_{02}\|_{**,\alpha;\sigma_0}\) \\
&\le C\delta_0^{-2} \(o(1)\left\|\psi^{\out}[\psi_{01}] - \psi^{\out}[\psi_{02}]\right\|_{*',\ar;\sigma_0} + \|\psi_{01} - \psi_{02}\|_{**,\alpha;\sigma_0}\),
\end{align*}
which yields \eqref{eq:outer36}.
\end{proof}

\section{Inner problem and the completion of the proof of Theorem \ref{thm:main}} \label{sec:inner}
This section is devoted to the analysis of the inner problem \eqref{eq:inner}, which will allow us to complete the proof of the main theorem.

In Subsection \ref{subsec:inner1}, we establish the unique solvability and a priori estimate of the non-uniform inhomogeneous parabolic equation
\begin{equation}\label{eq:inhom-i}
\begin{cases}
\displaystyle p W_{1,0}^{p-1} \psi_t = \frac{(N+2)\kappa_N}{4} \(\Delta \psi + p W_{1,0}^{p-1} \psi\) + W_{1,0}^{p-1} \mch &\text{in } \R^N \times (t_0,\infty),\\
\psi(\cdot,t_0) = e_0(t_0) Z_0 &\text{in } \R^N
\end{cases}
\end{equation}
under the assumption that
\begin{equation}\label{eq:ortho1}
\int_{\R^N} \mch(\by,t) \(W_{1,0}^{p-1}Z_n\)(\by) d\by = 0 \quad \text{for } t \in [t_0,\infty) \text{ and } n = 1, \ldots, N+1.
\end{equation}
Here,
\begin{itemize}
\item[-] $e_0(t_0)$ is a real number to be determined by $\mch$;
\item[-] $Z_1, \ldots, Z_{N+1}$ are the functions in \eqref{eq:ZN+1} and \eqref{eq:Zi};
\item[-] $Z_0$ is the unique positive radial function in $\R^N$ such that
\[\begin{cases}
L_0[Z_0] = W_{1,0}^{1-p} \(\Delta Z_0 + pW_{1,0}^{p-1}Z_0\) = \mfm_0 Z_0 \quad \text{in } \R^N,\\
\displaystyle \mfm_0 > 0,\ \int_{\R^N} W_{1,0}^{p-1} Z_0^2 = 1.
\end{cases}\]
It is clear that $\mfm_0 = p-1$ and $Z_0$ is a constant multiple of $W_{1,0}$; refer to \cite[Appendix]{BE}.
\end{itemize}

Given parameters $(\lambda,\xi)$, and functions $\psi^{\tin}$ and $\psi_0$ satisfying \eqref{eq:mu0}--\eqref{eq:lx} and \eqref{eq:outer31}, respectively,
let $\psi^{\out} = \psi^{\out}[\lambda, \xi, \dot{\lambda}, \dot{\xi}, \hps^{\, \tin}, \psi_0]$ be the unique solution to \eqref{eq:outer3} found in Proposition \ref{prop:outer3}.
In Subsection \ref{subsec:inner2}, we reduce
\begin{equation}\label{eq:ortho2}
\int_{\R^N} \eta_{2\mu_0^{\vep_1}}(x) \left[\({\bmu \over \mu}\)^{N-2 \over 2} (1+P)^{-p}(x) \mce_2[\mx](y,t)
+ \mck_1\left[\psi^{\tin}\right] + \mck_2\left[\psi^{\tin}, \psi^{\out}\right]\right] Z_n(\by) d\by = 0
\end{equation}
for $n = 1, \ldots, N+1$ and $t \in [t_0,\infty)$ into a system of nonlinear ODEs \eqref{eq:ODExi} and \eqref{eq:ODEla} of $(\lambda,\xi)$ and solve it;
refer to \eqref{eq:eta}, \eqref{eq:E2}, \eqref{eq:psimxo}, \eqref{eq:psi0}, \eqref{eq:mck1}, \eqref{eq:mck2}, and Subsection \ref{subsec:approx1} for the definition of the notations.

In Subsection \ref{subsec:inner3}, we solve \eqref{eq:inner} by applying \eqref{eq:ortho2}, the unique solvability of \eqref{eq:inhom-i} under \eqref{eq:ortho1}, and the contraction mapping theorem.

In Subsection \ref{subsec:inner4}, we complete the proof of Theorem \ref{thm:main}.

\subsection{Linear theory} \label{subsec:inner1}
In this subsection, we examine \eqref{eq:inhom-i} provided $a \in (\sigma_0,N-2)$ and $b \simeq 3$.

\begin{prop}\label{prop:lin}
Suppose that $\|\mch\|_{\sharp,a+2,b;\sigma_0} < \infty$ and \eqref{eq:ortho1} holds.
Then one can find $\psi = \psi[\mch]$ and $e_0 = e_0[\mch] \in \R$ satisfying \eqref{eq:inhom-i}.
These $\psi$ and $e_0$ are linear in $\mch$, and
\begin{equation}\label{eq:lin0}
\|\psi\|_{\sharp',a,b;\sigma_0(\R^N)} + \|e_0\|_{b;\sigma_0} \le C_4 \|\mch\|_{\sharp,a+2,b;\sigma_0}
\end{equation}
for a constant $C_4 > 0$ depending only on $N$, $a$, $b$ and $\sigma_0$. Here, the norms of $\psi$, $e_0$, and $\mch$ are given in Definitions \ref{defn:norm3} and \ref{defn:norm0}.
\end{prop}
\noindent As an intermediate step to proving the proposition, we will analyze a uniformly parabolic equation on $\S^N$
\begin{equation}\label{eq:inhom-i2}
\begin{cases}
\displaystyle p\, \vph_t = \frac{\kappa_N p}{N} \(\Delta_{\S^N} \vph + N\vph\) + \wtmch - \tc_0(t) &\text{on } \S^N \times (t_0,s_0), \\
\vph(\ty,t_0) = 0 &\text{on } \S^N
\end{cases}
\end{equation}
for $s_0 > \frac{3t_0}{2}$ large enough. Here, $\wtmch$ is a function satisfying
\[\left\|\wtmch\right\|_{\wtsh,a,b;t_0,s_0} := \left\|\mu_0^{-b} [d_{\,\S^N}(\ty)]^{N-a} \wtmch(\ty,t) \right\|_{L^{\infty}(\S^N \times [t_0,s_0))} < \infty,\]
where $d_{\,\S^N}(\ty) := \arccos(\ty_{N+1})$ is the geodesic distance on $\S^N$ between $\tyn := (0,\ldots,0,1)$ and $\ty$, and
\begin{equation}\label{eq:ortho11}
\int_{\S^N} \wtmch(\ty,t) \ty_n dS_{\ty} = 0 \quad \text{for } t \in [t_0,s_0) \text{ and } n = 1, \ldots, N+1.
\end{equation}
Also, $\tc_0$ is defined as
\begin{equation}\label{eq:tc}
\tc_0(t) = \frac{1}{\left|\S^N\right|} \int_{\S^N} \wtmch(\ty,t) dS_{\ty} \quad \text{for } t \in [t_0,s_0).
\end{equation}
By integrating the equation in \eqref{eq:inhom-i2} over $\S^N$, and applying \eqref{eq:tc} and the initial condition on $\vph$, we find
\begin{equation}\label{eq:ortho120}
\int_{\S^N} \vph(\ty,t) dS_{\ty} = 0 \quad \text{for } t \in [t_0,s_0).
\end{equation}
Moreover, testing $\ty_n$ on \eqref{eq:inhom-i2} for $n = 1, \ldots, N+1$, and exploiting \eqref{eq:ortho11}
and the fact that $N$ is the second eigenvalue of $-\Delta_{\S^N}$ with eigenfunctions $\ty_1, \ldots, \ty_{N+1}$, we obtain
\begin{equation}\label{eq:ortho12}
\int_{\S^N} \vph(\ty,t) \ty_n dS_{\ty} = 0 \quad \text{for } t \in [t_0,s_0) \text{ and } n = 1, \ldots, N+1.
\end{equation}

\begin{lemma}\label{lemma:lin}
Assume that $\left\|\wtmch\right\|_{\wtsh,a,b;t_0,s_0} < \infty$ for $s_0 > \frac{3t_0}{2}$ large, and \eqref{eq:ortho11} holds.
If $\vph$ is a solution to \eqref{eq:inhom-i2} such that $\vph,\, D_{\ty}\vph,\, D^2_{\ty}\vph,\, \vph_t \in L^2(\S^N \times (t_0,s_0))$,
then there exists a constant $C > 0$ depending only on $N$, $a$ and $b$ such that
\begin{equation}\label{eq:lin01}
\|\vph\|_{\wtsh,a+2,b;t_0,s_0} + \sup_{t \in [t_0,s_0)} \mu_0^{-b}(t) |\tc_0(t)| \le C\left\|\wtmch\right\|_{\wtsh,a,b;t_0,s_0}.
\end{equation}
\end{lemma}
\begin{proof}
To deduce the lemma, we will argue as in the proof of \cite[Lemma 4.1]{SWZ}; see Remark \ref{rmk:main} (3). The proof is divided into two steps.

Throughout the proof, $C$ denotes a universal constant independent of $s_0$.

\medskip \noindent \textsc{Step 1.} We insist that
\begin{equation}\label{eq:lin02}
\sup_{t \in [t_0,s_0)} \mu_0^{-b}(t) |\tc_0(t)| \le C \left\|\wtmch\right\|_{\wtsh,a,b;t_0,s_0} \quad \text{and} \quad \left\|\vph\right\|_{\wtsh,a+2,b;t_0,s_0} < \infty.
\end{equation}

The first inequality in \eqref{eq:lin02} is obvious. Let us consider the second one.
If $\delta > 0$ is any small number, parabolic regularity theory ensures the existence of $K_1 > 0$ depending on $\vph$ and $\wtmch$ such that
\[\|\vph\|_{C^0\(\(\S^N \setminus B_{\S^N}(\tyn,\delta)\) \times [t_0,s_0)\)} \le K_1.\]
Also, the geodesic distance $d_{\,\S^N}(\ty) = \arccos(\ty_{N+1})$ on $\S^N$ between $\tyn$ and $\ty$ satisfies
\[\Delta_{\S^N} [d_{\,\S^N}(\ty)]^{a-(N-2)} 
= ((N-2)-a)\(-a+O(\delta^2)\) [d_{\,\S^N}(\ty)]^{a-N} \quad \text{in } B_{\S^N}(\tyn,\delta).\]
Therefore, if we set
\[v(\ty,t) = K_2 \left\|\wtmch\right\|_{\wtsh,a,b;t_0,s_0} \mu_0^b(t) [d_{\,\S^N}(\ty)]^{a-(N-2)} \quad \text{in } B_{\S^N}(\tyn,\delta) \times (t_0,s_0),\]
then
\[\begin{cases}
\displaystyle (v \pm \vph)_t \ge \frac{\kappa_N}{N} \left[\Delta_{\S^N}(v \pm \vph) + N(v \pm \vph)\right] &\text{in } B_{\S^N}(\tyn,\delta) \times (t_0,s_0), \\
(v \pm \vph)(\ty,t_0) = v(\ty,t_0) \ge 0 &\text{in } B_{\S^N}(\tyn,\delta), \\
v \pm \vph \ge 0 &\text{on } \pa B_{\S^N}(\tyn,\delta) \times (t_0,s_0)
\end{cases}\]
provided $K_2 \gg K_1$. The parabolic maximum principle implies the second inequality in \eqref{eq:lin02}.

\medskip \noindent \textsc{Step 2.} We assert that
\begin{equation}\label{eq:lin03}
\left\|\vph\right\|_{\wtsh,a+2,b;t_0,s_0} \le C \left\|\wtmch\right\|_{\wtsh,a,b;t_0,s_0}.
\end{equation}
To the contrary, suppose that there exist increasing sequences $\{t_{\ell}\}_{\ell \in \N}$, $\{s_{\ell}\}_{\ell \in \N}$ of positive numbers for which \eqref{eq:outer1a11} holds,
and sequences $\{\vph_{\ell}\}_{\ell \in \N}$, 
$\{\wtmch_{\ell}\}_{\ell \in \N}$ of functions which satisfy
\begin{equation}\label{eq:lin04}
\begin{cases}
\displaystyle p\, (\vph_{\ell})_t = \frac{\kappa_N p}{N} \(\Delta_{\S^N} \vph_{\ell} + N\vph_{\ell}\) + \wtmch_{\ell} - \tc_{\ell}(t) &\text{on } \S^N \times (t_{\ell},s_{\ell}), \\
\vph_{\ell}(\ty,t_{\ell}) = 0 &\text{on } \S^N,
\end{cases}
\end{equation}
\begin{equation}\label{eq:ortho11l}
\int_{\S^N} \wtmch_{\ell}(\ty,t) \ty_n dS_{\ty} = 0 \quad \text{for } t \in [t_{\ell},s_{\ell}) \text{ and } n = 1, \ldots, N+1,
\end{equation}
and
\begin{equation}\label{eq:lin05}
\left\|\vph_{\ell}\right\|_{\wtsh,a+2,b;t_{\ell},s_{\ell}} = 1, \quad \left\|\wtmch_{\ell}\right\|_{\wtsh,a,b;t_{\ell},s_{\ell}} \to 0 \quad \text{as } \ell \to \infty.
\end{equation}
Here, $\tc_{\ell}$ is defined as
\begin{equation}\label{eq:tcl}
\tc_{\ell}(t) := \frac{1}{\left|\S^N\right|} \int_{\S^N} \wtmch_{\ell}(\ty,t) dS_{\ty} \quad \text{for } t \in [t_{\ell},s_{\ell}).
\end{equation}
By \eqref{eq:ortho120}, \eqref{eq:ortho12} and \eqref{eq:lin05}, we have
\begin{equation}\label{eq:ortho12l}
\int_{\S^N} \vph_{\ell}(\ty,t) dS_{\ty} = \int_{\S^N} \vph_{\ell}(\ty,t) \ty_n dS_{\ty} = 0
\quad \text{for } t \in [t_{\ell},s_{\ell}) \text{ and } n = 1, \ldots, N+1,
\end{equation}
and there is a sequence $\{(\ty_{\ell},\tau_{\ell})\}_{\ell \in \N}$ such that $\{\tau_{\ell}\}_{\ell \in \N}$ is increasing and
\begin{equation}\label{eq:lin06}
\begin{cases}
\displaystyle \ty_{\ell} = (\ty_{\ell 1}, \ldots, \ty_{\ell (N+1)}) \in \S^N,
\quad \lim_{\ell \to \infty} \ty_{\ell} = \ty_{\infty} \in \S^N, \quad \tau_{\ell} \in (t_{\ell},s_{\ell}), \\
\displaystyle \frac{1}{2} \le \mu_0^{-b}(\tau_{\ell}) [d_{\,\S^N}(\ty_{\ell})]^{N-2-a} |\vph_{\ell}(\ty_{\ell},\tau_{\ell})| \le 1 \quad \text{for all } \ell \in \N.
\end{cases}
\end{equation}

Let us check that $\ty_{\infty} = \tyn$. If not, we would have a small number $\delta > 0$ such that $d_{\,\S^N}(\ty_{\ell}) \ge \delta$ for all large $\ell \in \N$.
Let $\chi: \S^N \to [0,1]$ be a smooth function such that $\chi(\ty) = 1$ if $d_{\,\S^N}(\ty) \ge \frac{\delta}{2}$, $\chi(\ty) = 0$
if $d_{\,\S^N}(\ty) \le \frac{\delta}{4}$, and $|\nabla \chi| \le \frac{4}{\delta}$ on $\S^N$.
We test \eqref{eq:lin04} against $\chi^2 \vph_{\ell}$ for each $\tau \in (t_{\ell},s_{\ell})$, integrate the result over $(t_{\ell},\tau)$,
and apply \eqref{eq:lin02}, the initial condition on $\vph_{\ell}$, \eqref{eq:lin05}, and Gr\"onwall's inequality. Then we obtain
\[\int_{\left\{d_{\,\S^N}(\ty) \ge \frac{\delta}{2}\right\}} \vph_{\ell}^2(\ty,\tau) dS_{\ty}
\le o(1) \left[\int_{t_{\ell}}^{\tau} \mu_0^{2b}(t)dt\right] e^{C(\tau-t_{\ell})} \quad \text{for } \tau \in [t_{\ell},s_{\ell}).\]
Combining this with parabolic estimates yields that $\liminf_{\ell \to \infty}(\tau_{\ell} - t_{\ell}) = \infty$, as otherwise
\begin{align*}
&\ \mu_0^{-b}(\tau_{\ell}) \|\vph_{\ell}\|_{C^0(\{d_{\,\S^N}(\ty) \ge \delta\} \times [\tau_{\ell}-1,\tau_{\ell}])} \\
&\le C\mu_0^{-b}(\tau_{\ell}) \left[\|\vph_{\ell}\|_{L^2(\{d_{\,\S^N}(\ty) \ge \frac{\delta}{2}\} \times [\tau_{\ell}-2,\tau_{\ell}])}
+ \left\|\wtmch_{\ell}\right\|_{L^{\infty}(\{d_{\,\S^N}(\ty) \ge \frac{\delta}{2}\} \times [\tau_{\ell}-2,\tau_{\ell}])}
+ \|\tc_{\ell}\|_{L^{\infty}([\tau_{\ell}-2,\tau_{\ell}])} \right] \\
&\le o(1) \left[\left\{\int_{\tau_{\ell}-2}^{\tau_{\ell}} \int_{t_{\ell}}^{\tau} \(\frac{\tau_{\ell}}{t}\)^b e^{C(\tau-t_{\ell})} dt d\tau \right\}^{1 \over 2} + 1\right] = o(1),
\end{align*}
contradicting \eqref{eq:lin06}. It follows that the limit equation of the sequence $\{\mu_0^{-b}(t+\tau_{\ell})\, \vph_{\ell}(\ty,t+\tau_{\ell})\}_{\ell \in \N}$ is
\begin{equation}\label{eq:lin07}
\begin{cases}
\displaystyle (\bvp_{\infty})_t = \frac{\kappa_N}{N} \(\Delta_{\S^N} \bvp_{\infty} + N\bvp_{\infty}\) &\text{on } \S^N \times (-\infty,0), \\
\displaystyle |\bvp_{\infty}(\ty,t)| \le [d_{\,\S^N}(\ty)]^{a-(N-2)} &\text{for } (\ty,t) \in \S^N \times (-\infty,0),
\end{cases}
\end{equation}
and it holds that $\bvp_{\infty}(\ty_{\infty},0) \ne 0$ and
\begin{equation}\label{eq:ortho12inf}
\int_{\S^N} \bvp_{\infty}(\ty,t) dS_{\ty} = \int_{\S^N} \bvp_{\infty}(\ty,t) \ty_n dS_{\ty} = 0
\quad \text{for } t \in (-\infty,0) \text{ and } n = 1, \ldots, N+1.
\end{equation}
By the growth condition on $\bvp_{\infty}$ in \eqref{eq:lin07} and parabolic regularity,
$\bvp_{\infty}$ is in fact smooth on $\S^N \times (-\infty,0)$ and $\|\bvp_{\infty}\|_{C^{\infty}(\S^N \times (-\infty,0))} \le C$. Also, \eqref{eq:ortho12inf} gives
\[\mfB(\bvp_{\infty}(\cdot,t)) \ge 0 \quad \text{for } t \in (-\infty,0) \quad \text{where } \mfB(\bvp) := \int_{\S^N} \(|\nabla_{\S^N} \bvp|^2 - N\bvp^2\) dS_{\ty}.\]
Testing $(\bvp_{\infty})_t$ on \eqref{eq:lin07} and the equation of $(\bvp_{\infty})_t$, respectively, we get
\[\int_{-t}^0 \int_{\S^N} (\bvp_{\infty})_t^2 = \frac{\kappa_N}{2N} \left[\mfB(\bvp_{\infty}(\cdot,-t)) - \mfB(\bvp_{\infty}(\cdot,0))\right]
\le C\|\bvp_{\infty}\|_{C^1((\S^N \times (-\infty,0)))} \le C\]
and
\[\pa_t \int_{\S^N} (\bvp_{\infty})_t^2 = - \mathfrak{B}(\bvp_{\infty}(\cdot,t)) \le 0\]
for all $t \in (-\infty,0]$. Hence $(\bvp_{\infty})_t = 0$ on $\S^N \times (-\infty,0]$.
By plugging it into \eqref{eq:lin07} and employing \eqref{eq:ortho12inf} again, we conclude that $\bvp_{\infty} = 0$ on $\S^N \times (-\infty,0]$.
This is a contradiction, and we must have that $\ty_{\infty} = \tyn$.
Especially, $\nu_{\ell} := d_{\,\S^N}(\ty_{\ell}) \to 0$ and $(\ty_{\ell 1}, \ldots, \ty_{\ell N}) \to 0 \in \R^N$ as $\ell \to \infty$.

Now, we regard $\vph_{\ell}$ as the function in $B^N(0,\frac{1}{2}) \times (t_{\ell},s_{\ell})$, abusing the notation
\[\vph_{\ell}(y,t) = \vph_{\ell}\(\(y,\sqrt{1-|y|^2}\),t\) \quad
\begin{cases}
\text{for } (y,t) \in B^N\(0,\frac{1}{2}\) \times (t_{\ell},s_{\ell}) \\
\text{so that } \(\(y,\sqrt{1-|y|^2}\),t\) \in \S^N \times (t_{\ell},s_{\ell}),
\end{cases}\]
and define
\begin{multline*}
\tvp_{\ell}(Y,\tau) = \mu_0^{-b}(\tau_{\ell}) \nu_{\ell}^{N-2-a} \vph_{\ell}\(\ty_{\ell 1}+\nu_{\ell} Y_1, \ldots, \ty_{\ell N}+\nu_{\ell} Y_N, \tau_{\ell} + \nu_{\ell}^2 \tau\) \\
\text{for } Y = (Y_1,\ldots,Y_N) \in B^N\(0,\frac{1}{4\nu_{\ell}}\) \text{ and } \tau \in \(\frac{t_{\ell}-\tau_{\ell}}{\nu_{\ell}^2}, 0\right].
\end{multline*}
As $\ell \to \infty$, $\tvp_{\ell}$ tends to a function $\tvp_{\infty}$ in $\R^N \times (-\infty,0]$ which satsifies
\begin{equation}\label{eq:lin08}
\begin{cases}
\displaystyle (\tvp_{\infty})_{\tau} = \frac{\kappa_N}{N} \Delta \tvp_{\infty} \quad \text{in } \R^N \times (-\infty,0), \\
\displaystyle |\tvp_{\infty}(Y,\tau)| \le |Y-Y_0|^{a-(N-2)} \quad \text{for } (Y,\tau) \in \R^N \times (-\infty,0] \text{ and some } Y_0 \in \S^{N-1},
\end{cases}
\end{equation}
and $\frac{1}{2} \le |\tvp_{\infty}(0,0)| \le 1$.
However, a standard comparison argument shows that the only solution to \eqref{eq:lin08} is the trivial one.
This is a contradiction, and \eqref{eq:lin03} must hold.

\medskip
Inequality \eqref{eq:lin01} is an immediate consequence of \eqref{eq:lin02} and \eqref{eq:lin03}.
\end{proof}

\begin{cor}\label{cor:lin}
Assume that $\left\|\wtmch\right\|_{\wtsh,a,b} := \left\|\wtmch\right\|_{\wtsh,a,b;t_0,\infty} < \infty$,
and \eqref{eq:ortho11} and \eqref{eq:tc} hold for $s_0 = \infty$.
Then one can find a solution $\vph$ to the equation
\begin{equation}\label{eq:inhom-i21}
\begin{cases}
\displaystyle p\, \vph_t = \frac{\kappa_N p}{N} \(\Delta_{\S^N} \vph + N\vph\) + \wtmch - \tc_0(t) &\text{on } \S^N \times (t_0,\infty), \\
\vph(\ty,t_0) = 0 &\text{on } \S^N,
\end{cases}
\end{equation}
which automatically satisfies \eqref{eq:ortho120} and \eqref{eq:ortho12} for $s_0 = \infty$.
Also, $\vph = \vph[\wtmch]$ and $\tc_0 = \tc_0[\wtmch]$ are linear in $\wtmch$, and
\begin{equation}\label{eq:inhom-i22}
\|\vph\|_{\wtsh,a+2,b} + \sup_{t \in [t_0,\infty)} \mu_0^{-b}(t) |\tc_0(t)| \le C\left\|\wtmch\right\|_{\wtsh,a,b}
\end{equation}
for a constant $C > 0$ depending only on $N$, $a$ and $b$.
\end{cor}
\begin{proof}
Given an increasing sequence $\{s_{\ell}\}_{\ell \in \N}$ such that $\frac{3t_0}{2} < s_{\ell} \to \infty$ as $\ell \to \infty$, let us consider the problem
\begin{equation}\label{eq:lin11}
\begin{cases}
\displaystyle p\, (\vph_{\ell})_t = \frac{\kappa_N p}{N} \(\Delta_{\S^N} \vph_{\ell} + N\vph_{\ell}\) + \wtmch_{\ell} - \tc_{\ell}(t) &\text{on } \S^N \times (t_0,s_{\ell}), \\
\vph_{\ell}(\ty,t_0) = 0 &\text{on } \S^N.
\end{cases}
\end{equation}
Here,
\[\wtmch_{\ell}(\ty,t) := \min\left\{\wtmch(\ty,t),\ell\right\} - \sum_{n=1}^{N+1} \td_{n\ell}(t) \ty_n \in L^{\infty}(\S^N \times (t_0,s_{\ell})),\]
\[\td_{n\ell}(t) := \frac{N+1}{\left|\S^N\right|} \int_{\S^N} \min\left\{\wtmch(\ty,t),\ell\right\} \ty_n dS_{\ty} \quad \text{for } t \in (t_0,s_{\ell}) \text{ and } n = 1, \ldots, N+1,\]
and $\tc_{\ell}$ is the quantity defined by \eqref{eq:tcl} in which $t_{\ell}$ is replaced with $t_0$.
Clearly, \eqref{eq:ortho11l} holds provided $t_{\ell} = t_0$, and $\td_{n\ell}(t) \to 0$ and $\wtmch_{\ell} \to \wtmch$ a.e. as $\ell \to \infty$.

By applying the Galerkin method, we can construct a solution $\vph_{\ell}$ to \eqref{eq:lin11}
such that $\vph_{\ell}$, $D_{\ty}\vph_{\ell}$, $D^2_{\ty}\vph_{\ell}$, $(\vph_{\ell})_t \in L^2(\S^N \times (t_0,s_{\ell}))$.
Lemma \ref{lemma:lin} tells us that such $\vph_{\ell}$ is unique and
\begin{equation}\label{eq:lin12}
\|\vph_{\ell}\|_{\wtsh,a+2,b;t_0,s_{\ell}} + \sup_{t \in [t_0,s_{\ell})} \mu_0^{-b}(t) |\tc_{\ell}(t)|
\le C \left\|\wtmch_{\ell}\right\|_{\wtsh,a,b;t_0,s_{\ell}} \le C \left\|\wtmch\right\|_{\wtsh,a,b}
\end{equation}
for a constant $C > 0$ independent of $\ell \in \N$. By \eqref{eq:ortho11l}, we also have \eqref{eq:ortho12l} provided $t_{\ell} = t_0$.

Fix $\vep \in (0,\frac{a}{N-a})$. Given any $s > t_0$, parabolic regularity and \eqref{eq:lin12} yield that
\[\left\||\vph_{\ell}| + |D_{\ty}\vph_{\ell}| + |D^2_{\ty}\vph_{\ell}| + |(\vph_{\ell})_t|\right\|_{L^{1+\vep}(\S^N \times (t_0,s))} \le C\]
for some $C > 0$ independent of $\ell \in \N$. Therefore, we have a function $\vph_{\infty}$ on $\S^N \times (t_0,\infty)$ such that
\[\vph_{\ell} \to \vph_{\infty} \ \begin{cases}
\text{weakly in } W^{1,1+\vep}_{\text{loc}}(\S^N \times (t_0,\infty)),\\
\text{strongly in } L^{1+\vep}_{\text{loc}}(\S^N \times (t_0,\infty)),\\
\text{a.e. on } \S^N \times (t_0,\infty)
\end{cases} \quad \text{as } \ell \to \infty,\]
along a subsequence. In particular, $\vph_{\infty}$ satisfies \eqref{eq:inhom-i21} with $\tc_0(t) = \lim_{\ell \to \infty} \tc_{\ell}(t)$.
Setting $\vph = \vph_{\infty}$, we also infer from \eqref{eq:ortho12l} and \eqref{eq:lin12} that \eqref{eq:ortho120}--\eqref{eq:ortho12} for $s_0 = \infty$ and \eqref{eq:inhom-i22} are valid.

\medskip
The linear dependence of $\vph$ and $\tc_0$ in $\wtmch$ is obvious. This completes the proof.
\end{proof}

\begin{cor}\label{cor:lin2}
Assume that $\|\mch\|_{\sharp,a+2,b} < \infty$ and \eqref{eq:ortho1} holds.
Then one can find a solution $\phi$ to the equation
\begin{equation}\label{eq:inhom-i3}
\begin{cases}
\begin{aligned}
p W_{1,0}^{p-1} \phi_t &= \frac{(N+2)\kappa_N}{4} \(\Delta \phi + p W_{1,0}^{p-1} \phi\) \\
&\hspace{93pt} + W_{1,0}^{p-1} \mch - c_0(t) W_{1,0}^{p-1}Z_0
\end{aligned}
&\text{in } \R^N \times (t_0,\infty),\\
\phi(\by,t_0) = 0 &\text{in } \R^N
\end{cases}
\end{equation}
for a suitable $c_0$, which also satisfies
\begin{equation}\label{eq:ortho13}
\int_{\R^N} \phi(\by,t) \(W_{1,0}^{p-1}Z_n\)(\by) d\by = 0 \quad \text{for } t \in [t_0,\infty) \text{ and } n = 0, \ldots, N+1.
\end{equation}
Also, $\phi = \phi[\mch]$ and $c_0 = c_0[\mch]$ are linear in $\mch$, and
\begin{equation}\label{eq:lin13}
\|\phi\|_{\sharp',a,b(\R^N)} + \sup_{t \in [t_0,\infty)} \mu_0^{-b}(t) |c_0(t)| \le C\|\mch\|_{\sharp,a+2,b}
\end{equation}
for a constant $C > 0$ depending only on $N$, $a$ and $b$.
\end{cor}
\begin{proof}
Let $\Pi: \R^N \to \S^N_n$ and $\Pi_*f: \S^N_n \to \R$ be the inverse of the stereographic projection
and the weighted push-forward of $f: \R^N \to \R$ in Definition \ref{defn:stereo}, respectively.

We set $\wtmch = \Pi_*\mch$ for each fixed $t \in [t_0,\infty)$.
If $d_{\,\S^N}(\ty)$ is the geodesic distance on $\S^N$ between $\tyn = (0,\ldots,0,1)$ and $\ty = \Pi(\by)$ for $\by \in \R^N$, then
\[d_{\,\S^N}(\ty) = \arccos\(\frac{|\by|^2-1}{|\by|^2+1}\) = \frac{2}{|\by|} \left[1+O\(\frac{1}{|\by|^2}\)\right] \quad \text{for } |\by| \text{ large}.\]
Hence, $\|\mch\|_{\sharp,a+2,b} < \infty$ is reduced to $\left\|\wtmch\right\|_{\wtsh,a,b} < \infty$.
Besides, each $\Pi_*Z_{\ell}$ being a constant multiple of $\ty_{\ell}$, \eqref{eq:ortho1} is transformed into \eqref{eq:ortho11} with $s_0 = \infty$.

Given the solution $(\vph, \tc_0)$ to \eqref{eq:inhom-i21} deduced from Corollary \ref{cor:lin}, we set $(\phi, c_0)$ by
\[\vph = \Pi_*\phi \quad \text{and} \quad \tc_0(t) := c_0(t) (\Pi_*Z_0) \quad \text{for } t \in [t_0,\infty),\]
noting that $\Pi_*Z_0$ is a positive constant. By virtue of conformality of the map $\Pi$, it satisfies \eqref{eq:inhom-i3}--\eqref{eq:lin13}.
\end{proof}

\begin{cor} 
Suppose that all the assumptions on Proposition \ref{prop:lin} hold.
Let $\phi$ be a solution to \eqref{eq:inhom-i3} satisfying \eqref{eq:lin13}.
Then there exists a constant $C > 0$ depending only on $N$, $a$, $b$ and $\sigma_0$ such that
\begin{equation}\label{eq:lin21}
\|\phi\|_{\sharp',a,b;\sigma_0(\R^N)} + \|c_0\|_{b;\sigma_0} \le C\|\mch\|_{\sharp,a+2,b;\sigma_0}.
\end{equation}
\end{cor}
\begin{proof}
From the relation
\[c_0(t) = \int_{\R^N} W_{1,0}^{p-1}Z_0 \mch(\by,t) d\by \quad \text{for } t \in [t_0,\infty),\]
it is easy to see that $[c_0]_{C^{\sigma_0/2}_t}$ is controlled by $\|\mch\|_{\sharp,a+2,b;\sigma_0}$ provided $a > \sigma_0$.
Using this fact and \eqref{eq:lin13}, one can estimate $\|\phi\|_{\sharp',a,b;\sigma_0(\R^N)}$ as in Step 1 in the proof of Corollary \ref{cor:outer2}. We skip the details.
\end{proof}

\begin{proof}[Proof of Proposition \ref{prop:lin}]
Let $\phi$ be the solution to \eqref{eq:inhom-i3} satisfying \eqref{eq:ortho13} and \eqref{eq:lin21} found in Corollary \ref{cor:lin2}.
If we set $\psi = \phi + e_0(t)Z_0$, it satisfies \eqref{eq:inhom-i} provided $e_0$ satisfies
\[\dot{e}_0 - \kappa_N e_0 = p^{-1}c_0 \quad \text{on } [t_0,\infty),\]
which is explicitly given by
\[e_0(t) = - \frac{1}{p} \int_t^{\infty} \exp\(\kappa_N (t-\tau)\) c_0(\tau) d\tau \quad \text{for } t \in [t_0,\infty).\]
The function $e_0$ is linear in $\mch$. Moreover, by \eqref{eq:lin21},
\[\|e_0\|_{b;\sigma_0} \le C \sup_{t \in [t_0,\infty)} \mu_0^{-b}(t) |c_0(t)| \le C\|\mch\|_{\sharp,a+2,b},\]
which together with the condition $a < N-2$ yields
\[\|\psi\|_{\sharp',a,b;\sigma_0(\R^N)} \le \|\phi\|_{\sharp',a,b;\sigma_0(\R^N)} + \|e_0 Z_0\|_{\sharp',a,b;\sigma_0(\R^N)} \le C\|\mch\|_{\sharp,a+2,b;\sigma_0}\]
as desired.
\end{proof}

\subsection{Choice of parameters} \label{subsec:inner2}
As mentioned at the beginning of this section, we reduce \eqref{eq:ortho2} into a system of nonlinear ODEs of $(\lambda, \xi_1, \ldots, \xi_N)$.
Recall $\nu_1$ and $\nu_2$ in \eqref{eq:lx}, $\delta_4$ in \eqref{eq:delta4},
and $a$, $b$, $\alpha$, $\beta$, $\rho$, $\sigma_0$, and $\psi^{\out} = \psi^{\out}[\lambda, \xi, \dot{\lambda}, \dot{\xi}, \hps^{\, \tin}, \psi_0]$ in Proposition \ref{prop:outer3}.

\begin{lemma}\label{lemma:xi}
Assume that \eqref{eq:outer31} and \eqref{eq:mu0}--\eqref{eq:lx} hold. For $n = 1, \ldots, N$, \eqref{eq:ortho2} is equivalent to
\begin{equation}\label{eq:ODExi}
\begin{medsize}
\displaystyle \begin{pmatrix}
\dot{\xi}_1 \\ \vdots \\ \dot{\xi}_n
\end{pmatrix}
- \frac{N-4}{6h(z_0)t} \begin{pmatrix}
R_{11}(z_0) & \cdots & R_{1N}(z_0) \\
\vdots & \ddots & \vdots \\
R_{N1}(z_0) & \cdots & R_{NN}(z_0)
\end{pmatrix}
\begin{pmatrix}
\xi_1 \\ \vdots \\ \xi_n
\end{pmatrix}
= \mu_0^4 \begin{pmatrix}
\Theta_{11} \\ \vdots \\ \Theta_{N1}
\end{pmatrix}
+ \mu_0^{b+\nu_1} \begin{pmatrix}
\Theta_{12} \\ \vdots \\ \Theta_{N2}
\end{pmatrix}
+ \mu_0^{\delta_4}(t_0) \mu_0^{\beta+1} \begin{pmatrix}
\Theta_{13} \\ \vdots \\ \Theta_{N3}
\end{pmatrix}
\end{medsize}
\end{equation}
on $[t_0,\infty)$. Here, $\Theta_n = \Theta_{n1}, \Theta_{n2}, \Theta_{n3}$ is a function on $[t_0,\infty)$ such that
\begin{itemize}
\item[-] $\|\Theta_n\|_{0;\sigma_0} \le C_5\delta_0^{-\zeta_5}$ for some $C_5,\, \zeta_5 > 0$ where the norm is defined in \eqref{eq:nunorm};
\item[-] It depends on the parameters $\lambda,\, \dot{\lambda},\, \xi,\, \dot{\xi}$ and the functions $\psi^{\tin},\, \psi_0$. Besides, it satisfies
\begin{equation}\label{eq:Thetad}
\begin{cases}
\left\|\Theta_n[\lambda_1] - \Theta_n[\lambda_2]\right\|_{0;\sigma_0}
\le C\delta_0^{-\zeta} \|\lambda_1-\lambda_2\|_{\nu_1;\sigma_0},\\
\left\|\Theta_n[\xi_1] - \Theta_n[\xi_2]\right\|_{0;\sigma_0}
\le C\delta_0^{-\zeta} \|\xi_1-\xi_2\|_{\nu_2;\sigma_0},\\
\left\|\Theta_n[\dot{\lambda_1}] - \Theta_n[\dot{\lambda_2}]\right\|_{0;\sigma_0}
\le C\delta_0^{-\zeta} \|\dot{\lambda}_1-\dot{\lambda}_2\|_{\nu_1+2;\sigma_0},\\
\left\|\Theta_n[\dot{\xi}_1] - \Theta_n[\dot{\xi}_2]\right\|_{0;\sigma_0}
\le C\delta_0^{-\zeta} \|\dot{\xi}_1-\dot{\xi}_2\|_{\nu_2+2;\sigma_0},\\
\left\|\Theta_n[\psi^{\tin}_1] - \Theta_n[\psi^{\tin}_2]\right\|_{0;\sigma_0}
\le C\delta_0^{-\zeta} \left\|\psi^{\tin}_1 - \psi^{\tin}_2\right\|_{\sharp',a,b;\sigma_0\(\mcb_{\bmx}\)}
\end{cases}
\end{equation}
and
\begin{equation}\label{eq:Thetad2}
\left\|\Theta_n[\psi_{01}] - \Theta_n[\psi_{02}]\right\|_{0;\sigma_0}
\le C\delta_0^{-\zeta} \mu_0^{-\delta_4}(t_0) \left\|\psi_{01} - \psi_{02}\right\|_{**,\alpha;\sigma_0}
\end{equation}
for $C,\, \zeta > 0$ large depending only on $(M,g_0)$, $N$, $h$, $z_0$, $a$, $b$, $\alpha$, and $\sigma_0$.
\end{itemize}

\end{lemma}
\begin{proof}
Throughout the proof, we will use the notation $\Theta_n$ to refer functions which behave as in the statement of the lemma.
They may vary from line to line and even in the same line.

\medskip
Fixing $n = 1,\ldots,N$, we will examine
\begin{equation}\label{eq:Bn}
\begin{aligned}
B_{n1} &:= \({\bmu \over \mu}\)^{N-2 \over 2} \int_{\R^N} \eta_{2\mu_0^{\vep_1}}(x) (1+P)^{-p}(x) \mce_2[\mx](y,t) Z_n(\by) d\by,\\
B_{n2} &:= \int_{\R^N} \eta_{2\mu_0^{\vep_1}}(x) \mck_1\left[\psi^{\tin}\right] Z_n(\by) d\by,
\quad B_{n3} := \int_{\R^N} \eta_{2\mu_0^{\vep_1}}(x) \mck_2\left[\psi^{\tin}, \psi^{\out}\right] Z_n(\by) d\by
\end{aligned}
\end{equation}
for $t \in [t_0,\infty)$, respectively. Here, $x = \mu y+\xi = \bmu\by+\xi$.

\medskip \noindent \textsc{Estimate on $B_{n1}$}: To estimate $B_{n1}$, we treat the quantity $\mce_2[\mx](y,t)$ in \eqref{eq:E2} term by term.

First, we have
\begin{align*}
&\ \bmu^{-1} \(\dot{\lambda} - \bmu^{-1} \dot{\bmu} \lambda\) \({\bmu \over \mu}\)^{N-2 \over 2} \int_{\R^N} \eta_{2\mu_0^{\vep_1}}(x) (1+P)^{-p}(x) \(pW_{1,0}^{p-1}Z_{N+1}\)(y) Z_n(\by) d\by \\
&= \mu_0^{\nu_1+1} \Theta_n \(1+ \mu_0^{\nu_1-1}\Theta_n\)
\left[\int_{\R^N} \(\eta_{2\mu_0^{\vep_1}}(x)-1\) \(pW_{1,0}^{p-1}Z_{N+1}Z_n\)(\by) d\by \right.\\
&\hspace{130pt} + \int_{\R^N} \eta_{2\mu_0^{\vep_1}}(x) \left\{\(pW_{1,0}^{p-1}Z_{N+1}\)(y) - \(pW_{1,0}^{p-1}Z_{N+1}\)(\by)\right\} Z_n(\by) d\by \\
&\hspace{130pt} \left. + \int_{\R^N} \eta_{2\mu_0^{\vep_1}}(x) \left\{(1+P)^{-p}(x)-1\right\} \(pW_{1,0}^{p-1}Z_{N+1}\)(y) Z_n(\by) d\by\right] \\
&= \mu_0^{\nu_1+1} \Theta_n \(\mu_0^{(N+1)(1-\vep_1)} \Theta_n + \mu_0^{\nu_1-1} \Theta_n + \mu_0^2 \Theta_n\) = \mu_0^{2\nu_1} \Theta_n,
\end{align*}
and similarly,
\[\({\bmu \over \mu}\)^{N-2 \over 2} \mu^{-1} \dot{\xi} \cdot \int_{\R^N} \eta_{2\mu_0^{\vep_1}}(x) (1+P)^{-p}(x) \(p W_{1,0}^{p-1} \nabla W_{1,0}\)(y) Z_n(\by) d\by \\
= c_5 \bmu^{-1} \dot{\xi}_n + \mu_0^{\nu_1+\nu_2} \Theta_n\]
where
\begin{equation}\label{eq:c5}
c_5 := p \int_{\R^N} W_{1,0}^{p-1}Z_1^2 > 0.
\end{equation}
Second, we infer from \eqref{eq:F0} and the parity that
\[2\bmu\lambda \({\bmu \over \mu}\)^{N-2 \over 2} \int_{\R^N} \eta_{2\mu_0^{\vep_1}}(x) (1+P)^{-p}(x) \mcf_0(y) Z_n(\by) d\by
= 2\bmu\lambda \int_{\R^N} (\mcf_0Z_n)(\by) d\by + \mu_0^{2\nu_1} \Theta_n = \mu_0^{2\nu_1} \Theta_n.\]
Third, \eqref{eq:F1} and the identity
\[\frac{2}{N-2} \int_{\R^N} (y \cdot \nabla W_{1,0}(y)) W_{1,0}^p(y)dy = - \int_{\R^N} |\nabla W_{1,0}|^2\]
imply
\begin{multline*}
\({\bmu \over \mu}\)^{N-2 \over 2} \int_{\R^N} \eta_{2\mu_0^{\vep_1}}(x) (1+P)^{-p}(x) \mcf_1[\mx](y) Z_n(\by) d\by \\
= \int_{\R^N} (\mcf_1[\bmu,\xi] Z_n)(\by) d\by + \mu_0^{\nu_1+\nu_2} \Theta_n
= - \frac{(N+2)\kappa_N}{12N} \bmu \xi_j R_{jn}(z_0) \int_{\R^N} |\nabla W_{1,0}|^2 + \mu_0^{\nu_1+\nu_2} \Theta_n.
\end{multline*}
Lastly,
\begin{multline*}
\({\bmu \over \mu}\)^{N-2 \over 2} \int_{\R^N} \eta_{2\mu_0^{\vep_1}}(x) (1+P)^{-p}(x) \left[\mu^3 a^{\{1\}} + \mu^{2\nu_1} a^{\{2\}} \right. \\
\left. + \mu\dot{\mu} a^{\{0\}} + \mu^{\nu_1-2} \dot{\lambda}\, a^{\{-2\}} + \mu\dot{\xi} \cdot \mba^{\{-1\}}\right] Z_n(\by) d\by
= \mu_0^3 \Theta_n.
\end{multline*}

From the above calculations, \eqref{eq:mu0}, and \eqref{eq:c5}, we conclude that
\begin{equation}\label{eq:ODExi1}
B_{n1} = c_5 \bmu^{-1} \left[\dot{\xi}_n - \frac{N-4}{6h(z_0)t} R_{jn}(z_0) \xi_j \right] + \mu_0^3 \Theta_n.
\end{equation}

\medskip \noindent \textsc{Estimate on $B_{n2}$}: To estimate $B_{n2}$, we handle the quantity $\mck_1[\psi^{\tin}]$ in \eqref{eq:mck1} term by term.

As an illustration, we consider the integrals involving the first and third terms of $\mck_1[\psi^{\tin}]$:
The mean value theorem, \eqref{eq:sharp'-norm}, \eqref{eq:outer31}, and the assumption $a > \sigma_0$ lead to
\[\int_{\R^N} \eta_{2\mu_0^{\vep_1}}(x) \left[W_{1,0}^{p-1}(\by) - W_{1,0}^{p-1}(y)\right] \(\psi^{\tin}\)_t Z_n(\by) d\by = \mu_0^{b+\nu_1-1} \Theta_n\]
and
\begin{equation}\label{eq:ODExi4}
\begin{aligned}
&\ \int_{\R^N} \eta_{2\mu_0^{\vep_1}}(x) \left[\({\mu \over \bmu}\)^2 (1+P)^{1-p}(x) (\Delta_{g_0(x)} \psi^{\tin})(\by,t)
- \Delta \psi^{\tin}(\by,t)\right] Z_n(\by) d\by \\
&= \({\mu \over \bmu}\)^2 \int_{\R^N} \eta_{2\mu_0^{\vep_1}}(x) \left[(1+P)^{1-p}(x)-1\right] \(\Delta_{g_0(x)} \psi^{\tin}\)(\by,t) Z_n(\by) d\by \\
&\ + \left[\({\mu \over \bmu}\)^2-1\right] \int_{\R^N} \eta_{2\mu_0^{\vep_1}}(x) \(\Delta_{g_0(x)} \psi^{\tin}\)(\by,t) Z_n(\by) d\by \\
&\ + \int_{\R^N} \eta_{2\mu_0^{\vep_1}}(x) \left[\(\Delta_{g_0(\bmu\cdot+\xi)}-\Delta\) \psi^{\tin}\right](\by,t) Z_n(\by) d\by \\
&= \mu_0^{b+1} \Theta_n
+ \mu_0^{b+\nu_1-1} \Theta_n + \(\mu_0^{a+b+1-\sigma_0} \Theta_n + \mu_0^{b+2} \Theta_n\) = \mu_0^{b+\nu_1-1} \Theta_n.
\end{aligned}
\end{equation}

In a similar manner, one can compute the integrals involving the remaining terms of $\mck_1[\psi^{\tin}]$, deriving
\begin{equation}\label{eq:ODExi2}
B_{n2} = \mu_0^{b+\nu_1-1} \Theta_n.
\end{equation}

\medskip \noindent \textsc{Estimate on $B_{n3}$}: We consider two terms of $\mck_2[\psi^{\tin}, \psi^{\out}]$ in \eqref{eq:mck2} separately.

First, by \eqref{eq:*'s-norm} and \eqref{eq:outer32},
\[\({\bmu \over \mu}\)^{N-2 \over 2} \int_{\R^N} \eta_{2\mu_0^{\vep_1}}(x) (1+P)^{-1}(x) W_{1,0}^{p-1}(y) \mu^{N-2 \over 2} \psi^{\out}(x,t) Z_n(\by) d\by
= \mu_0^{\delta_4}(t_0) \mu_0^{\beta} \Theta_n.\]

Second, the mean value theorem yields
\begin{align*}
&\ \({\bmu \over \mu}\)^{N-2 \over 2} \int_{\R^N} \eta_{2\mu_0^{\vep_1}}(x) \left[\left\{\(W_{1,0}+\bmu^2 Q_0\)(y)
+ \({\mu \over \bmu}\)^{N-2 \over 2} \psi^{\tin} + (1+P)^{-1}(x) \mu^{N-2 \over 2} \psi^{\out}(x,t)\right\}^{p-1} \right. \\
&\hspace{105pt} \left. - \(W_{1,0}+\bmu^2 Q_0\)^{p-1}(y) \right] \eta_{\mu_0^{\vep_1}/2}\,
\mu^{N-2 \over 2} \left\{\mu^{-{N-2 \over 2}} \(W_{1,0}+\bmu^2 Q_0\)(y) \right\}_t Z_n(\by) d\by
\\
&= \mu_0^{b+2} \Theta_n + \mu_0^{\delta_4}(t_0) \mu_0^{\beta+2} \Theta_n.
\end{align*}

Consequently,
\begin{equation}\label{eq:ODExi3}
B_{n3} = \mu_0^{b+2} \Theta_n + \mu_0^{\delta_4}(t_0) \mu_0^{\beta} \Theta_n.
\end{equation}

\medskip
Notice that \eqref{eq:ortho2} is equivalent to the equation $B_{n1} + B_{n2} + B_{n3} = 0$.
By combining \eqref{eq:ODExi1}, \eqref{eq:ODExi2}, and \eqref{eq:ODExi3}, we establish \eqref{eq:ODExi}.

\medskip
A closer look at the above computations with \eqref{eq:outer35}--\eqref{eq:outer36} gives \eqref{eq:Thetad}--\eqref{eq:Thetad2}. The details are omitted.
\end{proof}

\begin{lemma}
Assume that \eqref{eq:outer31} and \eqref{eq:mu0}--\eqref{eq:lx} hold.
For $n = N+1$, \eqref{eq:ortho2} is equivalent to
\begin{equation}\label{eq:ODEla}
\dot{\lambda}(t) + \frac{3}{2t} \lambda(t) = \mu_0^4 \Theta_{(N+1)1}
+ \mu_0^{b+\nu_1} \Theta_{(N+1)2} + \mu_0^{a+b+1-\sigma_0} \Theta_{(N+1)3} + \mu_0^{\delta_4}(t_0) \mu_0^{\beta+1} \Theta_{(N+1)4}
\end{equation}
on $[t_0,\infty)$. Here, $\Theta_{N+1} = \Theta_{(N+1)1}, \ldots, \Theta_{(N+1)4}$ is a function on $[t_0,\infty)$
that behaves as the function $\Theta_n$ described in the statement of Lemma \ref{lemma:xi}.
\end{lemma}
\begin{proof}
The proof is analogous to that of the previous lemma.
Let $B_{(N+1)1}$, $B_{(N+1)2}$, and $B_{(N+1)3}$ be the quantities obtained by taking $n = N+1$ in \eqref{eq:Bn}.
It suffices to estimate each of them.

\medskip \noindent \textsc{Estimate on $B_{(N+1)1}$}:
In view of \eqref{eq:c1c2} and \eqref{eq:mu0}--\eqref{eq:mu} (see also \eqref{eq:bmu}), we know
\[\frac{N+2}{2} \bmu^2 h(z_0) c_2 - \bmu^{-1} \dot{\bmu} c_1 = \frac{3}{2t} c_1.\]
By \eqref{eq:E2}, \eqref{eq:orthoN+1}, and the above identity, the dominating term of $B_{(N+1)1}$ turns out to be
\begin{multline*}
\int_{\R^N} \left[\bmu^{-1} \(\dot{\lambda} - \bmu^{-1} \dot{\bmu}\lambda\) pW_{1,0}^{p-1}(\by) + 2\bmu\lambda \mcf_0(\by)\right] Z_{N+1}(\by) d\by \\
= \bmu^{-1} \(\dot{\lambda} - \bmu^{-1} \dot{\bmu} \lambda\) c_1 + \frac{N+2}{2} \bmu\lambda h(z_0) c_2
= \bmu^{-1} c_1 \(\dot{\lambda} + \frac{3}{2t} \lambda\).
\end{multline*}
Treating the other parts of $B_{(N+1)1}$, we obtain
\[B_{(N+1)1} = \bmu^{-1} c_1 \(\dot{\lambda} + \frac{3}{2t} \lambda\) + \mu_0^3 \Theta_{N+1}.\]

\medskip \noindent \textsc{Estimate on $B_{(N+1)2}$ and $B_{(N+1)3}$}:
Arguing as in the derivation of \eqref{eq:ODExi2} and \eqref{eq:ODExi3}, we deduce
\[B_{(N+1)2} = \mu_0^{b+\nu_1-1} \Theta_{N+1} + \mu_0^{a+b-\sigma_0} \Theta_{N+1}\]
and
\[B_{(N+1)3} = \mu_0^{b+2} \Theta_{N+1} + \mu_0^{\delta_4}(t_0) \mu_0^{\beta} \Theta_{N+1}.\]
Observe that the above estimate on $B_{(N+1)2}$ has one more term than that on $B_{n2}$ in \eqref{eq:ODExi2}.
It is because the integrand of $B_{(N+1)2}$ decays slower than that of $B_{n2}$ as $|\by| \to \infty$; compare \eqref{eq:ZN+1} and \eqref{eq:Zi}.
\end{proof}

We next solve the system \eqref{eq:ODExi} and \eqref{eq:ODEla} of nonlinear ODEs, thereby determining the parameters $(\lambda,\xi)$.
\begin{prop}\label{prop:inner1}
Suppose that \eqref{eq:outer31} holds and $t_0 > 0$ is large enough.
Given $\nu_1 = \nu_2 = 2-\vep_0$ in \eqref{eq:lx}, $\delta_4$ in \eqref{eq:delta4},
and $a$, $b$, $\alpha$, $\beta$, $\sigma_0$, and $\psi^{\out}$ in Proposition \ref{prop:outer3},
we further assume that $\delta_4 \le \vep_0$, $a \in [\sigma_0+\vep_0,N-2)$, and $b = \beta = 3-\vep_0$.
Then there exists a solution $(\lambda[\hps^{\, \tin}, \psi_0], \xi[\hps^{\, \tin}, \psi_0])$
to the system \eqref{eq:ODExi} and \eqref{eq:ODEla} of ODEs which satisfies \begin{equation}\label{eq:lx2}
\|\lambda\|_{\nu_1;\sigma_0} + \|\dot{\lambda}\|_{\nu_1+2;\sigma_0} + \|\xi\|_{\nu_2;\sigma_0} + \|\dot{\xi}\|_{\nu_2+2;\sigma_0} \le C \delta^{-\zeta} \mu_0^{\delta_4}(t_0)
\end{equation}
where $C,\, \zeta > 0$ are constants depending only on $(M,g_0)$, $N$, $h$, $z_0$, $a$, $b$, $\alpha$, $\sigma_0$, and $\vep_0$; compare with \eqref{eq:lx}. Additionally,
\begin{equation}\label{eq:inner1a}
\left\|\lambda\big[\hps^{\, \tin}_1\big] - \lambda\big[\hps^{\, \tin}_2\big]\right\|_{\nu_1}
+ \left\|\xi\big[\hps^{\, \tin}_1\big] - \xi\big[\hps^{\, \tin}_2\big]\right\|_{\nu_2}
\le C \delta_0^{-\zeta} \mu_0^{\delta_4}(t_0) \left\|\hps^{\, \tin}_1 - \hps^{\, \tin}_2\right\|_{\sharp',a,b;\sigma_0\(\mcb_{\bmx}\)}
\end{equation}
and
\begin{equation}\label{eq:inner1b}
\left\|\lambda\left[\psi_{01}\right] - \lambda\left[\psi_{02}\right]\right\|_{\nu_1}
+ \left\|\xi\left[\psi_{01}\right] - \xi\left[\psi_{02}\right]\right\|_{\nu_2}
\le C \delta_0^{-\zeta} \left\|\psi_{01}-\psi_{02}\right\|_{**,\alpha;\sigma_0}.
\end{equation}
\end{prop}
\begin{proof}
The proof is decomposed into three steps.

\medskip \noindent \textsc{Step 1.} Let $\mcm$ be an $N \times N$ symmetric matrix and
$\ovmfh(t) := (\mfh_1(t), \ldots, \mfh_N(t))$ a function on $[t_0,\infty)$ such that $\|\ovmfh\|_{\nu;\sigma} < \infty$ for some numbers $\nu > 2$ and $\sigma \in (0,1)$.
We assert that there exists a solution $\xi(t) = (\xi_1(t), \ldots, \xi_N(t))$ to the ODE system
\begin{equation}\label{eq:xieq}
\dot{\xi}(t) + \frac{1}{t} \mcm \xi(t) = \ovmfh(t) \quad \text{on } [t_0,\infty),
\quad \xi(t) \to 0 \quad \text{as } t \to \infty \quad \text{and} \quad \|\dot{\xi}\|_{\nu-\ep;\sigma} \le C \|\ovmfh\|_{\nu;\sigma}
\end{equation}
where $\xi$ and $\ovmfh$ are regarded as column vectors, $\ep > 0$ is any small number,
and $C > 0$ is a constant depending only on $N$, $\mcm$, $\nu$, $\ep$, and $\sigma$.

Indeed, since $\mcm$ is symmetric, there exist an orthogonal matrix $\mcq$ and a diagonal matrix $\mcd = \text{diag}(\vsi_1, \ldots, \vsi_N)$
such that $\mcm = \mcq^T \mcd \mcq$. If we set $\xi(t)$ by
\begin{equation}\label{eq:xieq2}
\xi(t) = \mcq^T \tilde{\xi}(t) \quad \text{where} \quad \tilde{\xi}_i(t) := \begin{cases}
\displaystyle t^{-\vsi_i} \int_{t_0}^t s^{\vsi_i} \(\mcq \ovmfh\)_i(s) ds &\text{if } \vsi_i \ge \frac{\nu}{2}-1, \\
\displaystyle - t^{-\vsi_i} \int_t^{\infty} s^{\vsi_i} \(\mcq \ovmfh\)_i(s)ds &\text{if } \vsi_i < \frac{\nu}{2}-1,
\end{cases}
\end{equation}
then it satisfies \eqref{eq:xieq}. Note that one can select $\ep = 0$ unless $\vsi_i = \frac{\nu}{2}-1$ for some $i = 1, \ldots, N$.

\medskip \noindent \textsc{Step 2.} Let $\mfh(t) = (\mfh_1(t), \ldots, \mfh_{N+1}(t))$ be a function on $[t_0,\infty)$ and
\begin{equation}\label{eq:mcm}
\mcm = -\frac{N-4}{6h(z_0)} \begin{pmatrix}
R_{11}(z_0) & \cdots & R_{1N}(z_0) \\
\vdots & \ddots & \vdots \\
R_{N1}(z_0) & \cdots & R_{NN}(z_0)
\end{pmatrix}.
\end{equation}
By adjusting $\vep_0$ suitably, we may assume that $\vsi_i \ne \frac{\nu_1}{2} = \frac{\nu_2}{2}$ for all $i = 1, \ldots, N$, and so we may take $\ep = 0$.
Let $\xi$ be the solution to \eqref{eq:xieq} chosen in the previous step,
\begin{equation}\label{eq:laeq}
\lambda(t) = t^{-{3 \over 2}} \int_{t_0}^t s^{3 \over 2} \mfh_{N+1}(s)ds
\quad \text{so that } \dot{\lambda}(t) + \frac{3}{2t} \lambda(t) = \mfh_{N+1}(t) \quad \text{on } [t_0,\infty),
\end{equation}
and
\begin{equation}\label{eq:mct}
\mct^{\tpa}[\mfh] = (\Xi, \Lambda) := (\dot{\xi}, \dot{\lambda}).
\end{equation}
By \eqref{eq:xieq} and \eqref{eq:laeq},
\begin{equation}\label{eq:mct2}
\begin{aligned}
\|\mct^{\tpa}[\mfh]\|_{\nu_2+2,\, \nu_1+2;\sigma_0} &= \|(\Xi,\Lambda)\|_{\nu_2+2,\, \nu_1+2;\sigma_0} := \|\Xi\|_{\nu_2+2;\sigma_0} + \|\Lambda\|_{\nu_1+2;\sigma_0}\\
&\le C_6 (\|(\mfh_1, \ldots, \mfh_n)\|_{\nu_2+2;\sigma_0} + \|\mfh_{N+1}\|_{\nu_1+2;\sigma_0}) = C_6\|\mfh\|_{\nu_2+2,\,\nu_1+2;\sigma_0}
\end{aligned}
\end{equation}
for some $C_6 > 0$.

\medskip \noindent \textsc{Step 3.} Parameters $(\lambda,\xi)$ solve the system \eqref{eq:ODExi} and \eqref{eq:ODEla} of nonlinear ODEs if and only if
\begin{equation}\label{eq:inner11}
(\Xi, \Lambda) = \mct^{\tpa}\left[f^{\tpa}\left[\Xi, \Lambda, \hps^{\, \tin}, \psi_0\right]\right]
\end{equation}
where
\[f^{\tpa}_n\left[\Xi, \Lambda, \hps^{\, \tin}, \psi_0\right]
:= \begin{cases}
\displaystyle \mu_0^4 \Theta_{n1} + \mu_0^{b+\nu_1} \Theta_{n2} + \mu_0^{\delta_4}(t_0) \mu_0^{\beta+1} \Theta_{n3} &\text{for } n = 1,\ldots,N,\\
\displaystyle \mu_0^4 \Theta_{n1} + \mu_0^{b+\nu_1} \Theta_{n2} + \mu_0^{a+b+1-\sigma_0} \Theta_{n3} + \mu_0^{\delta_4}(t_0) \mu_0^{\beta+1} \Theta_{n4} &\text{for } n = N+1
\end{cases}\]
and
\[f^{\tpa}\left[\Xi, \Lambda, \hps^{\, \tin}, \psi_0\right]
:= \(f^{\tpa}_1\left[\Xi, \Lambda, \hps^{\, \tin}, \psi_0\right], \ldots, f^{\tpa}_{N+1}\left[\Xi, \Lambda, \hps^{\, \tin}, \psi_0\right]\)\]
on $[t_0,\infty)$.

\medskip
We claim that there exists a point $(\Xi, \Lambda)$ on the set
\[\mcd^{\tpa} := \left\{\mfh: \|\mfh\|_{\nu_2+2,\, \nu_1+2;\sigma_0} \le 4\delta_0^{-\zeta_5} \mu_0^{\delta_4}(t_0) C_5C_6\right\}\]
for which \eqref{eq:inner11} holds, where $C_5,\, \zeta_5 > 0$ are the numbers in the statement of Lemma \ref{lemma:xi}.

To see this, we infer from \eqref{eq:mct2} and the conditions $a \ge \sigma_0+\vep_0$ and $b+1 = \beta+1 = \nu_2+2$ that
\begin{align*}
\left\|\mct^{\tpa}\left[f^{\tpa}\left[\Xi, \Lambda, \hps^{\, \tin}, \psi_0\right]\right]\right\|_{\nu_2+2,\, \nu_1+2;\sigma_0}
&\le C_6 \left\|f^{\tpa}\left[\Xi, \Lambda, \hps^{\, \tin}, \psi_0\right]\right\|_{\nu_2+2,\, \nu_1+2;\sigma_0} \\
&\le 4\delta_0^{-\zeta_5} \mu_0^{\delta_4}(t_0) C_5C_6.
\end{align*}
In addition, \eqref{eq:Thetad} says
\begin{align*}
&\ \left\|\mct^{\tpa}\left[f^{\tpa}\left[\Xi_1, \Lambda_1, \hps^{\, \tin}, \psi_0\right]\right]
- \mct^{\tpa}\left[f^{\tpa}\left[\Xi_2, \Lambda_2, \hps^{\, \tin}, \psi_0\right]\right] \right\|_{\nu_2+2,\, \nu_1+2;\sigma_0} \\
&\le C \left\|f^{\tpa}\left[\Xi_1, \Lambda_1, \hps^{\, \tin}, \psi_0\right]
- f^{\tpa}\left[\Xi_2, \Lambda_2, \hps^{\, \tin}, \psi_0\right]\right\|_{\nu_2+2,\, \nu_1+2;\sigma_0} \\
&\le C \delta_0^{-\zeta} \mu_0^{\delta_4}(t_0) \|(\Xi_1,\Lambda_1) - (\Xi_2,\Lambda_2)\|_{\nu_2+2,\, \nu_1+2;\sigma_0}
\le \frac{1}{2} \|(\Xi_1,\Lambda_1) - (\Xi_2,\Lambda_2)\|_{\nu_2+2,\, \nu_1+2;\sigma_0}
\end{align*}
provided $t_0 > 0$ large enough. By the contraction mapping theorem, the assertion holds.

Checking \eqref{eq:lx2} for the associated parameters $(\lambda,\xi)$ is a simple task.
In addition, a closer look at the above computations with \eqref{eq:Thetad}--\eqref{eq:Thetad2} gives \eqref{eq:inner1a}--\eqref{eq:inner1b}. The proof is finished.
\end{proof}

\subsection{Unique solvability of \eqref{eq:inner}} \label{subsec:inner3}
We are in position to solve the inner problem \eqref{eq:inner}. We define
\begin{multline}\label{eq:fin}
f^{\tin}_{\mx}\left[\psi^{\tin}, \psi^{\out}\right] \\
= \eta_{2\mu_0^{\vep_1}}(x) W_{1,0}^{1-p} \left[\({\bmu \over \mu}\)^{N-2 \over 2} (1+P)^{-p}(x) \mce_2[\mx](y,t)
+ \mck_1\left[\psi^{\tin}\right] + \mck_2\left[\psi^{\tin}, \psi^{\out}\right]\right];
\end{multline}
cf. \eqref{eq:ortho2}. If a function $\psi^{\tin}$ satisfies
\begin{equation}\label{eq:inner2}
\begin{cases}
\displaystyle p W_{1,0}^{p-1} \(\psi^{\tin}\)_t = \frac{(N+2)\kappa_N}{4} \(\Delta \psi^{\tin} + p W_{1,0}^{p-1} \psi^{\tin}\) + W_{1,0}^{p-1} f^{\tin}_{\mx}\left[\psi^{\tin}, \psi^{\out}\right]
&\text{in } \R^N \times (t_0,\infty), \\
\displaystyle \psi^{\tin}(\cdot,t_0) = e_0(t_0) Z_0 &\text{in } \R^N
\end{cases}
\end{equation}
for some $e_0(t_0) \in \R$, then it solves \eqref{eq:inner} in $B^N(0,\mcb_{\bmx}) \times (t_0,\infty)$
where $\mcb_{\bmx} = B^N(-\bmu^{-1}\xi, 2\bmu^{-1}\mu_0^{\vep_1})$.
Furthermore, estimate \eqref{eq:inner3} given below clearly implies the first inequality of \eqref{eq:outer31}.

\begin{prop}\label{prop:inner2}
Suppose that the second inequality in \eqref{eq:outer31} holds and $t_0 > 0$ is large enough.
Given numbers $a$, $b$, $\alpha$, $\beta$, and $\sigma_0$ in Proposition \ref{prop:inner1},
we further assume that $a \in [\sigma_0+\vep_0,\min\{N-4,\alpha+2\}]$.
Then one can find $\psi^{\tin} = \psi^{\tin}[\psi_0]$ and $e_0 = e_0[\psi_0]$ satisfying \eqref{eq:inner2},
\begin{equation}\label{eq:inner3}
\left\|\psi^{\tin}\right\|_{\sharp',a,b;\sigma_0\(\R^N\)} + \mu_0^{-b}(t_0) |e_0(t_0)| \le C,
\end{equation}
and
\begin{multline}\label{eq:inner31}
\left\|\psi^{\tin}[\psi_{01}] - \psi^{\tin}[\psi_{02}] \right\|_{\sharp',a,b;\sigma_0\(\mcb_{\bmx}\)} + \mu_0^{-b}(t_0)|e_0[\psi_{01}](t_0) - e_0[\psi_{02}](t_0)| \\
\le C\delta_0^{-\zeta} \|\psi_{01}-\psi_{02}\|_{**,\alpha;\sigma_0}
\end{multline}
for a constant $C,\, \zeta > 0$ depending only on $(M,g_0)$, $N$, $h$, $z_0$, $a$, $b$, $\alpha$, and $\sigma_0$.
Here, $(\lambda,\xi)$ and $\psi^{\out}$ are the ones determined in Propositions \ref{prop:inner1} and \ref{prop:outer3}, respectively.
\end{prop}
\noindent As in the proof of Proposition \ref{prop:outer3}, we start by estimating $f^{\tin}_{\mx}\left[\psi^{\tin}, \psi^{\out}\right]$ in the $\sharp$-norm, defined in \eqref{eq:sharp-norm}.

\begin{lemma}
We have
\begin{equation}\label{eq:iest1}
\left\|\eta_{2\mu_0^{\vep_1}}(x) W_{1,0}^{1-p} \({\bmu \over \mu}\)^{N-2 \over 2} (1+P)^{-p}(x) \mce_2[\mx](y,t)\right\|_{\sharp,a+2,b;\sigma_0} \le C_{71}
\end{equation}
where $C_{71} > 0$ is a constant depending only on $(M,g_0)$, $N$, $h$, $z_0$, $\alpha$, and $\sigma_0$.
\end{lemma}
\begin{proof}
It follows from \eqref{eq:E2} that the dominating term of $\mce_2[\mx]$ is a constant multiple of $\bmu\lambda W_{1,0}(y)$. For this term, we have
\begin{align*}
&\ \mu_0^{-b} \({\bmu \over \mu}\)^{N-2 \over 2} \bmu\lambda \left[\(1+|\by|^{a+2}\) \left|\eta_{2\mu_0^{\vep_1}}(x) (1+P)^{-p}(x) W_{1,0}(y)\right| \right. \\
&\left. \hspace{80pt} + \(1+|\by|^{a+2+\sigma_0}\) \left[\eta_{2\mu_0^{\vep_1}}(\bmu\cdot+\xi)
(1+P)^{-p}(\bmu\cdot+\xi) W_{1,0}\({\bmu \over \mu} \cdot\)\right]_{C^{\sigma_0}_{\R^N}}(\by,t) \right] \\
&\le C \frac{\mu_0^{\nu_1+1-b}}{1+|\by|^{N-4-a}} \le C
\end{align*}
and
\begin{multline*}
\mu_0^{-b} \(1+|\by|^{a+2-\sigma_0}\) \left[\bmu\lambda \({\bmu \over \mu}\)^{N-2 \over 2} \eta_{2\mu_0^{\vep_1}}(\bmu\cdot+\xi)
(1+P)^{-p}(\bmu\cdot+\xi) W_{1,0}\({\bmu \over \mu} \cdot\)\right]_{C^{\sigma_0/2}_t}(\by,t) \\
\le C \frac{\mu_0^{\nu_1+3-b}}{1+|\by|^{N-4-a+\sigma_0}} \le C\mu_0^2(t_0)
\end{multline*}
where we employed the condition $a \in (0,N-4]$; see Definition \ref{defn:semi-norm} for the definition of the local H\"older semi-norms.

Handling the other terms of $\mce_2[\mx]$ in an analogous way, we establish \eqref{eq:iest1}.
\end{proof}

\begin{lemma}
We have
\begin{equation}\label{eq:iest2}
\left\|\eta_{2\mu_0^{\vep_1}}(x) W_{1,0}^{1-p} \mck_1\left[\psi^{\tin}\right]\right\|_{\sharp,a+2,b;\sigma_0}
\le C_{72} \mu_0^{2\vep_1}(t_0) \left\|\psi^{\tin}\right\|_{\sharp',a,b;\sigma_0\(\mcb_{\bmx}\)}
\end{equation}
where $C_{72} > 0$ is a constant depending only on $(M,g_0)$, $N$, $h$, $z_0$, $\alpha$, and $\sigma_0$.
\end{lemma}
\begin{proof}
Arguing as in the proof of Lemma \ref{lemma:oest2}, we can estimate the first three terms of $\mck_1[\psi^{\tin}]$ in \eqref{eq:mck1}.

Let us consider the fourth term. From the inequalities
\begin{align*}
&\ \mu_0^{-b} \(1+|\by|^{a+2}\) \eta_{2\mu_0^{\vep_1}}(x) \left|\({\mu \over \bmu}\)^2 (1+P)^{1-p}(x) (\Delta_{g_0(x)} \psi^{\tin})(\by,t) - \Delta \psi^{\tin}(\by,t)\right| \\
&\le C \mu_0^{-b} \(1+|\by|^{a+2}\) \left[\(\mu_0^2|\by|^2 + \mu_0^{\nu_1-1}\) \left|\nabla^2 \psi^{\tin}(\by,t)\right| + \mu_0^2|\by| |\nabla \psi^{\tin}(\by,t)|\right] \\
&\le C \mu_0^{2\vep_1}(t_0) \left\|\psi^{\tin}\right\|_{\sharp',a,b\(\mcb_{\bmx}\)},
\end{align*}
we obtain its weighted $L^{\infty}$-bound; cf. \eqref{eq:ODExi4}. By further inspection, we deduce a weighted H\"older estimate.

Handling the remaining terms of $\mck_1[\psi^{\tin}]$ in an analogous way, we establish \eqref{eq:iest2}.
\end{proof}

\begin{lemma}
We have
\begin{equation}\label{eq:iest3}
\left\|\eta_{2\mu_0^{\vep_1}}(x) W_{1,0}^{1-p} \mck_2\left[\psi^{\tin}, \psi^{\out}\right]\right\|_{\sharp,a+2,b;\sigma_0}
\le C_{73} \(\mu_0^2(t_0) \left\|\psi^{\tin}\right\|_{\sharp',a,b;\sigma_0\(\mcb_{\bmx}\)} + \left\|\psi^{\out}\right\|_{*',\ar;\sigma_0}\)
\end{equation}
where $C_{73} > 0$ is a constant depending only on $(M,g_0)$, $N$, $h$, $z_0$, $\alpha$, and $\sigma_0$.
\end{lemma}
\begin{proof}
Let us consider the first term of $\mck_2[\psi^{\tin}]$. Applying the condition $a \in (0,\alpha+2]$, we deduce
\begin{align*}
&\ \mu_0^{-b} \(1+|\by|^{a+2}\) \eta_{2\mu_0^{\vep_1}}(x) \left|\({\bmu \over \mu}\)^{N-2 \over 2} (1+P)^{-1}(x) W_{1,0}^{p-1}(y) \mu^{N-2 \over 2} \psi_{\mx}^{\out}(x,t)\right| \\
&\le C \left\|\psi^{\out}\right\|_{*',\ar;\sigma_0} \sup_{(\by,t) \in \mcb_{\bmx} \times [t_0,\infty)} {\mu_0^{\beta-b} \over 1+|\by|^{\alpha-a+2}} \le C \left\|\psi^{\out}\right\|_{*',\ar;\sigma_0},
\end{align*}
from which we obtain a weighted $L^{\infty}$-bound. By further inspection, we deduce a weighted H\"older estimate.

Employing the mean value theorem, we can handle the second term of $\mck_1[\psi^{\tin}]$.
Inequality \eqref{eq:iest3} then readily follows.
\end{proof}

\begin{proof}[Completion of the proof of Proposition \ref{prop:inner2}]
Let $\mct^{\tin}[\mch] := \psi[\mch]$ and $e_0[\mch]$ be the solution to \eqref{eq:inhom-i} found in Proposition \ref{prop:lin}.
A function $\psi^{\tin}$ solves \eqref{eq:inner2} if it satisfies
\[\psi^{\tin} = \mct^{\tin}\left[f^{\tin}_{\mx}\left[\psi^{\tin}, \psi^{\out}\left[\psi^{\tin}; \psi_0\right]\right]\right]\]
where $f^{\tin}_{\mx}$ is the map defined in \eqref{eq:fin},
and the notation $\psi^{\out}[\psi^{\tin}; \psi_0]$ emphasizes the dependence of $\psi^{\out}$ on $\psi^{\tin}$ and $\psi_0$.

We claim that the operator $\mct^{\tin} \circ f^{\tin}_{\mx}$ has a fixed point on the set
\[\mcd^{\tin} := \left\{\psi^{\tin}: \left\|\psi^{\tin}\right\|_{\sharp',a,b;\sigma_0\(\R^N\)} \le 2C_4C_{71}\right\}\]
where $C_4$ and $C_{71}$ are the numbers appearing in \eqref{eq:lin0} and \eqref{eq:iest1}.

By using \eqref{eq:lin0}, \eqref{eq:fin}, \eqref{eq:iest1}, \eqref{eq:iest2}, \eqref{eq:iest3}, and \eqref{eq:outer32}, and taking a large $t_0 > 0$ if needed, we get
\begin{equation}\label{eq:iest4}
\begin{aligned}
&\ \left\|\mct^{\tin}\left[f^{\tin}_{\mx}\left[\psi^{\tin}, \psi^{\out}\left[\psi^{\tin};
\psi_0\right]\right]\right]\right\|_{\sharp',a,b;\sigma_0\(\R^N\)} \\
&\le C_4 \left\|f^{\tin}_{\mx}\left[\psi^{\tin}, \psi^{\out}\left[\psi^{\tin}; \psi_0\right]\right]\right\|_{\sharp,a+2,b;\sigma_0} \\
&\le C_4 \left[C_{71} + \(C_{72} \mu_0^{2\vep_1}(t_0) + C_{73} \mu_0^2(t_0)\) \left\|\psi^{\tin}\right\|_{\sharp',a,b;\sigma_0\(\mcb_{\bmx}\)}
+ C_{73} \left\|\psi^{\out}\right\|_{*',\ar;\sigma_0} \right] \le 2C_4C_{71}.
\end{aligned}
\end{equation}
Moreover, by the linearity of $\mct^{\tin}$, \eqref{eq:lin0}, \eqref{eq:outer35}, and \eqref{eq:inner1a} (see also the derivation of \eqref{eq:iest2}),
\begin{align*}
&\ \left\|\mct^{\tin}\left[f^{\tin}_{(\mx)[\psi^{\tin}_1]} \left[\psi^{\tin}_1, \psi^{\out}\left[\psi^{\tin}_1; \psi_0\right]\right]\right]
- \mct^{\tin}\left[f^{\tin}_{(\mx)[\psi^{\tin}_2]} \left[\psi^{\tin}_2, \psi^{\out}\left[\psi^{\tin}_2; \psi_0\right]\right]\right]\right\|_{\sharp',a,b;\sigma_0\(\R^N\)} \\
&\le C \left\|f^{\tin}_{(\mx)[\psi^{\tin}_1]} \left[\psi^{\tin}_1, \psi^{\out}\left[\psi^{\tin}_1; \psi_0\right]\right]
- f^{\tin}_{(\mx)[\psi^{\tin}_2]} \left[\psi^{\tin}_2, \psi^{\out}\left[\psi^{\tin}_2; \psi_0\right]\right]\right\|_{\sharp,a+2,b;\sigma_0} \\
&\le C \delta_0^{-\zeta} \mu_0^{\min\{\delta_4, 2\vep_1\}}(t_0) \left\|\psi^{\, \tin}_1 - \psi^{\, \tin}_2\right\|_{\sharp',a,b;\sigma_0\(\mcb_{\bmx}\)}
\le \frac{1}{2} \left\|\psi^{\, \tin}_1 - \psi^{\, \tin}_2\right\|_{\sharp',a,b;\sigma_0\(\mcb_{\bmx}\)}
\end{align*}
where the notation $(\mx) = (\mx)[\psi^{\tin}]$ stresses the dependence of $(\mx)$ in $\psi^{\tin}$.
By the contraction mapping theorem, the assertion holds.

\medskip
We now have a unique solution $\psi^{\tin}$ to \eqref{eq:inner2} with the desired $\sharp'$-norm bound given in \eqref{eq:inner3}.
Also, the bound on $e_0(t_0)$ in \eqref{eq:inner3} immediately follows from \eqref{eq:lin0} and \eqref{eq:iest4}.
A further inspection based on \eqref{eq:outer36} and \eqref{eq:inner1b} gives \eqref{eq:inner31}, concluding the proof.
\end{proof}

\subsection{Completion of the proof of Theorem \ref{thm:main}} \label{subsec:inner4}
Let $z_0$ be a point on $M$ such that $h(z_0) > 0$.
In Subsection \ref{subsec:approx1}, we selected sufficiently small numbers $\vep_0,\, \vep_1,\, \sigma_0 \in (0,1)$ and set $\nu_1 = \nu_2 = 2-\vep_0$. Take
\[a = N-4, \quad b = 3-\vep_0, \quad \alpha = N-5+\vep_0,
\quad \beta = 3-\vep_0, \quad \delta_2 = \vep_0,\]
and any $\delta_4 \in (0,\vep_0]$ satisfying \eqref{eq:delta4}. Choose also $\psi_0 = 0$ on $M$.
From the discussion in Section \ref{sec:inout} and Propositions \ref{prop:outer3}, \ref{prop:inner1}, and \ref{prop:inner2},
we find a solution $u_{z_0}$ to \eqref{eq:Yamabefp} of the form $u_{z_0} = u_{\mx} + \psi_{\mx}$ on $M \times [t_0,\infty)$
where $u_{\mx} = u_{\mx}^{(2)}$ and $\psi_{\mx}$ are given in \eqref{eq:u2} and \eqref{eq:psimxo}, respectively.
Estimates \eqref{eq:outer32}, \eqref{eq:lx}, and \eqref{eq:inner3} (or \eqref{eq:outer31}) imply that $u_{z_0} > 0$ on $M \times [t_0,\infty)$ and \eqref{eq:approx} holds.
Consequently, the proof of Theorem \ref{thm:main} is completed.

\medskip
Note that we have a freedom to choose the initial value $\psi_0 = \psi^{\out}(\cdot,t_0)$ on the outer problem \eqref{eq:outer} or \eqref{eq:outer3},
provided the second inequality of \eqref{eq:outer31} holds.
In Subsection \ref{subsec:cor}, we will further analyze this observation to establish the $k$-codimensional stability stated in Corollary \ref{cor:main}.

\section{Proof of Theorem \ref{thm:main2} and Corollary \ref{cor:main}} \label{sec:multi}
Throughout the section, $k \in \N$ is fixed and $l$ can take any integer between $1$ and $k$.

\subsection{Proof of Theorem \ref{thm:main2}}
In this subsection, we provide the outline of the proof of Theorem \ref{thm:main2}, pointing out the changes needed with respect to the one bubble case.

\medskip
By choosing the number $\delta_0 > 0$ sufficiently small, we may assume that $\mbz_0 = (z_0^{(1)}, \ldots, z_0^{(k)})$ is an element of the configuration set
\[\mcc := \left\{\(z_0^{(1)}, \ldots, z_0^{(k)}\) \in M^k: d_{g_0}\(z_0^{(l)},z_0^{(m)}\) \ge c\delta_0 \text{ for } 1 \le l \ne m \le k \right\}\]
where $c > 0$ is a number determined by $(M,g_0)$, $N$, $h$, and $k$.
Given $l = 1, \ldots, k$, we write
\begin{equation}\label{eq:mul}
d^{(l)} := \frac{1}{\sqrt{h\(z_0^{(l)}\)}} \quad \text{and} \quad
\mu^{(l)}(t) = d^{(l)}\mu_0(t) + \lambda^{(l)}(t) =: \bmu^{(l)}(t) + \lambda^{(l)}(t) \quad \text{for } t \in [t_0,\infty)
\end{equation}
where $\mu_0$ is the function in \eqref{eq:mu0}, and $\lambda^{(l)}$ is a higher-order term.
We assume that \eqref{eq:lx} holds for $(\lambda,\xi) = (\lambda^{(l)}, \xi^{(l)})$.

By substituting $z_0^{(l)}$ and $\bmu^{(l)}$ for $z_0$ and $\bmu$, respectively,
we define the analogue $P^{(l)}$ of $P$ in \eqref{eq:P_21}, and the analogue $\Psi_0^{(l)}$ of $\Psi_0$ in the paragraph after Lemma \ref{lemma:Qdecay}.
In \eqref{eq:u2}--\eqref{eq:v2}, we put $z_0$, $(\mx)$ and $\Psi_0$ in place of $z_0^{(l)}$, $(\mu^{(l)},\xi^{(l)})$, respectively,
to define $u_{\mu^{(l)},\xi^{(l)}}^{(2)}$ and $v_{\mu^{(l)},\xi^{(l)}}^{(2)}$.
Then we set the refined approximate solution
\begin{equation}\label{eq:umbmx}
u_{\mbmx}(z,t) = u_{\mbmx}^{(2)}(z,t) = \sum_{l=1}^k u_{\mu^{(l)},\xi^{(l)}}^{(2)}(z,t) \quad \text{on } M \times [t_0,\infty)
\end{equation}
where $\mbmu := (\mu^{(1)}, \ldots, \mu^{(k)}) \in (0,\infty)^k$ and $\mbxi := (\xi^{(1)}, \ldots, \xi^{(k)}) \in (\R^N)^k$.
It holds that
\[\mu^{N+2 \over 2} \mcs\(u^{(2)}_{\mbmx}\)(y,t) = \mce_2\left[\mu^{(l)},\xi^{(l)}\right](y,t)\]
for $z = \exp_{z_0^{(l)}}(x) = \exp_{z_0^{(l)}}(\mu^{(l)}y+\xi^{(l)}) \in B_{g_0}(z_0^{(l)},\delta_0)$
and $t \in [t_0,\infty)$, where $\mce_2$ is the function in \eqref{eq:E2}.
Indeed, the interaction between two different bubbles is estimated by
\[\mbone_{N=5} \mu_0^3 a^{\{1\}} + \mu_0^{N-2} a^{\{N-4\}} \(\lesssim \mbone_{N=5} \mu_0^3 a^{\{1\}} + \mu_0^4 a^{\{2\}}\),\]
which is smaller than the main order terms of $\mce_2$.
An analogous formula to \eqref{eq:S2out} is also true.

\medskip
The norms in Subsection \ref{subsec:norms} must be adjusted accordingly.
For instance, every $u_{\mx}$ in the norms has to be replaced with $u_{\mbmx}$,
and the weight $w_{\alpha,\gamma}$ in \eqref{eq:wag} needs to be redefined as
\begin{multline*}
w_{\alpha,\ga}(z,t) = \max \left\{\sum_{l=1}^k \eta_{\delta_0}\(\left|x^{(l)}\right|\)
\frac{\mu_0^{-\ga}}{1 + \left|\(\mu^{(l)}\)^{-1} \(x^{(l)}-\xi^{(l)}\)\right|^{\alpha+\ga}}, \right. \\
\left. 2^{3(\alpha+\ga)} \min\left\{h\(z_0^{(1)}\),\ldots,h\(z_0^{(k)}\)\right\}^{-\frac{\alpha+\ga}{2}} \delta_0^{-(\alpha+\ga)} \mu_0^{\alpha}\right\}
\end{multline*}
where $x^{(l)} := \exp_{z_0^{(l)}}^{-1}(z)$.

\medskip
Let $\psi_{\mbmx}$ stand for the remainder term such that $u = u_{\mbmx} + \psi_{\mbmx}$ is a solution to \eqref{eq:Yamabefp2}.
We assume that it can be written as
\begin{equation}\label{eq:psimbmx}
\psi_{\mbmx} = \psi_{\mbmx}^{\out} + \psi_{\mbmx}^{\tin} \quad \text{on } M \times [t_0,\infty)
\end{equation}
where $\psi_{\mbmx}^{\tin}$ has the form
\begin{equation}\label{eq:psitin}
\psi_{\mbmx}^{\tin}(z,t) = \sum_{l=1}^k \eta_{\mu_0^{\vep_1}} \(\left|x^{(l)}\right|\)\,
\(1+P^{(l)}\(x^{(l)}\)\) \(\bmu^{(l)}\)^{-{N-2 \over 2}} \(\psi_{\mbmx}^{\tin}\)^{(l)}\(\by^{(l)}, t\)
\end{equation}
for $(z,t) \in M \times [t_0,\infty)$ and $x^{(l)} = \bmu^{(l)}\by^{(l)}+\xi^{(l)}$.
The inner-outer gluing procedure consists of finding a solution $\psi_{\mbmx}^{\out}$ to the outer problem on $M \times [t_0,\infty)$
and a solution $(\psi_{\mbmx}^{\, \tin})^{(l)}$ of an inner problem on $B^N(-(\bmu^{(l)})^{-1}\xi^{(l)}, 2(\bmu^{(l)})^{-1}\mu_0^{\vep_1}) \times [t_0,\infty)$ for $l = 1, \ldots, k$.

\medskip
Minor modifications of the arguments in Sections \ref{sec:outer} and \ref{sec:inner}
yield a priori estimates for inhomogeneous equations associated to the outer and inner problems.
For example, the only change required in the proof of Lemma \ref{lemma:outer1a} is to replace \eqref{eq:outer1a8} with
\[\int_0^1 \int_{B_{g_0}\(z_0^{(m)},R_1 \mu_{\ell}^{(m)} (\cdot+\tau_{\ell})\)} \phi_\ell^2 u_\ell^{p-1} dv_{g_0} dt \ge \frac{1}{4k}\]
for some $m \in \{1, \ldots, k\}$. The existence of such $m$ is guaranteed by the pigeonhole principle.
To prove Lemma \ref{lemma:outer2}, we divide the manifold $M$
into $k+1$ subsets $B_{g_0}(z_0^{(1)},\frac{\delta_0}{4}), \ldots, B_{g_0}(z_0^{(k)},\frac{\delta_0}{4})$, and $M \setminus \bigcup_{l=1}^k B_{g_0}(z_0^{(l)},\frac{\delta_0}{6})$,
and then examine the behavior of $\psi$ on each subset.

By employing the a priori estimates for inhomogeneous equations
and adopting the arguments in Sections \ref{sec:outer} and \ref{sec:inner} once again,
we establish the existence of $\psi_{\mbmx}^{\out}$, $(\lambda^{(l)}, \xi^{(l)})$,
and $(\psi_{\mbmx}^{\, \tin})^{(l)}$ for $l = 1, \ldots, k$ as well as estimates on their respective norms.
Combining this information leads us to complete the proof of Theorem \ref{thm:main2}.

\subsection{Proof of Corollary \ref{cor:main}}\label{subsec:cor}
By modifying the argument in the proof of \cite[Corollary 1.1]{CdPM} suitably, one can prove the corollary. Here we give a sketch the proof.

As in the previous subsection, a quantity with the superscript $(l)$ indicates that it is related to the $l$-th blow-up point $z_0^{(l)}$.

\medskip
Fix $l = 1, \ldots, k$ and let $\mcm$ be the matrix in \eqref{eq:mcm} with $z_0 = z_0^{(l)}$.
If $\vsi_1, \ldots, \vsi_N$ are the eigenvalues of $\mcm$, then the assumption on the Ricci curvature implies that
\begin{equation}\label{eq:vsi}
\min\{\vsi_1, \ldots, \vsi_N\} \ge 1 > 1-\frac{\vep_0}{2} = \frac{\nu_2+2}{2}-1.
\end{equation}
Set $(\overline{f^{\tpa}})^{(l)} = ((f^{\tpa}_1)^{(l)}, \ldots, (f^{\tpa}_N)^{(l)})$
and an orthogonal matrix $\mcq$ such that $\mcm = \mcq^T \mcd \mcq$ with $\mcd = \text{diag}(\vsi_1, \ldots, \vsi_N)$.
In light of \eqref{eq:vsi}, \eqref{eq:xieq2}, \eqref{eq:laeq}, \eqref{eq:mct}, and \eqref{eq:inner11},
we can select the parameters $(\lambda^{(l)}, \xi^{(l)})$ by
\[\lambda^{(l)}(t) = t^{-{3 \over 2}} \int_{t_0}^t s^{3 \over 2} \(f^{\tpa}_{N+1}\)^{(l)}(s)ds\]
and
\begin{equation}\label{eq:xil}
\xi^{(l)}(t) = \mcq^T \tilde{\xi}^{(l)}(t) \quad \text{where} \quad
\tilde{\xi}^{(l)}_i(t) = t^{-\vsi_i} \int_{t_0}^t s^{\vsi_i} \left[\mcq \(\overline{f^{\tpa}}\)^{(l)}\right]_i(s) ds
\quad \text{for } i = 1, \ldots, N.
\end{equation}
They satisfy \eqref{eq:lx}, and
\begin{equation}\label{eq:xil2}
\mu^{(l)}(t_0) = d^{(l)} \mu_0(t_0), \quad \xi^{(l)}(t_0) = 0 \quad \text{and} \quad \xi^{(l)}(t) \to 0 \quad \text{as } t \to \infty
\end{equation}
by \eqref{eq:mul}.\footnote{Suppose that $\vsi_1 < \frac{\nu_2+2}{2}-1 = 1-\frac{\vep_0}{2}$.
If we define $\xi^{(l)}$ by \eqref{eq:xil}, then we cannot obtain the estimate $\|\dot{\xi}\|_{\nu_2+2;\sigma} \le C \|\ovmfh\|_{\nu_2+2;\sigma}$ needed in the proof of Proposition \ref{prop:inner1}.}

Let $\psi_0$ be the initial condition for the outer problem,
and $e_0^{(l)}[\psi_0]Z_0$ the initial condition for the inner problem of $(\psi_{\mbmx}^{\, \tin})^{(l)}$.
Recalling \eqref{eq:umbmx} and \eqref{eq:psimbmx}--\eqref{eq:psitin}, we choose the initial datum for \eqref{eq:Yamabefp2} of the form
\begin{align*}
u(z,t_0) &= u_{\mbmx}(z,t_0) + \sum_{l=1}^k \eta_{\mu_0^{\vep_1}(t_0)} \(\left|x^{(l)}\right|\)
\(1+P^{(l)}\(x^{(l)}\)\) \(\bmu^{(l)}(t_0)\)^{-{N-2 \over 2}} e_0^{(l)} [0](t_0)Z_0(z) \\
&\ + F[\psi_0](z) \quad \text{for } z \in M
\end{align*}
where
\begin{multline*}
F[\psi_0](z) := \psi_0(z) \\
+ \sum_{l=1}^k \eta_{\mu_0^{\vep_1}(t_0)} \(\left|x^{(l)}\right|\)
\(1+P^{(l)}\(x^{(l)}\)\) \(\bmu^{(l)}(t_0)\)^{-{N-2 \over 2}} \left[e_0^{(l)} [\psi_0](t_0) - e_0^{(l)} [0](t_0)\right]Z_0(z).
\end{multline*}
According to \eqref{eq:xil2}, the solution $u(z,t)$ to \eqref{eq:Yamabefp2} blows up precisely at $z_0^{(1)}, \ldots, z_0^{(k)}$ on $M$.
Besides, $F$ is a $C^1$-function on $C^{2,\sigma_0}(M)$ (see \eqref{eq:inner31}), $F[0]=0$,
and $DF[\Psi_0] = \text{Id}$ on the subspace $W := \cap_{l=1}^k \text{ker} (D_{\psi_0}e_0^{(l)}[0])$ of $C^{2,\sigma_0}(M)$.
Hence the inverse function theorem says that there is a manifold of codimension $\text{codim}(W) \le k$ in a neighborhood of $0 \in C^{2,\sigma_0}(M)$
such that each element is expressed as $F[\psi_0]$ for some $\psi_0 \in C^{2,\sigma_0}(M)$ near $0$.
Calling such a manifold $\mcm_{\mbz_0}$ and taking $\sigma = \sigma_0$, we conclude the proof.

\bigskip \noindent \small{\textbf{Acknowledgement.}
S. Kim was supported by Basic Science Research Program through the National Research Foundation of Korea (NRF) funded by the Ministry of Education
(NRF2020R1C1C1A01010133, NRF2020R1A4A3079066) and associate member program of Korea Institute for Advanced Study (KIAS).
M. Musso has been supported by EPSRC research Grant EP/T008458/1.}

\end{document}